\documentclass[11pt]{amsart}%
\usepackage{amsfonts}
\usepackage{amsmath}
\usepackage{amssymb}
\usepackage{graphicx}
\usepackage{ulem}
\usepackage[colorlinks=true, linkcolor=blue, citecolor=red, urlcolor=red,
backref=page]{hyperref}%
\providecommand{\U}[1]{\protect\rule{.1in}{.1in}}
\newtheorem{theorem}{Theorem}[section]
\theoremstyle{plain}

\newtheorem{lemma}{Lemma}

\newtheorem{proposition}{Proposition}

\theoremstyle{definition}
\newtheorem{definition}{Definition}
\newtheorem{example}{Example}
\numberwithin{equation}{section}
\theoremstyle{remark}
\newtheorem{remark}[theorem]{Remark}
\def\R{\mathbb{R}}
\def\T{\mathbb{T}}
\def\Z{\mathbb{Z}}
\def\N{\mathbb{N}}

\begin{document}
\title{Uniqueness for $2$D Euler and transport equations via extrapolation}
\author{\'{O}scar Dom\'{\i}nguez}
\address{Department of Mathematics, 
CUNEF University, 
Madrid, Spain}
\email{oscar.dominguez@cunef.edu}
\author{Mario Milman}
\address{Instituto Argentino de Matematica, Buenos Aires, Argentina}
\email{mario.milman@icloud.com}
\urladdr{https://sites.google.com/site/mariomilman}

\thanks{\emph{Acknowledgements.} \'Oscar Dom\'{i}nguez is supported by the AEI grant RYC2022-037402-I. Part of this research was carried out while he was a postdoc at the Institut Camille Jordan, Lyon, supported by the Labex Milyon and the French National Research Agency (ANR-10-LABX-0070), (ANR-11-IDEX-0007). It is our pleasure to thank Prof. Petru Mironescu for some helpful comments that helped to improve the presentation of the paper.}

\begin{abstract}
Using extrapolation theory, we develop a new framework to prove the uniqueness
of weak solutions for a wide class of active scalar equations, including $2$D Euler and Navier-Stokes equations, SQG, IPM, Keller-Segel models, aggregation equations, ... We  apply our
methodology to unify and extend the classical results of Yudovich and Vishik
for $2$D Euler equations. In particular, we establish the uniqueness for the Euler flow whose vorticity belongs to  new scales of function spaces  that contain both Yudovich spaces and $\text{BMO}$.  We give a self contained presentation.

\end{abstract}
\maketitle
\tableofcontents

\section{Introduction}

\subsection{Uniqueness of Euler flows}

\label{Section1.1} The uniqueness of weak solutions of the Cauchy problem for
the Euler equations of an ideal incompressible fluid on $\Omega=\mathbb{R}%
^{2}$ or $\mathbb{T}^{2}$ still presents challenging open questions. In $2$D
the Euler equations can be formulated in terms of a transport equation
\begin{equation}
\left\{
\begin{array}
[c]{c}%
\omega_{t}+v\cdot\nabla\omega=0,\qquad\text{on }(0,\infty)\times\Omega,\\
v=k\ast\omega,\qquad\text{on }(0,\infty)\times\Omega,\\
\omega(0,x)=\omega_{0}(x),\qquad x\in\Omega,
\end{array}
\right.  \label{intro1}%
\end{equation}
where $v:[0,\infty)\times\Omega\rightarrow\mathbb{R}^{2}$ is the
\emph{velocity} of the flow, $\omega=\operatorname{curl}(v)=\partial_{x_1}%
v_{2}-\partial_{x_2}v_{1}$ is the \emph{vorticity}, and $k$ is standard
Biot-Savart kernel\footnote{The expression of $k$ is particularly simple in
$\mathbb{R}^{2},$ namely, $k(x)=\frac{x^{\perp}}{2\pi\left\vert x\right\vert
^{2}}$ for $x=(x_{1},x_{2})\neq\mathbf{0}$. Here $x^{\perp}=(-x_{2},x_{1})$.}.
The equation $v=k\ast\omega$ is the so-called \emph{Biot-Savart law}, allowing
to recover $v$ from $\omega$. It follows (cf. \cite{Y63, L96, C98}) that%
\begin{equation}
\nabla v=\mathcal{K}\omega, \label{intro2}%
\end{equation}
where $\mathcal{K}$ is a certain Calder\'{o}n--Zygmund operator (in short CZO).

The uniqueness of weak solutions to \eqref{intro1} with bounded vorticity was
established in a celebrated paper by Yudovich \cite{Y63}, via a clever
argument involving the energy method combined with sharp $L^{p}$-inequalities
for CZO. To extend the uniqueness result to
unbounded vorticities, Yudovich \cite{Y95} proposed new function spaces,
nowadays called \emph{Yudovich spaces}, that are somewhat bigger than
$L^{\infty}$ (or $L^{p_{0}}\cap L^{\infty}$ in the case of unbounded domains
$\Omega$). More precisely, these spaces are constructed by means of fixing a
non-decreasing doubling function $\Theta:(1,\infty)\rightarrow(0,\infty)$
(\textquotedblleft a growth function\textquotedblright) and $p_{0}\geq1,$ and
letting $Y_{p_{0}}^{\Theta}(\Omega)$  be the set of all $\omega\in
\cap_{p > p_{0}}L^{p}(\Omega)$ such that
\begin{equation}
\left\Vert \omega\right\Vert _{Y_{p_{0}}^{\Theta}(\Omega)}:=\sup_{p > p_{0}%
}\frac{\left\Vert \omega\right\Vert _{L^{p}(\Omega)}}{\Theta(p)}<\infty.
\label{intro:YudoSpace}%
\end{equation}
In particular, if\footnote{Given two non-negative quantities $A$ and $B$, the notation $A \lesssim B$ means that there exists a constant $C$, independent of all essential parameters, such that $A \leq C B$. We write $A \approx B$ if $A \lesssim B \lesssim A$.} $\Theta(p)\approx1,$ it is easy to see that (cf. Appendix
\ref{SectionA3},  Proposition \ref{ejemplomarkao})%
\begin{equation}
Y_{p_{0}}^{\Theta}(\Omega)=L^{p_{0}}(\Omega)\cap L^{\infty}(\Omega),
\label{introYLinfty}%
\end{equation}
while if $\Theta(p)\rightarrow\infty,$ as $p\rightarrow\infty,$ $Y_{p_{0}%
}^{\Theta}(\Omega)$ contains unbounded functions. Although the definitions do not give
explicit descriptions of the $Y_{p_{0}}^{\Theta}$ spaces, it is relatively
easy to construct examples of elements belonging to them. It is worthwhile to
mention that the definition of $Y_{p_{0}}^{\Theta}(\Omega)$ is independent of
$p_{0}$ if $\Omega = \T^2$ (or more generally, if $\Omega$ has finite measure), but if $\Omega=\mathbb{R}^{2}$ (cf.
\eqref{introYLinfty}) we only have the trivial embeddings\footnote{Given two
normed spaces $X$ and $Y$, we write $X\hookrightarrow Y$ if $X\subset Y$ and
the natural embedding from $X$ to $Y$ is continuous.}
\[
Y_{p_{0}}^{\Theta}(\Omega)\hookrightarrow Y_{p_{1}}^{\Theta}(\Omega
),\qquad\text{if}\qquad p_{1}>p_{0}.
\]

The uniqueness results for $Y^{\Theta}_{p_{0}}$ vorticities are formulated in
terms of functions associated with a growth function $\Theta$, which we shall
term \emph{Yudovich functions}\footnote{In \cite{Y95}, the function
$y_{\Theta}$ is defined in slightly different way for $r \in(0, 1)$. However,
this minor modification will not play a role in what follows, cf.
\eqref{intro:osg}.} $y_{\Theta},$%
\begin{equation}
y_{\Theta}(r):=\inf_{p>p_{0}}\{\Theta(p)r^{1/p}\},\qquad r>0.
\label{intro:Auxy}%
\end{equation}
In the theory of \cite{Y95}, Yudovich functions appear naturally in the
crucial differential inequality connecting $Y^{\Theta}_{p_{0}}$ solutions with
the energy method. Indeed, let $(v_{i},\omega_{i}),i=1,2,$ be two solutions of
(\ref{intro1})  and let
\begin{equation}\label{EnergyDef}
E(t)= \frac{1}{2} \left\Vert v_1(t, \cdot) - v_2 (t, \cdot)\right\Vert _{L^2(\Omega)}^{2}.
\end{equation} From the definition of
weak solution, it is plain to see that
\begin{equation}\label{EnergyEstim}
	\frac{d E(t)}{dt}  \leq \int_\Omega |\nabla v_1(t, x)| |v_1(t, x)-v_2(t, x)|^2 \, dx. 
\end{equation}
Applying  the norm of $Y^{\Theta}_{p_{0}} (\Omega)$, combined with the
Biot-Savart law (\ref{intro2}) together with the sharp $L^{p}$-norm estimates
for CZO, one can get
\begin{equation}\label{KeyYudo}
\frac{dE(t)}{dt} \lesssim \left\Vert \omega_1(t)\right\Vert _{Y^{\Theta}_{p_{0}} (\Omega)}E(t)
\, y_{\Theta_{1}}(E(t)^{-1}), 
\end{equation}
where $y_{\Theta_{1}}$ is the Yudovich function \eqref{intro:Auxy} associated with the growth
function $\Theta_{1}$ given by
\begin{equation}
\Theta_{1}(p):=p\,\Theta(p). \label{intro:AuxTheta1}%
\end{equation}
Uniqueness (i.e., $v_1\equiv v_2$) is then achieved under the following Osgood condition
on $y_{\Theta_{1}}$(cf. \cite{Y95})
\begin{equation}
\int_{0}^{1}\frac{dr}{ry_{\Theta_{1}}(\frac{1}{r})}=\infty. \label{intro:osg}%
\end{equation}
Moreover, Yudovich \cite[Theorem 2]{Y95} also proved that the velocity $v$ satisfies 
\begin{equation}\label{18a}
	|v(t, x)-v(t, y)| \lesssim \|\omega(t)\|_{Y^{\Theta}_{p_0}(\Omega)} |x-y| y_{\Theta_1} \bigg(\frac{1}{|x-y|^3} \bigg)
\end{equation} 
for every $x, y \in \Omega$ and a.e. $t > 0$.

It is easy to verify that (\ref{intro:osg}) holds for $\Theta(p) \approx1$,
then from \eqref{introYLinfty} we obtain the classical uniqueness result for bounded vorticity 
of \cite{Y63}. It also holds for $\Theta(p)\approx\log p,$ whose corresponding
space $Y^{\Theta}_{p_{0}}$ can be seen to include unbounded vorticities of the
form 
\begin{equation}\label{PY}
\omega(x) \approx\left\vert \log\left\vert \log\left\vert x\right\vert
\right\vert \right\vert, \qquad x \to 0.
\end{equation}
 However, (\ref{intro:osg}) places a severe
restriction on $\Theta$ and, indeed, it fails for linear growth $\Theta
(p)\approx p,$ which corresponds to vorticities $\omega$ in the Orlicz space
$e^{L}$ of exponentially integrable functions. It follows that, by this method, uniqueness cannot be guaranteed for
vorticities of the form $\omega(x) \approx\left\vert \log\left\vert
x\right\vert \right\vert ,$ the prototype of an unbounded function in BMO (cf.
\cite{JN61}). We mention that an alternative and elementary approach to the well-posedness of $2$D Euler equations in $Y^\Theta_{p_0}$ (more precisely, its localized version) has been recently proposed by Crippa and Stefani \cite{CS21}.

Working in the full plane, the expected uniqueness result for vorticities with
logarithmic singularities was obtained somewhat later by Vishik \cite{V99},
using a different method. Interestingly, Vishik's method is also constructive
and relies on the introduction of the \textquotedblleft Vishik
spaces\textquotedblright\ $B_{\Pi}(\mathbb{R}^{d})$ associated with
\textquotedblleft growth functions\textquotedblright\ $\Pi$, that control the
growth of partial sums of the $L^{\infty}$-norm of dyadic frequency
localizations $\{\Delta_{j}f\}_{j\geq0}$ of  $f\in\mathcal{S}^{\prime
}(\mathbb{R}^{d})$. Specifically,
\begin{equation}
\Vert f\Vert_{B_{\Pi}(\mathbb{R}^{d})}:=\sup_{N\geq0}\frac{1}{\Pi(N)}%
\,\sum\limits_{j=0}^{N}\left\Vert \Delta_{j}f\right\Vert _{L^{\infty
}(\mathbb{R}^{d})}<\infty. \label{introVishikSpace}%
\end{equation}
Note that $\Pi(N)\approx1$ gives $B_{\Pi}(\mathbb{R}^{d})=B_{\infty,1}%
^{0}(\mathbb{R}^{d})$, a classical Besov space.

In this approach the Osgood condition on $\Pi,$ controlling
uniqueness, is given by%
\begin{equation}
\int_{1}^{\infty}\frac{dr}{\Pi(r)}=\infty. \label{IntroPi}%
\end{equation}
Under this assumption, uniqueness of the Euler flow of is guaranteed provided
that\footnote{The result holds for $d$ arbitrary, under the assumption $p_0 \in (1, d)$.}
\begin{equation}
\omega\in L^{\infty}([0,T];B_{\Pi}(\mathbb{R}^{2})\cap L^{p_{0}}%
(\mathbb{R}^{2})) \label{IntroViSpa}%
\end{equation}
for some $p_{0}\in(1,2)$. Moreover, in this case the following estimate for the modulus of continuity of $v$ holds (compare with \eqref{18a})
\begin{equation}\label{18c}
	|v(t, x)-v(t, y)| \lesssim |x-y| \big(\|\omega\|_{L^\infty ([0, T]; L^{p_0})} + \|\omega\|_{L^\infty([0, T]; B_\Pi)} \Pi(|\log |x-y||) \big). 
\end{equation}

In particular, $\Pi(p)\approx p$ satisfies \eqref{IntroPi} and uniqueness for
vorticities satisfying\footnote{$\text{bmo}(\mathbb{R}^{d})$ refers to the
local (a.k.a. inhomogeneous) version of $\text{BMO}(\mathbb{R}^{d})$. Recall
that $\text{bmo}(\mathbb{R}^{d}) \subsetneq\text{BMO}(\mathbb{R}^{d})$.}
\begin{equation}
\label{VUC}\omega\in L^{\infty}([0, T]; \text{bmo} (\mathbb{R}^{2})\cap
L^{p_{0}}(\mathbb{R}^{2})), \qquad p_{0} \in(1, 2),
\end{equation}
can be deduced from the embeddings
\begin{equation}
\text{bmo}(\mathbb{R}^{2})\hookrightarrow B_{\infty,\infty}^{0}(\mathbb{R}%
^{2})\hookrightarrow B_{\Pi}(\mathbb{R}^{2}). \label{intro112}%
\end{equation}

Another interesting approach to uniqueness of Euler equations in $\text{BMO}$ is due to Azzam and Bedrossian
\cite{AB}, cf. also \cite{Co} for an extension to $B^0_{\infty, \infty}$. The method of \cite{AB} (built on $\dot{H}^{-1}$ norms) may be considered as a
further refinement of the energy method of \cite{Y63}, that takes into account
not only the integrability, but also the inherited regularity properties of
$\text{BMO}$ functions.

In what concerns larger classes of vorticities, a challenging open problem in
the theory is to decide whether uniqueness of solutions for \eqref{intro1} can
be achieved for vorticities in $L^{p}.$ In recent remarkable work \cite{V2018}, Vishik (cf. also the lecture notes \cite{ADele}) shows that, for any
$2<p<\infty,$ there exist $\omega_{0}\in L^{1}(\R^2)\cap L^{p}(\R^2),$ and a force $f$,
such that there exist infinitely many weak solutions of
\[
\omega_{t}+v\cdot\nabla\omega=f,
\]
with $\omega\in L^{1} (\R^2)\cap L^{p}(\R^2)$, uniformly in time.  Without force (i.e., $f \equiv 0$), a very recent work of Bru\'e, Colombo and Kumar \cite{BCK} asserts non-uniqueness with some integrable (periodic) initial vorticity $\omega_0$. For general $L^p$ without force, the  (non-)uniqueness problem is widely open.


In view of the above discussion, an outstanding question in the area is to obtain sufficiently large classes of vorticities between  $L^p$  and $\text{BMO}$ that still guarantee   uniqueness for the $2$D Euler flow. In this paper, we address this question introducing the \emph{sharp Yudovich space} $Y^{\#\Theta}_{p_0}(\Omega)$ (the precise definition will be postponed to Section \ref{SeIntroBMOY},  Definition \ref{DefInYS}). This space has some remarkable features. Indeed, we show that the Yudovich functions
$y_{\Theta_{1}}$ (cf. \eqref{intro:Auxy} and \eqref{intro:AuxTheta1}) associated with $Y_{p_{0}}^{\Theta}$
and $Y_{p_{0}}^{\#\Theta}$ are the same, but now, for every growth\footnote{As usual, the space $\text{BMO}(\Omega)$ in \eqref{Intro1newnew} should be replaced by $\text{BMO}(\Omega) \cap L^{p_0}(\Omega)$ if $\Omega = \R^2$.}  $\Theta$, 
\begin{equation}\label{Intro1newnew}
	Y^{\Theta}_{p_0}(\Omega) \cup \text{BMO}(\Omega) \subset Y^{\#\Theta}_{p_0}(\Omega).
\end{equation}
As a consequence, we establish the  following uniqueness assertion for $Y^{\#\Theta}_{p_0}$.

\begin{theorem}
\label{Thm3Intro} Let $\Omega = \R^2, \T^2$.  Assume that the growth function $\Theta$ satisfies \eqref{intro:osg}. Then a  weak solution
$\omega$ of \eqref{intro1}, such that
\begin{equation}
\label{CThmIn}\omega\in L^{\infty}([0, \infty); Y^{\#\Theta}_{p_{0}}(\Omega))
\qquad\text{for some} \qquad p_{0} \in(2, \infty),
\end{equation}
is uniquely determined by its initial value $\omega_{0}$.
\end{theorem} 
The previous result tells us that the uniqueness vorticity classes $Y^\Theta_{p_0}$ and $\text{BMO}$ considered in \cite{Y95} and \cite{V99}  can be considerably enlarged. Note that \eqref{CThmIn} holds
provided that $\omega_{0}\in Y_{p_{0}}^{\Theta} (\Omega)$ (cf. \eqref{Intro1newnew}), while Theorem \ref{Thm3Intro} with $\Theta(p) \approx 1$ recovers uniqueness in 
\begin{equation}
Y_{p_{0}}^{\#\Theta} (\Omega)=L^{p_{0}} (\Omega)\cap\text{BMO} (\Omega). \label{introYLinftysharp}%
\end{equation}
Hence, we improve the classical result \eqref{introYLinfty} in the sense that
$L^{\infty}(\Omega)$ is now replaced by $\text{BMO}(\Omega)$. 
However, the applicability of Theorem \ref{Thm3Intro} goes further beyond $Y^{\Theta}_{p_0}$ and $\text{BMO}$. Indeed, the
admissible growth $\Theta(p)\approx\log p$ in Theorem \ref{Thm3Intro} provides uniqueness for vorticities of type $$\omega(x)\approx(1+|\log|x||)\log(1+|\log|x||),$$ cf. Section \ref{SectionCounterExamples}, 
Example \ref{Example4}. Note that these vorticities do not belong to
$Y_{p_{0}}^{\Theta}$  (cf. \eqref{PY})  nor to $\text{BMO}$ (they do not even belong to the
larger space\footnote{Recall that, as a consequence of the John--Nirenberg inequality \cite{JN61},  $\text{BMO}$ is (locally) contained in $e^L$. } $e^{L}$). In fact, they grow to infinity faster than both $\log (1 +|\log |x||)$ and $1 + |\log |x||$, which are the prototypes of vorticity in Yudovich's and Vishik's methods, respectively. Furthermore, for a general growth $\Theta$, we propose a simple approach to construct elements in $Y^{\#\Theta}_{p_0}(\Omega)$ that are not in the classical scale $Y^\Theta_{p_0}(\Omega)$; cf. Section \ref{SectionCounterExamples}. The limiting case for the $Y_{p_{0}}^{\#\Theta}$ scale is once again the
linear growth $\Theta(p)\approx p.$ This suggests a possible route to settle
the $L^{p}$ uniqueness problem by means of showing a counterexample belonging
to $Y_{p_{0}}^{\#\Theta}$ with $\Theta(p)\approx p$.

In the forthcoming subsections, we will  motivate and explain how the spaces $Y^{\#\Theta}_{p_0}(\Omega)$ arise naturally in the study of Euler equations and allow us to exploit the techniques of extrapolation theory. 

\subsection{A novel methodology via extrapolation}\label{SubME}

In this paper, we apply the extrapolation theory of Jawerth--Milman
\cite{JM91} to develop a new methodology\footnote{Our presentation was
inspired by the Colombian writer Garc\'{\i}a M\'{a}rquez and his novel
\textquotedblleft Chronicle of a death foretold", vintage publishers, 2003,
that starts with a murder and then develops the plot backwards. In particular,
further explanations and documentation on extrapolation theory will be given
in due time.} that allows us to significantly enlarge the
known uniqueness classes of vorticities for Euler equations and, at the same
time, establish new uniqueness results for a wide class of active scalar
equations. In particular, our method can be applied to establish in a unified fashion uniqueness of weak solutions to $2$D Euler and  Navier-Stokes equations, Keller-Segel models, aggregation equations, as well as more singular equations than $2$D Euler equations such as SQG equations and $2$D incompressible porous media equation (IPM). 


 Indeed, using extrapolation we  are able to construct new
uniqueness spaces that contain BMO, but they are still easy to work with. Our new understanding of Yudovich spaces (cf. Section \ref{SectionIntro1.3}) motivated the construction of extrapolation spaces, where the role played by the pair $(L^{p_{0}},L^{\infty})$ in Yudovich's theory is now replaced by the larger
interpolation pair $(L^{p_{0}},\text{BMO}).$ In concrete terms, new spaces
$Y_{p_{0}}^{\#\Theta}(\Omega)$ (cf. Definition \ref{DefInYS} below) are introduced by
means of replacing $\left\Vert \omega\right\Vert _{L^{p}(\Omega)}$ in
\eqref{intro:YudoSpace} by $\| M_{\Omega}^{\#}\omega\|
_{L^{p}(\Omega)},$ where $M_{\Omega}^{\#}\omega$ denotes the
\emph{Str\"{o}mberg--Jawerth--Torchinsky maximal operator} \cite{JT85}
\begin{equation}
M_{\Omega}^{\#}\omega(x):=\sup_{\Omega\supset Q\ni x}\inf_{c\in\mathbb{R}%
}((f-c)\chi_{Q})^{\ast}(|Q|/\alpha). \label{IntroDefMaxS}%
\end{equation}
Here, $f^{\ast}$ is the \emph{non-increasing rearrangement} of $f$, $Q$ is a cube
with sides parallel to the axes of coordinates, and $\alpha>0$ is a
sufficiently small fixed parameter. While this change worsens the $L^{p}$ norm
by a factor of\footnote{Recall that the classical Fefferman--Stein inequality \cite{FS72} says, loosely speaking, that $L^p$-norms of $\omega$ and $M^{\#}_\Omega \omega$ are comparable, but the equivalence constant deteriorates as $p$ when $p \to \infty$. } $``p"$, this is compensated by the fact that, a quantitative version of the BMO boundedness of CZO, allows us to derive sharp norm estimates for maximal functions of type (\ref{IntroDefMaxS}). In this fashion, extrapolation methods can be applied to exploit the operators involved in the Biot-Savart law:   In
our approach, the interplay between two extrapolation spaces arises. On the
one hand, the extrapolation space related to the vorticities (say $Y_{p_{0}%
}^{\#\Theta}$) and, on the other hand, the extrapolation space in the Besov
scale (e.g. associated with the pair $(L^{\infty},\dot{W}_{\infty}^{1}))$. The
latter controls the modulus of continuity of the associated velocity $v$ according to the Biot-Savart law and informs the
corresponding Osgood condition required to establish uniqueness  (cf. (\ref{intro:osg}) and \eqref{IntroPi}). Once the optimal Osgood condition has been settled, we are able to achieve the desired uniqueness results via an extrapolation refinement of the energy method of \cite{Y95}, in the sense that the role played by $Y^{\Theta}_{p_0}$ in \eqref{KeyYudo} is now replaced by the general extrapolation space that arises  from the Biot-Savart estimate. As already mentioned above, this scheme is very general in the sense that it is not specific of Euler equations, but also works with a large
class of active scalar equations, including Navier-Stokes equations, aggregation equations, Keller-Segel models, SQG, IPM, ...

\subsection{Extrapolation and the role of the Yudovich functions: a preview}

\label{SectionIntro1.3}

To explain the connection of $Y_{p_{0}}^{\#\Theta}$, BMO, and CZO with
extrapolation we need to develop some background information. This subsection
contains some basic definitions with more details and documentation to follow
(cf. Appendix \ref{sec:A0}). Two important aspects regarding the level of
generality that we require should be emphasized here:

\begin{itemize}
\item Since it is crucial for us to consider a variety of scales of function
spaces (e.g., $L^{p}$, BMO, Besov, Sobolev, Yudovich, Vishik) it will be
necessary to formulate the definitions in a sufficiently general context.

\item To deal with different types of scales of function spaces, that measure
different characteristics of their elements (smoothness, integrability,
oscillations, etc.), while at the same time achieving a unified description, the
Peetre $K$-functional \emph{associated with each scale} (cf. \eqref{Peetre} below) will be an invaluable
tool. Indeed, using the $K$-functional the format of the formulae for the
interpolation norms is the same for \emph{all} the scales under consideration, since
it is the particular $K$-functional that contains the quantitative information
associated with the given scale. This explains the \textquotedblleft
universal\textquotedblright\ characterization of extrapolation spaces (cf.
(\ref{form:universal2}) below.)
\end{itemize}

In an informal manner, in interpolation we start with a pair $\bar{A}=(A_{0},A_{1})$ of
compatible\footnote{Roughly speaking, this means that $A_{0}+A_{1}$
\textquotedblleft makes sense\textquotedblright.} Banach spaces and we wish to extract as information as possible on intermediate spaces from the end-points $A_0$ and $A_1$. Let
$(\theta,p)\in(0,1)\times\lbrack1,\infty],$ the \emph{real interpolation
space} $\bar{A}_{\theta,p}$ is the set of all $f\in A_{0}+A_{1}$ such that
\begin{equation*}
\left\Vert f\right\Vert _{\bar{A}_{\theta,p}}:=\bigg\{\int_{0}^{\infty
}[t^{-\theta}K(t,f;A_{0},A_{1})]^{p}\frac{dt}{t}\bigg\}^{1/p}<\infty
\end{equation*}
(with the usual modification if $p=\infty$), where
\begin{equation}\label{Peetre}
K(t,f;A_{0},A_{1}):=\Vert f\Vert_{A_{0}+tA_{1}}=\inf\,\{\left\Vert
f_{0}\right\Vert _{A_{0}}+t\left\Vert f_{1}\right\Vert _{A_{1}}:f=f_{0}%
+f_{1},\quad f_{i}\in A_{i},\quad i=0,1\}
\end{equation}
is the \emph{Peetre $K$-functional} (cf. \cite{BS88, BL76}). It is convenient to normalize the
spaces in order to have a continuous scale with respect to $\theta$. This is
achieved by letting $\bar{A}_{\theta,p}^{\blacktriangleleft}=c_{\theta
,p}\,\bar{A}_{\theta,p}$, where\footnote{If $p=\infty,$ we let $c_{\theta,p}=1.$ Then $\bar{A}_{\theta,\infty}^{\blacktriangleleft} = \bar{A}_{\theta,\infty}$ with equality of norms.}
$c_{\theta,p}:=(\theta(1-\theta)p)^{1/p}$ (cf. \cite{JM91}), and
\begin{equation}
\Vert f\Vert_{\bar{A}_{\theta,p}^{\blacktriangleleft}}:=c_{\theta,p}\,\Vert
f\Vert_{\bar{A}_{\theta,p}}. \label{IntroNormInt}%
\end{equation}

We may view extrapolation as a converse of interpolation: Starting with a family of intermediate spaces we wish to extract information about the possible end-points.  The rigorous definition of the special extrapolation method applied in this paper is as follows.
Let $\Theta$ be a given a growth function, the $\Delta$-\emph{extrapolation space}
$\Delta_{\theta\in(0,1)}\big\{\frac{\bar{A}_{\theta,p(\theta)}%
^{\blacktriangleleft}}{\Theta(\frac{1}{1-\theta})}\big\}$ is defined as the
set of all $f\in\cap_{\theta\in(0,1)}\bar{A}_{\theta,p(\theta)}^{\blacktriangleleft}$ such that
(cf. \cite{JM89, JM91})%
\begin{equation}
\left\Vert f\right\Vert _{\Delta_{\theta\in(0,1)}\big\{\frac{\bar{A}%
_{\theta,p(\theta)}^{\blacktriangleleft}}{\Theta(\frac{1}{1-\theta})}%
\big\}}:=\sup_{\theta\in(0,1)}\frac{\left\Vert f\right\Vert _{\bar{A}%
_{\theta,p(\theta)}^{\blacktriangleleft}}}{\Theta(\frac{1}{1-\theta})}<\infty.
\label{intro:extra}%
\end{equation}

The characterization of these spaces hinges upon the fact that the second
index $p(\theta)$ of interpolation norms can be replaced by $\infty,$ in other words the second index
is not important at the level of normalized norms (cf. \cite{JM91})
\begin{equation}
\Delta_{\theta\in(0,1)}\bigg\{\frac{(A_{0},A_{1})_{\theta,p(\theta
)}^{\blacktriangleleft}}{\Theta(\frac{1}{1-\theta})}\bigg\}=\Delta_{\theta
\in(0,1)}\bigg\{\frac{(A_{0},A_{1})_{\theta,\infty}^{\blacktriangleleft}%
}{\Theta(\frac{1}{1-\theta})}\bigg\}. \label{form:universal}%
\end{equation}
Thus, the commutation of the underlying suprema (\textquotedblleft
Fubini\textquotedblright!) yields\footnote{See the discussion in Appendix
\ref{sec:A1}.}
\begin{equation}
\Delta_{\theta\in(0,1)}\bigg\{\frac{(A_{0},A_{1})_{\theta,p(\theta
)}^{\blacktriangleleft}}{\Theta(\frac{1}{1-\theta})}\bigg\}=\bigg\{f:\sup
_{t\in(0,\infty)}\frac{K(t,f;A_{0},A_{1})}{t\varphi_{\Theta}\big(\frac{1}%
{t}\big)}<\infty\bigg\}, \label{form:universal2}%
\end{equation}
where
\begin{equation}
\varphi_{\Theta}(t):=\inf_{\theta\in(0,1)}\bigg\{\Theta\bigg(\frac{1}{1-\theta
}\bigg)\,t^{1-\theta}\bigg\}. \label{laphijm}%
\end{equation}

For the remaining of this subsection, we shall place ourselves under the
conditions of \cite{Y95} and focus our attention on the spaces $Y_{p_{0}}^{\Theta}(\Omega), \, \Omega= \T^2$ (although there are analogous statements  for $\Omega = \R^2$ or even general domains $\Omega \subset \R^d$). It is
known that, with equivalence of norms independent of $\theta$,
\begin{equation}
(L^{p_{0}},L^{\infty})_{\theta,p(\theta)}^{\blacktriangleleft}=L^{p(\theta
)},\qquad\frac{1}{p(\theta)}=\frac{1-\theta}{p_{0}} \label{IntroLp}%
\end{equation}
(cf. \cite[eq. (31), p. 61]{KM}). Using this fact, together with the
monotonicity properties of the Lebesgue scale and \eqref{form:universal},
yields
\begin{equation}
Y_{p_{0}}^{\Theta}   =\Delta_{\theta\in(0,1)}\bigg\{\frac{(L^{p_{0}%
} ,L^{\infty})_{\theta,p(\theta)}^{\blacktriangleleft}}{\Theta(\frac{p_{0}%
}{1-\theta})}\bigg\}  =\Delta_{\theta\in(0,1)}\bigg\{\frac{(L^{p_{0}},L^{\infty}
)_{\theta,\infty}^{\blacktriangleleft}}{\Theta(\frac{p_{0}}{1-\theta}%
)}\bigg\}. \label{IntroYExtra}%
\end{equation}
Furthermore, comparing (\ref{intro:Auxy}) and (\ref{laphijm}) via the change
$p\leftrightarrow\frac{p_{0}}{1-\theta},$ we find that, with constants
depending only on $p_{0}$,
\begin{equation}
\varphi_{\Theta}(t)\approx y_{\Theta}(t). \label{laphidy}%
\end{equation}
To give a full explicit characterization of $Y_{p_{0}}^{\Theta}$  we consider first the case $p_{0}=1$. Recall that (cf. \eqref{A9new}) $$K(t,f;L^{1},L^{\infty})=\int_{0}%
^{t}f^{\ast}(s)\,ds,$$ is closely related to the classical maximal function $f^{\ast\ast}(t):=\frac{1}{t}\int%
_{0}^{t}f^{\ast}(s)\,ds$. This information combined with
\eqref{form:universal}, \eqref{form:universal2}, (\ref{IntroYExtra}) and
(\ref{laphidy}) yields\footnote{Since $\left\vert \Omega\right\vert <\infty,$
it is not hard to see that we can restrict to $t$ $\in(0,1).$ Indeed, observe that  $\int_0^t f^*(s) \, ds = \int_0^1 f^*(s) \, ds$ if $t > 1$ and $\sup_{t>1}%
\frac{1}{t\varphi_{\Theta}(\frac{1}{t})}<\infty$ because $\Theta$ is a
non-decreasing function.}
\begin{equation}
\left\Vert f\right\Vert _{Y_{p_{0}}^{\Theta}}\approx\sup
_{t \in (0, \infty)}\,\frac{\int_{0}^{t}f^{\ast}(s)\,ds}{t\varphi_{\Theta}(\frac{1}{t}%
)}\approx\sup_{t\in(0,1)}\,\frac{f^{\ast\ast}(t)}{y_{\Theta}(\frac{1}{t}%
)}\approx\sup_{t\in(0,1)}\,\frac{f^{\ast}(t)}{y_{\Theta}(\frac{1}{t})},
\label{explica2}%
\end{equation}
where the last equivalence follows from \eqref{ClaimInitial} below. We now show that for the growth functions considered in this paper we have $Y_{p_{0}}^{\Theta}=Y_{1}^{\Theta}$, $p_{0}>1$. Indeed, we can use the formula for the $K$-functional for the pair $(L^{p_0},L^{\infty})$ given in Example {\ref{Ex2.a}} and then proceed as in the proof of \eqref{CharSec2.1new} below. Alternatively, the result follows directly from \eqref{app2} applied to $(L^{p_0},L^{\infty})$. As a consequence,  Yudovich spaces can be identified with the more familiar \emph{Marcinkiewicz spaces}, that have been extensively studied in the literature (cf. \cite{KPS82,
BS88}).
\begin{example}
Let $\Theta(p) \approx p,$ then $Y^{\Theta}_{p_{0}}=e^{L}$ (cf. \cite{JM89, JM91}) with
\[
\left\Vert f\right\Vert _{Y^{\Theta}_{p_{0}}} \approx
\sup_{t\in(0,1)} \, \frac{f^{\ast\ast}(t)}{1-\log t} \approx\sup_{t\in(0,1)}
\, \frac{f^{\ast}(t)}{1-\log t}.
\]
\end{example}
\begin{example}\label{ExIntroNew}
More generally, suppose that $\Theta$ is a growth function such that
$p\in(p_{0},\infty)\mapsto e^{p_{0}/p}\Theta(p)$ is a quasi-decreasing
function\footnote{Some examples are $\Theta(p) \approx p^\alpha (\log p)^{\alpha_1} (\log_2 p)^{\alpha_2} \cdots (\log_m p)^{\alpha_m}$, where $\alpha, \alpha_i \in \R$ and $\log_m p = \underbrace{\log \ldots \log}_\text{m times} p$ for $m \geq 2$.} (i.e., equivalent to a decreasing function), then
\[
\Vert f\Vert_{Y_{p_{0}}^{\Theta}}\approx\sup_{t\in
(0,1)}\,\frac{f^{\ast\ast}(t)}{\Theta(1-\log t)}\approx\sup_{t\in(0,1)}%
\,\frac{f^{\ast}(t)}{\Theta(1-\log t)}.
\]
This is a consequence of \eqref{explica2} and Lemma \ref{Lemma1} below.
\end{example}

\subsection{A priori Biot-Savart estimates via extrapolation}

\label{SectionIntroAPrior}

It will be instructive to revisit the main uniqueness result of \cite{Y95}
using our method. As a first step we use the Biot-Savart law to obtain a
priori estimates for the modulus of continuity of the solutions.

To simplify the exposition, we work again with $\Omega = \T^2$ (but similar results hold for $\Omega = \R^2$ or $\Omega$ a smooth domain in $\R^2$, or even in $\R^d$). Fix $p_{0}>2$.  By \eqref{IntroLp}, (\ref{intro2}), the sharp $L^{p}$ norm
inequalities for CZO, and the definition of $Y_{p_{0}}^{\Theta},$ we derive
\begin{equation}
\left\Vert \nabla v\right\Vert _{(L^{p_{0}},L^{\infty})_{\theta,p(\theta
)}^{\blacktriangleleft}}\leq c_{p_{0}}\,\frac{1}{1-\theta}\,\Theta
\bigg(\frac{1}{1-\theta}\bigg)\left\Vert \omega\right\Vert _{Y_{p_{0}}%
^{\Theta}}. \label{update}%
\end{equation}
The well-understood interpolation theory of Sobolev spaces (cf. \cite{DS79}) allows us to
rewrite the left-hand side as 
\begin{equation}
\left\Vert \nabla v\right\Vert _{(L^{p_{0}},L^{\infty})_{\theta,p(\theta
)}^{\blacktriangleleft}}\approx\left\Vert v\right\Vert _{(\dot{W}_{p_{0}}%
^{1} ,\dot{W}_{\infty}^{1})_{\theta,p(\theta)}^{\blacktriangleleft}}
\label{IntroIntSobDS}%
\end{equation}
with equivalence constants independent of $\theta$. 
Then inserting this information in (\ref{update}), and rewriting the right-hand side using the definition of $\Theta_{1}$ (cf. \eqref{intro:AuxTheta1}),
we obtain%
\[
\left\Vert v\right\Vert _{(\dot{W}_{p_{0}}^{1},\dot{W}_{\infty}^{1}
)_{\theta,p(\theta)}^{\blacktriangleleft}}\leq c_{p_{0}}\Theta_{1} \bigg(\frac
{1}{1-\theta} \bigg)\left\Vert \omega\right\Vert _{Y_{p_{0}}^{\Theta}}.
\]
Consequently,%
\begin{equation}
\left\Vert v\right\Vert _{\Delta_{\theta\in(0,1)}\Big\{\frac{(\dot{W}_{p_{0}%
}^{1},\dot{W}_{\infty}^{1})_{\theta,p(\theta)}^{\blacktriangleleft}}%
{\Theta_{1}(\frac{1}{1-\theta})}\Big\}}=\sup_{\theta\in(0,1)}\frac{\left\Vert
v\right\Vert _{(\dot{W}_{p_{0}}^{1},\dot{W}_{\infty}^{1})_{\theta,p(\theta
)}^{\blacktriangleleft}}}{\Theta_{1}(\frac{1}{1-\theta})} \lesssim \left\Vert \omega\right\Vert _{Y_{p_{0}}^{\Theta}}. \label{puntojap}%
\end{equation}
The Sobolev embedding theorem (recall $p_{0}>d = 2$, where $d$ denotes the dimension of the ambient space) combined with (\ref{puntojap}%
) yields%
\[
\left\Vert v\right\Vert _{\Delta_{\theta\in(0,1)}\Big\{\frac{(L^{\infty}
,\dot{W}_{\infty}^{1})_{\theta,p(\theta)}^{\blacktriangleleft}}{\Theta
_{1}(\frac{1}{1-\theta})}\Big\}}\lesssim \left\Vert \omega\right\Vert
_{Y_{p_{0}}^{\Theta}}.
\]
The extrapolation norm that appears on the left-hand side can be computed
explicitly using (\ref{form:universal2}), \eqref{laphidy} and the well-known
fact that (cf. \eqref{A11novel})
\begin{equation}
K(t,v;L^{\infty},\dot{W}_{\infty}^{1})\approx\sup_{\left\vert x-y\right\vert
\leq t}\left\vert v(x)-v(y)\right\vert . \label{IntroKFunctModLinfty}%
\end{equation}
It then follows that 
\begin{equation}
\left\vert v(x)-v(y)\right\vert \lesssim\left\vert x-y\right\vert
y_{\Theta_{1}}\bigg(\frac{1}{\left\vert x-y\right\vert }\bigg)\left\Vert
\omega\right\Vert _{Y_{p_{0}}^{\Theta}}, \label{Intro:EstimModYtheta}%
\end{equation}
where%
\begin{equation}
y_{\Theta_{1}}(t)=\inf_{\theta\in(0,1)}\,\bigg\{\Theta_{1}\bigg(\frac
{1}{1-\theta}\bigg)\,t^{1-\theta}\bigg\}. \label{Introytheta}%
\end{equation}
This is the main estimate of the modulus of continuity given in \eqref{18a}. To obtain the corresponding result of \cite{Y63}, let $\Theta (p) \approx 1.$
Then, $Y_{p_{0}}^{\Theta}=L^{\infty}$ (cf. (\ref{introYLinfty}%
)), and $y_{\Theta_{1}}(t)\approx\log t$ if $t>1,$ yielding%
\begin{equation}\label{18b}
\left\vert v(x)-v(y)\right\vert \lesssim\left\vert x-y\right\vert
|\log|x-y||\left\Vert \omega\right\Vert _{L^{\infty}}
\end{equation}
if $|x-y| < 1$. See also \cite[Section 4.1]{L96}.


\subsection{Uniqueness of weak solutions}

\label{SectionIntro4}

Inequalities of type \eqref{Intro:EstimModYtheta} and \eqref{18b} contain precious information about uniqueness of the particle path. Specifically, it is well known (see e.g. \cite{BCD11}) that  the flow map $\phi$ relative to $v$ according to
\begin{equation}
\frac{d}{dt}\phi(t,x)=v(t,\phi(t,x)),\qquad\phi(0,x)=x, \label{Flow}%
\end{equation}
is uniquely defined  provided that
\begin{equation}\label{18.4}
|v(t,x)-v(t,y)|\leq L(|x-y|)\,\Vert\omega(t)\Vert_{X}\qquad\text{for
a.e.}\quad t\in(0,\varepsilon_{L}), 
\end{equation}
where $L:(0,\varepsilon_{L})\rightarrow(0,\infty)$ is  a continuous
nondecreasing function, 
$\varepsilon_{L}\in(0,\infty]$, such that
\begin{equation}\label{18.2}
\Vert\omega(t)\Vert_{X}\in L_{\text{loc}}^{1}(0,\varepsilon_{L}),
\end{equation}
and $L$ satisfies the Osgood condition
\begin{equation}\label{18.3}
\int_{0}^{\varepsilon_{L}}\frac{dr}{L(r)}=\infty. 
\end{equation}
Hence, at least in some sense,  \eqref{18.4} allows to identify optimal candidates of function spaces $X$ (cf. \eqref{18.2}) and Osgood conditions  (cf. \eqref{18.3}) that guarantee uniqueness of the flow map $\phi$ in \eqref{Flow}. However, uniqueness of $\phi$ is still not strong enough to claim uniqueness of weak solutions $v$ to Euler equations, and further arguments\footnote{We owe this comment to a reviewer of an earlier version of the paper.} are needed to control the difference $v_1-v_2$  of two weak solutions $v_1$ and $v_2$. Indeed, in this paper we show a simple general mechanism to derive uniqueness of weak solutions in terms of  extrapolation spaces $X$ satisfying  \eqref{18.4} with \eqref{18.3}.

\subsection{Summary of our approach to uniqueness}

\label{SectionIntroSum} The method described above can be summarized as follows.

\vspace{1mm}

\textbf{Step 1:} Compute the extrapolation spaces involved.

\vspace{1mm} \textbf{Step 2:} Use the Biot-Savart laws to obtain a priori
estimates of the smoothness of solutions via extrapolation.

\vspace{1mm}

\textbf{Step 3: }Prove uniqueness via an extrapolation-based energy method relative to the extrapolation spaces and the Osgood conditions arising in Steps 1 and 2, respectively.

\subsection{Extending Yudovich's uniqueness theorem beyond $\text{BMO}$}

\label{SeIntroBMOY}

As already announced in Theorem \ref{Thm3Intro} and the subsequent discussion, in this paper we extend $\text{BMO}$ uniqueness to spaces that contain \emph{functions of
unbounded mean oscillation}. This will be done through the introduction of the
new (extrapolation) spaces $Y_{p_{0}}^{\#\Theta}$.

\begin{definition}
[Sharp Yudovich spaces]\label{DefInYS} Let\footnote{The definition can
obviously be given in a more general context.} $\Omega=\mathbb{R}^{2}%
,\,\mathbb{T}^{2}$ and let $p_{0}\in\lbrack1,\infty)$. Given a growth function
$\Theta$, the \emph{sharp Yudovich space} $Y_{p_{0}}^{\#\Theta}(\Omega)$ is defined
to be the set of all\footnote{$\Vert f\Vert_{(L^{p})^{\#}(\Omega)}:=\Vert
M^{\#}_\Omega f\Vert_{L^{p}(\Omega)}.$ In particular, $(L^{\infty})^{\#}%
(\Omega)=\text{BMO}(\Omega)$ and $(L^{p})^{\#}(\Omega)=L^{p}(\Omega
),\,p\in(1,\infty),$ suitably interpreted modulo constants; cf.
\cite[Corollaries 2.5 and 2.6]{JT85}.} $f\in\cap_{p>p_{0}}(L^{p})^{\#}%
(\Omega)$ such that
\[
\left\Vert f\right\Vert _{Y_{p_{0}}^{\#\Theta}(\Omega)}:=\sup_{p>p_{0}}%
\frac{\left\Vert f\right\Vert _{(L^{p})^{\#}(\Omega)}}{\Theta(p)}<\infty.
\]

\end{definition}

As in the case of Yudovich spaces, the definition of $Y_{p_{0}}^{\#\Theta
}(\Omega),$ with $\Omega=\mathbb{T}^{2},$ does not depend on $p_{0}$. When $\Omega=\mathbb{R}^{2}$, we only
have the trivial embeddings
\begin{equation*}
Y_{p_{0}}^{\#\Theta}(\Omega)\hookrightarrow Y_{p_{1}}^{\#\Theta}%
(\Omega),\qquad\text{if}\qquad p_{1}>p_{0}. 
\end{equation*}
For example, when $\Theta(p)\approx1$, $$Y_{p_{0}}^{\#\Theta}(\mathbb{T}%
^{2})=\text{BMO}(\mathbb{T}^{2}) \qquad \text{and} \qquad Y_{p_{0}}^{\#\Theta}(\mathbb{R}%
^{2})=\text{BMO}(\mathbb{R}^{2}) \cap L^{p_0}(\R^2),$$
cf. \eqref{introYLinftysharp}. On the other hand, since $\Vert
f\Vert_{(L^{p})^{\#}(\Omega)}\lesssim\Vert f\Vert_{L^{p}(\Omega)}%
,\,1<p_{0}<p\leq\infty$, we have
\begin{equation}
Y_{p_{0}}^{\Theta}(\Omega)\hookrightarrow Y_{p_{0}}^{\#\Theta}(\Omega).
\label{IntroEmYshY}%
\end{equation}

In Section \ref{Section2.1} we show that $Y_{p_{0}}^{\#\Theta}$ fits into the
abstract extrapolation framework proposed above. In this context, $\text{BMO}$
plays the same role as $L^{\infty}$ in connection with $Y_{p_{0}}^{\Theta}$.
In particular, we obtain characterizations of $Y_{p_{0}}^{\#\Theta}$
corresponding to those for $Y_{p_{0}}^{\Theta}$ (cf. Table \ref{Tabla2}).

\begin{table}[h]
\centering
\begin{tabular}
[c]{|c|c|}\hline
$Y^{\Theta}_{p_{0}}$ & $Y^{\#\Theta}_{p_{0}}$\\\hline
& \\
$\Big\{f \in\bigcap_{p > p_{0}} L^{p}: \sup_{p > p_{0}} \frac{\left\Vert
f\right\Vert _{L^{p}}}{\Theta(p)} < \infty\Big\}$ & $\Big\{f \in\bigcap_{p >
p_{0}} (L^{p})^{\#}: \sup_{p > p_{0}} \frac{\left\Vert f\right\Vert
_{(L^{p})^{\#}}}{\Theta(p)} < \infty\Big\}$\\
& \\\hline
& \\
$\Delta_{\theta\in(0,1)} \, \Big\{\frac{(L^{p_{0}},L^{\infty})_{\theta
,p(\theta)}^{\blacktriangleleft}}{\Theta(\frac{p_{0}}{1-\theta})} \Big\}$ &
$\Delta_{\theta\in(0,1)} \, \Big\{\frac{(L^{p_{0}},\text{BMO})_{\theta
,p(\theta)}^{\blacktriangleleft}}{\Theta(\frac{p_{0}}{1-\theta})} \Big\}$\\
& \\\hline
& \\
$\Big\{f : \sup_{t \in(0, \infty)}\frac{K(t^{1/p_{0}},f;L^{p_{0}},L^{\infty}%
)}{t^{1/p_{0}} y_{\Theta}(\frac{1}{t})}<\infty\Big\} $ & $\Big\{f : \sup_{t
\in(0, \infty)}\frac{K(t^{1/p_{0}},f;L^{p_{0}},\text{BMO})}{t^{1/p_{0}}
y_{\Theta}(\frac{1}{t})}<\infty\Big\} $\\
& \\\hline
& \\
$\Big\{f: \sup_{t \in(0, \infty)}\frac{(|f|^{p_{0}})^{\ast\ast}(t)^{1/p_{0}}%
}{y_{\Theta}(\frac{1}{t})} < \infty\Big\}$ & $\Big\{f: \sup_{t \in(0, \infty
)}\frac{(|M^{\#} f|^{p_{0}})^{**}(t)^{1/p_{0}}}{y_{\Theta}(\frac{1}{t})} <
\infty\Big\}$\\
& \\\hline
\end{tabular}
\caption{Yudovich vs. sharp Yudovich}%
\label{Tabla2}%
\end{table}

The usefulness of $Y_{p_{0}}^{\#\Theta}$ emerges when establishing a priori
estimates for the modulus of continuity of the velocity $v$ in terms of the vorticity $\omega$ (cf. \eqref{intro1}).

\begin{theorem}
\label{Theorem2}  Assume that $\omega\in Y_{p_{0}%
}^{\#\Theta}(\Omega)$ for some $p_{0}\in(2,\infty)$. Then
\begin{equation}
|v(x)-v(y)|\lesssim|x-y|\,y_{\Theta_{1}}\bigg(\frac{1}{|x-y|}\bigg)\,\Vert
\omega\Vert_{Y_{p_{0}}^{\#\Theta}(\Omega)}, \label{introMBMO}%
\end{equation}
where $y_{\Theta_{1}}$ is given by \eqref{Introytheta}.
\end{theorem}

In particular, \eqref{introMBMO} extends the classical Yudovich's estimate
\eqref{Intro:EstimModYtheta} from $Y^{\Theta}_{p_{0}}$ to $Y^{\#\Theta}%
_{p_{0}}$. Then, in Section \ref{SectionUY}, we show the
desired uniqueness result for $Y^{\#\Theta}_{p_{0}}$ stated in Theorem \ref{Thm3Intro}. The proof of this result (in fact, two different proofs will be presented, one of them has the advantage of working not only with Euler equations, but also with the more general class of active scalars, cf. \eqref{Trans} below) will be given in Section \ref{Section22}. 
%
%

\subsection{Vishik's uniqueness theorem for active scalar equations}\label{SVAS}

As already claimed  in Section \ref{SubME}, we show that Vishik spaces
$B_{\Pi}(\Omega)$ (cf.
\eqref{introVishikSpace}) are also special examples of extrapolation
constructions in the sense of \eqref{intro:extra}. In fact, we obtain several
characterizations of Vishik spaces by a variety of means (extrapolation, interpolation, Yudovich functions, and growths of classical
Besov norms) thus showing the full analogy between $B_\Pi$ and $Y^{\Theta}_{p_0}$. The results are collected in Table \ref{Tabla1} (where, for
simplicity, we let once again $p_{0}=1$ in the definition of $Y_{p_{0}}^{\Theta}$), the corresponding proofs may be found in Section \ref{SectionCharacVishik} below.

\begin{table}[h]
\centering
\begin{tabular}
[c]{|c|c|}\hline
$Y^{\Theta}_{1}$ & $B_{\Pi}$\\\hline
& \\
$\Big\{f \in\bigcap_{p > 1} L^{p}: \sup_{p > 1} \frac{\left\Vert
f\right\Vert _{L^{p}}}{\Theta(p)} < \infty\Big\}$ & $\Big\{f
\in\bigcap_{\alpha\in(-1, 0)} B^{\alpha}_{\infty, 1} : \sup_{\alpha
\in(-1, 0)} \frac{\|f\|_{B^{\alpha}_{\infty, 1}}}{\Pi(-\frac{1}%
{\alpha})} < \infty\Big\}$\\
& \\\hline
& \\
$\Delta_{\theta\in(0,1)} \, \bigg\{\frac{(L^{1},L^{\infty}%
)_{\theta,\frac{1}{1-\theta}}^{\blacktriangleleft}}{\Theta(\frac
{1}{1-\theta})} \bigg\}$ & $\Delta_{\theta\in(0, 1)} \bigg\{\frac
{(B^{-1}_{\infty, 1}, B^{0}_{\infty, 1})_{\theta, \frac
{1}{1-\theta}}^{\blacktriangleleft}}{\Pi(\frac{1}{1-\theta})} \bigg\}$\\
& \\\hline
& \\
$\Big\{f : \sup_{t \in(0, \infty)}\frac{K(t,f;L^{1},L^{\infty}%
)}{t y_{\Theta}(\frac{1}{t})}<\infty\Big\} $ & $\Big\{f: \sup_{t
\in(0, \infty)} \frac{K(t, f; B^{-1}_{\infty, 1}, B^{0}_{\infty,
1})}{t y_{\Pi}(\frac{1}{t})} < \infty\Big\} $\\
& \\\hline
& \\
$\Big\{f: \sup_{0<t<\infty}\frac{f^{\ast\ast}(t)}{y_{\Theta}(\frac{1}{t})} <
\infty\Big\}$ & $\Big\{ f : \sup_{N \geq0} \frac{1}{\Pi(N)} \, \sum
\limits_{j=0}^{N}\left\Vert \Delta_{j}f\right\Vert _{L^{\infty}} <
\infty\Big\}$\\
& \\\hline
\end{tabular}
\caption{Yudovich vs. Vishik}%
\label{Tabla1}%
\end{table}

Having at hand the information contained in Table \ref{Tabla1}, we are in a
position to apply the extrapolation approach to uniqueness developed in
Section \ref{SectionIntroSum}. Next we explain this in more detail.

So far, we have discussed the role played by a variety of function spaces in
connection with $2$D Euler equations. However, our method is robust enough to  deal
with the wider class of active scalar equations modeled by
\begin{equation}
\omega_{t}+V\omega\cdot\nabla\omega=0,\label{Trans}%
\end{equation}
where, $\omega=\omega(t,x),\,x\in \Omega, t>0,$
is a scalar function and the vector field $V$ satisfies suitable minimal assumptions. To simplify the exposition,  in the rest of this subsection we concentrate on the case $\Omega = \T^2$, but all the stated results admit analogs for $\Omega = \T^d, \R^d$. We assume that $V : \dot{H}^{\beta-1}(\Omega; \R) \to L^2(\Omega; \R^2)$ is a linear operator such that $\nabla \cdot Vf = 0$ (i.e., divergence free), $V$ commutes with $\Delta_j$ and satisfies
\begin{equation}\label{E12n}
	\|V \Delta_j f\|_{L^\infty} \lesssim 2^{j (\beta-1)} \|\Delta_j f\|_{L^\infty}, \qquad \forall j \geq 0, \qquad \forall f \in \mathcal{S'}. 
\end{equation}
Perhaps the most challenging case in \eqref{E12n} occurs when $\beta=1$. This is indeed the case for the surface quasi-geostrophic (SQG) equations\footnote{As usual, $R_1$ and $R_2$ denote the two components of the Riesz transform in $\R^2$.} ($V \omega = (-R_2 \omega, R_1 \omega)$) and $2$D incompressible porous media (IPM) equations ($V \omega = (-R_2, R_1) R_1 \omega$). The reader is referred to \cite{B72}, and the references therein, for the physical background and applications. Instead here we just mention that they are strongly connected with $3$D Euler equations, with which they share a remarkable number of similarities (cf. \cite{CMT94} and \cite{Kis}).

In this paper we establish the following uniqueness result for \eqref{Trans} in terms of Vishik spaces. For the corresponding notation, cf.  Definitions \ref{DefVishik} and \ref{DefPClass} in Section \ref{SectionVishik}. 

\begin{theorem}
\label{Thm35a}
 Assume that the
growth function $\Pi\in\mathcal{P}_{\kappa}, \, \kappa > 0,$ satisfies 
\begin{equation}\label{introOsB}
\int_{1}^{\infty}\frac{dr}{r y_{\Pi}(r)} = \infty.
\end{equation}
Then a weak solution $\omega$ of \eqref{Trans} with $V$ satisfying \eqref{E12n} under $\beta =1$ (in particular, this applies to SQG and IPM), such that
\[
\omega\in L^{\infty}([0, T]; B^{1}_{\Pi})
\]
is uniquely determined by its initial value $\omega_0$.
\end{theorem}

In particular, Theorem \ref{Thm35a} goes further beyond the usual  uniqueness regularity $H^2$ in SQG or even the borderline regularity $\nabla \omega \in \text{BMO}$ (cf. \cite{AB}). On the other hand, since $C^1 \hookrightarrow B^1_\Pi$ when $\Pi(r) \approx r$, Theorem \ref{Thm35a} can be automatically applied to achieve uniqueness of IPM within the class of \emph{classical} solutions. This class of solutions is non empty, as has been shown recently by C\'ordoba and Mart\'inez-Zoroa \cite{DLu}, where smooth solutions to (forced) IPM exhibiting finite time singularities are constructed, and their explicit behaviour used to prove uniqueness. In contrast, Theorem \ref{Thm35a} can be applied to establish uniqueness of classical solution, without any additional information on its structure. 

The  method of proof of Theorem \ref{Thm35a}, which is based again on extrapolation, is  remarkably simple, and avoids the known sophisticated estimates that rely on heavy harmonic analysis techniques:  Littlewood--Paley theory and Bony's paraproduct decompositions (cf. \cite{BCD11, Co, V99}). It is not clear to us if these harmonic analysis techniques can be adapted to the general setting considered in this paper.

 As already mentioned above, Theorem \ref{Thm35a} works not only for SQG or IPM, but can be also applied to important variants of $2$D Euler equations. More precisely, we can cover a wide family of active scalar equations of type:
\begin{equation}\label{16.9}
	\omega_t + V \omega \cdot \nabla \omega = \Delta A( \omega),
\end{equation}
where $A: \R \to [0, \infty)$ is measurable and non-decreasing and $V$ satisfies \eqref{E12n} with $\beta =0$. Model examples of \eqref{16.9} are  $2$D Euler equations ($A \equiv 0$) and Navier-Stokes equations ($A \equiv \nu \, \text{id}$ and $\nu > 0$ is a fixed constant) where $V \omega = \nabla^\perp \psi$ and $\Delta \psi = \omega$.  Other important examples of \eqref{16.9} include the Keller-Segel model ($V \omega = -\nabla \Delta^{-1} \omega$) and the aggregation equations with Newtonian potentials.  

Next we place the classical Vishik uniqueness theorem (cf. \eqref{IntroViSpa}) within the general setting of active scalars of type \eqref{16.9}. In particular, this provides an extension of \cite{AB} and \cite{Co} that respectively dealt with $\omega \in \text{BMO}$ and $\omega \in B^0_{\infty, \infty}$, to the larger class  $\omega \in B_\Pi$.

\begin{theorem}
\label{IntroThm6} Suppose that $\Pi$ satisfies the conditions of Theorem  \ref{Thm35a}. 
Then a weak solution $\omega$ of \eqref{16.9} such that
\begin{equation}
\omega\in L^{\infty}([0,T];B_{\Pi}) \label{IntroOSB2}%
\end{equation}
is uniquely determined by its initial value $\omega_0$.
\end{theorem}

\begin{remark}
\label{Remark5}

Note that Theorem \ref{IntroThm6} contains as a distinguished special case Vishik's uniqueness theorem (cf. \eqref{IntroPi}-\eqref{IntroViSpa}) for Euler equations.  Indeed, from the trivial estimate\footnote{In fact, under natural assumptions on $\Pi$, we have $y_{\Pi
}(r)\approx\Pi(\log r)$; cf. Lemma \ref{Lemma1} below. Then, a simple change of variables
gives $\int_{1}^{\infty}\frac{dr}{\Pi(r)}\approx\int_{1}^{\infty}\frac
{dr}{ry_{\Pi}(r)}$.} $y_{\Pi}(r)\lesssim\Pi(\log r)$ it follows that
\[
\int_{1}^{\infty}\frac{dr}{\Pi(r)}\lesssim\int_{1}^{\infty}\frac{dr}{ry_{\Pi
}(r)}.
\]
Therefore, the validity of the Osgood condition \eqref{IntroPi}
implies \eqref{introOsB}. Furthermore, the
exchange $\Theta_{1}\leftrightarrow\Pi$ shows a connection with \eqref{intro:osg}.

\end{remark}

In order to facilitate the reading for non experts, we close the paper with an Atlas on Interpolation and Extrapolation (cf. Appendix \ref{sec:A0}), where we collect documentation and supplementary material. 

We believe that the techniques developed in this paper could be useful in other related contexts.

\subsection{Road map} We simply remark that the previous discussion and the table of contents show the local organization of the paper.

\section{The spaces $Y^{\#\Theta}_{p_{0}}$}

\label{Section2}

\subsection{Characterizations}

\label{Section2.1} We establish several characterizations of $Y^{\#\Theta
}_{p_{0}}$ (cf. Definition \ref{DefInYS}) in terms of extrapolation of 
interpolation scales and $K$-functionals involving $\text{BMO}$, and maximal functions.

\begin{theorem}
[Characterization via extrapolation]\label{ThmExtSh} Let $p_{0}\in
\lbrack1,\infty),\,\theta\in(0,1),$ and let $\Omega=\mathbb{R}^{2}%
,\mathbb{T}^{2}$. We have
\begin{equation}
(L^{p_{0}}(\Omega),\emph{BMO}(\Omega))_{\theta,p(\theta)}^{\blacktriangleleft
}=(L^{p(\theta)})^{\#}(\Omega),\qquad\frac{1}{p(\theta)}=\frac{1-\theta}%
{p_{0}}, \label{ClaimEx1}%
\end{equation}
with equivalence constants independent of $\theta$. As a consequence,
\begin{equation}
Y_{p_{0}}^{\#\Theta}(\Omega)=\Delta_{\theta\in(0,1)}\bigg\{\frac{(L^{p_{0}%
}(\Omega),\emph{BMO}(\Omega))_{\theta,p(\theta)}^{\blacktriangleleft}}%
{\Theta(\frac{p_{0}}{1-\theta})}\bigg\}. \label{ClaimEx2}%
\end{equation}

\end{theorem}

\begin{proof}
Using the Jawerth--Torchinsky formula (cf. \eqref{A10novel})
\begin{equation}
\label{JTFor}K(t, f; L^{p_{0}}(\Omega), \text{BMO}(\Omega)) \approx\bigg(\int_{0}^{t^{p_{0}}}
[(M^{\#}_\Omega f)^{*}(\xi)]^{p_{0}} \, d \xi\bigg)^{1/p_{0}},
\end{equation}
where the maximal function $M^{\#}_\Omega$ is defined by \eqref{IntroDefMaxS}, we
have
\begin{equation}
\label{PrS1.1}\|f\|_{(L^{p_{0}}(\Omega),\text{BMO}(\Omega))_{\theta,p(\theta)}} \approx
\bigg\{\int_{0}^{\infty} \bigg(\frac{1}{t} \int_{0}^{t} [(M_\Omega^{\#} f)^{*}%
(\xi)]^{p_{0}} \, d \xi\bigg)^{p(\theta)/p_{0}} \, dt\bigg\}^{1/p(\theta)}.
\end{equation}

Applying the sharp version of Hardy's inequality (note that $p(\theta
)=\frac{p_{0}}{1-\theta}>p_{0}$) stated in \cite[Appendix A.4, page 272]{Stein}, we can
estimate \eqref{PrS1.1} as follows
\begin{equation}
\Vert f\Vert_{(L^{p_{0}}(\Omega),\text{BMO}(\Omega))_{\theta,p(\theta)}}\lesssim
\theta^{-1/p_{0}}\Vert f\Vert_{(L^{p(\theta)})^{\#}(\Omega)}. \label{PrS1.2}%
\end{equation}
Conversely, one can invoke the reverse Hardy inequality
(cf. \cite{M97} and the references therein). Indeed
\begin{align}
\Vert f\Vert_{(L^{p(\theta)})^{\#}(\Omega)}  &  =\Vert(M^{\#}f)^{\ast}\Vert
_{L^{p(\theta)} (0, \infty)}=\Vert(M^{\#}f)^{p_{0}\ast}\Vert_{L^{p(\theta)/p_{0}}%
(0, \infty)}^{1/p_{0}}\nonumber\\
&  \lesssim\theta^{1/p(\theta)}\,\Vert\lbrack(M^{\#}f)^{p_{0}}]^{\ast\ast
}\Vert_{L^{p(\theta)/p_{0}} (0, \infty)}^{1/p_{0}}\approx\theta^{1/p(\theta)}\,\Vert
f\Vert_{(L^{p_{0}}(\Omega),\text{BMO}(\Omega))_{\theta,p(\theta)}}, \label{PrS1.3}%
\end{align}
where in the last estimate we have used \eqref{PrS1.1}. Combining
(\ref{PrS1.2}) and (\ref{PrS1.3}) we obtain%
\begin{align*}
\Vert f\Vert_{(L^{p_{0}} (\Omega),\text{BMO}(\Omega))_{\theta,p(\theta)}^{\blacktriangleleft}%
} &\lesssim c_{\theta p(\theta)}\theta^{-1/p_{0}}\Vert f\Vert_{(L^{p(\theta
)})^{\#} (\Omega)} \\
& \hspace{-2.5cm} \lesssim\theta^{-\theta/p_0}\,\Vert f\Vert_{(L^{p_{0}} (\Omega),\text{BMO}(\Omega)%
)_{\theta,p(\theta)}^{\blacktriangleleft}} \approx \Vert f\Vert_{(L^{p_{0}} (\Omega),\text{BMO}(\Omega)%
)_{\theta,p(\theta)}^{\blacktriangleleft}}.
\end{align*}
Now, (\ref{ClaimEx1}) follows readily since $c_{\theta p(\theta)}%
=(\theta(1-\theta)p(\theta))^{1/p(\theta)}\approx\theta^{1/p(\theta)}%
\approx\theta^{1/p_{0}}$.
\end{proof}


\begin{theorem}
[Characterization via $K$-functional and maximal function]\label{Thm2.2} Let
$p_{0}\in\lbrack1,\infty)$ and  $\Omega=\mathbb{R}^{2},\mathbb{T}^{2}$. 

\begin{enumerate}
\item[(i)] We have
\begin{equation}
\Vert f\Vert_{Y_{p_{0}}^{\#\Theta}(\Omega)}\approx\sup_{t\in(0,\infty)}%
\frac{K(t^{1/p_{0}},f;L^{p_{0}}(\Omega),\emph{BMO}(\Omega))}{t^{1/p_{0}%
}y_{\Theta}(\frac{1}{t})}\approx\sup_{t\in(0,\infty)}\frac{[(M_{\Omega}%
^{\#}f)^{p_{0}}]^{\ast\ast}(t)^{1/p_{0}}}{y_{\Theta}(\frac{1}{t})},
\label{CharSec2.1}%
\end{equation}
where $y_{\Theta}$ is given by \eqref{intro:Auxy}. 

\item[(ii)] Assume that the map
\begin{equation}
p\in(p_{0},\infty)\mapsto e^{p_{0}/p}\Theta(p)\qquad\text{is quasi-decreasing}. \label{HypTheta}%
\end{equation}
Then
\begin{align}
\Vert f\Vert_{Y_{p_{0}}^{\#\Theta}(\Omega)}  &  \approx\Vert f\Vert
_{(L^{p_{0}})^{\#}(\Omega)}+\sup_{t\in(0,e^{-1})}\frac{K(t,f;L^{p_{0}}%
(\Omega),\emph{BMO}(\Omega))}{t\,\Theta(-\log t)}\label{CharTheta1.1}\\
&  \approx\Vert f\Vert_{(L^{p_{0}})^{\#}(\Omega)}+\sup_{t\in(0,e^{-1})}%
\frac{[(M_{\Omega}^{\#}f)^{p_{0}}]^{\ast\ast}(t)^{1/p_{0}}}{\Theta(-\log
t)}.\nonumber
\end{align}
\end{enumerate}
\end{theorem}

In the proof of Theorem \ref{Thm2.2} we will make use of the following lemma,
which provides an explicit characterization of $y_{\Theta}$ in terms of
$\Theta$.

\begin{lemma}
\label{Lemma1} 
Let $y_\Theta$ be the Yudovich function relative to $\Theta$ (cf. \eqref{intro:Auxy}). 
\begin{enumerate}
\item[(i)] We have
\begin{equation}
y_{\Theta}(r)=\Theta(p_{0})r^{1/p_{0}},\qquad\text{if}\qquad r\in(0,1).
\label{Lemma1.1}%
\end{equation}
\item[(ii)] Suppose that $\Theta$ satisfies \eqref{HypTheta}. Then
\begin{equation*}
y_{\Theta}(r)\approx\Theta(\log r),\qquad\text{if}\qquad r>e^{p_{0}}.
\end{equation*}
\end{enumerate}
\end{lemma}

\begin{proof}
(i) The formula \eqref{Lemma1.1} is an immediate consequence of the fact that
for each $r\in(0,1),$ the map $p\in(p_{0},\infty)\mapsto\Theta(p)r^{1/p}$ is increasing.

(ii) Note that $y_{\Theta}(r)\leq\Theta(\log r)\,r^{1/\log r}\approx
\Theta(\log r)$. Conversely, given any $p>p_{0}$,
\[
\Theta(p)r^{1/p}\geq\Theta(\log r)\qquad\text{if}\qquad p>\log r.
\]
If $p\in(p_{0},\log r)$, then
\[
\Theta(p)r^{1/p}=e^{p_{0}/p}\Theta(p)\bigg(\frac{r}{e^{p_{0}}}\bigg)^{1/p}%
\gtrsim e^{p_{0}/\log r}\Theta(\log r)\bigg(\frac{r}{e^{p_{0}}}\bigg)^{1/\log
r}\approx\Theta(\log r).
\]
Hence,
\[
y_{\Theta}(r)=\inf_{p>p_{0}}\{\Theta(p)r^{1/p}\}\gtrsim\Theta(\log r).
\]

\end{proof}

\begin{proof}
[Proof of Theorem \ref{Thm2.2}](i) By the reiteration property of the $\Delta
$-extrapolation method (cf. \eqref{app1}), we have
\begin{equation}\label{212newnovel}
\Delta_{\theta\in(0,1)}\bigg\{\frac{(L^{p_{0}}(\Omega),\text{BMO}(\Omega))_{\theta,p(\theta
)}^{\blacktriangleleft}}{\Theta(\frac{p_{0}}{1-\theta})}\bigg\}=\Delta
_{\theta\in(0,1)}\bigg\{\frac{(L^{p_{0}}(\Omega),\text{BMO}(\Omega))_{\theta,\infty
}^{\blacktriangleleft}}{\Theta(\frac{p_{0}}{1-\theta})}\bigg\}.
\end{equation}
Then (cf. \eqref{la cara} and \eqref{4.5novel})
\begin{equation}\label{212newnovel2}
\Vert f\Vert_{\Delta_{\theta\in(0,1)}\Big\{\frac{(L^{p_{0}}(\Omega),\text{BMO}(\Omega)%
)_{\theta,p(\theta)}^{\blacktriangleleft}}{\Theta(\frac{p_{0}}{1-\theta}%
)}\Big\}}    \approx \sup_{t\in(0,\infty)}\frac{K(t^{1/p_{0}},f;L^{p_{0}}(\Omega),\text{BMO}(\Omega)%
)}{t^{1/p_{0}}y_{\Theta}(\frac{1}{t})}.
\end{equation}
Combining with \eqref{ClaimEx2} we obtain the first equivalence in
\eqref{CharSec2.1}. The second equivalence now follows from \eqref{JTFor}.

(ii) Suppose now that $\Theta$ satisfies \eqref{HypTheta}. By Lemma
\ref{Lemma1} and \eqref{JTFor}, we have
\begin{align}
\sup_{t\in(0,\infty)}\frac{K(t^{1/p_{0}},f;L^{p_{0}}(\Omega),\text{BMO}%
(\Omega))}{t^{1/p_{0}}y_{\Theta}(\frac{1}{t})}  &  \approx\nonumber\\
&  \hspace{-5.5cm}\sup_{t\in(0,e^{-p_{0}})}\frac{K(t^{1/p_{0}},f;L^{p_{0}%
}(\Omega),\text{BMO}(\Omega))}{t^{1/p_{0}}y_{\Theta}(\frac{1}{t})}+\sup
_{t\in(1,\infty)}\frac{K(t^{1/p_{0}},f;L^{p_{0}}(\Omega),\text{BMO}(\Omega
))}{t^{1/p_{0}}y_{\Theta}(\frac{1}{t})}\nonumber\\
&  \hspace{-5.5cm}\approx\sup_{t\in(0,e^{-p_{0}})}\frac{K(t^{1/p_{0}%
},f;L^{p_{0}}(\Omega),\text{BMO}(\Omega))}{t^{1/p_{0}}\Theta(-\log t)}%
+\sup_{t\in(1,\infty)}K(t^{1/p_{0}},f;L^{p_{0}}(\Omega),\text{BMO}%
(\Omega))\nonumber\\
&  \hspace{-5.5cm}\approx\sup_{t\in(0,e^{-1})}\frac{K(t,f;L^{p_{0}}%
(\Omega),\text{BMO}(\Omega))}{t\Theta(-\log t)}+\Vert M^{\#}_\Omega f\Vert_{L^{p_{0}%
}(\Omega)}. \label{212}%
\end{align}
Similarly, one can show that
\begin{equation}
\sup_{t\in(0,\infty)}\frac{[(M_{\Omega}^{\#}f)^{p_{0}}]^{\ast\ast}%
(t)^{1/p_{0}}}{y_{\Theta}(\frac{1}{t})}\approx\sup_{t\in(0,e^{-1})}%
\frac{[(M_{\Omega}^{\#}f)^{p_{0}}]^{\ast\ast}(t)^{1/p_{0}}}{\Theta(-\log
t)}+\Vert M^{\#}_\Omega f\Vert_{L^{p_{0}}(\Omega)}. \label{213}%
\end{equation}
Combining \eqref{CharSec2.1} with \eqref{212} and \eqref{213}, we complete the
proof of \eqref{CharTheta1.1}.
\end{proof}

When $\Omega=\mathbb{T}^{2}$, the characterizations provided by Theorem
\ref{Thm2.2} can be supplemented as follows.

\begin{theorem}
\label{Thm2.13} Let $p_{0}\in\lbrack1,\infty)$ and suppose that $\Theta$
satisfies \eqref{HypTheta}. Then
\begin{align}
\Vert f\Vert_{Y_{p_{0}}^{\#\Theta}(\mathbb{T}^{2})}  &  \approx\sup
_{t\in(0,e^{-1})}\frac{K(t,f;L^{p_{0}}(\mathbb{T}^{2}),\emph{BMO}%
(\mathbb{T}^{2}))}{t\Theta(-\log t)} \nonumber\\
&  \approx\sup_{t\in(0,e^{-1})}\frac{[(M_{\mathbb{T}^{2}}^{\#}f)^{p_{0}%
}]^{\ast\ast}(t)^{1/p_{0}}}{\Theta(-\log t)}\approx\sup_{t\in(0,e^{-1})}%
\frac{(M_{\mathbb{T}^{2}}^{\#}f)^{\ast}(t)}{\Theta(-\log t)}.\label{CharSec2.1new}
\end{align}

\end{theorem}

\begin{remark}
Note that the last expression in \eqref{CharSec2.1new} is independent of $p_0$. This is in accord with the fact that the definition of $Y_{p_{0}}^{\#\Theta}(\mathbb{T}^{2})$ does not depend on $p_0$. 
\end{remark}

\begin{proof}[Proof of Theorem \ref{Thm2.13}]
Using \eqref{JTFor} we can write
\[
K(e^{-1},f;L^{p_{0}}(\mathbb{T}^{2}),\text{BMO}(\mathbb{T}^{2}))\approx
\lbrack(M_{\mathbb{T}^{2}}^{\#}f)^{p_{0}}]^{\ast\ast}(e^{-1})\approx\Vert
f\Vert_{(L^{p_{0}})^{\#}(\mathbb{T}^{2})}.
\]
Accordingly, the first and second equivalences in \eqref{CharSec2.1new} are
consequences of \eqref{CharTheta1.1}.


We now prove the last equivalence in \eqref{CharSec2.1new}. We claim that
\begin{equation}
\int_{0}^{t}y_{\Theta}\bigg(\frac{1}{s}\bigg)^{p_{0}}\,ds\lesssim ty_{\Theta
}\bigg(\frac{1}{t}\bigg)^{p_{0}},\qquad t\in(0,1). \label{ClaimInitial}%
\end{equation}
Assuming momentarily the validity of (\ref{ClaimInitial}), we have
\begin{align*}
\lbrack(M_{\mathbb{T}^{2}}^{\#}f)^{p_{0}}]^{\ast\ast}(t)^{1/p_{0}}  &
=\bigg[\frac{1}{t}\int_{0}^{t}\bigg[\frac{(M_{\mathbb{T}^{2}}^{\#}f)^{\ast
}(s)}{y_{\Theta}(\frac{1}{s})} \, \, y_{\Theta}\bigg(\frac{1}%
{s}\bigg)\bigg]^{p_{0}}\,ds\bigg]^{1/p_{0}}\\
&  \leq\bigg[\frac{1}{t}\int_{0}^{t}y_{\Theta}\bigg(\frac{1}{s}\bigg)^{p_{0}%
}\,ds\bigg]^{1/p_{0}}\,\sup_{t\in(0,1)}\frac{(M_{\mathbb{T}^{2}}^{\#}f)^{\ast
}(t)}{y_{\Theta}(\frac{1}{t})}\\
&  \lesssim y_{\Theta}\bigg(\frac{1}{t}\bigg)\,\sup_{t\in(0,1)}\frac
{(M_{\mathbb{T}^{2}}^{\#}f)^{\ast}(t)}{y_{\Theta}(\frac{1}{t})}.
\end{align*}
Consequently,
\[
\sup_{t\in(0,1)}\frac{[(M_{\mathbb{T}^{2}}^{\#}f)^{p_{0}}]^{\ast\ast
}(t)^{1/p_{0}}}{y_{\Theta}(\frac{1}{t})}\lesssim\sup_{t\in(0,1)}%
\frac{(M_{\mathbb{T}^{2}}^{\#}f)^{\ast}(t)}{y_{\Theta}(\frac{1}{t})}.
\]
The converse estimate is a simple consequence of the fact that $(M_{\mathbb{T}%
^{2}}^{\#}f)^{\ast}$ is decreasing. This completes the proof of the last
equivalence in \eqref{CharSec2.1new}, under the assumption that
\eqref{ClaimInitial} holds true.

Next we turn to the proof of \eqref{ClaimInitial}. For this purpose we need
the following easily established fact:
\begin{equation}
y_{\Theta}(r)\approx\widetilde{y}_{\Theta}(r):=\inf_{p>2p_{0}}\{\Theta
(p)r^{1/p}\},\qquad\text{for}\qquad r>1. \label{Claim1}%
\end{equation}
Indeed, the estimate $y_{\theta}(r)\leq\widetilde{y}_{\Theta}(r)$ is obvious.
Conversely, let $p>p_{0},$ then $2p>2p_{0},$ and using the doubling property
of $\Theta$ we see that, for $r>1,$
\[
\Theta(p)r^{\frac{1}{p}}\approx\Theta(2p)r^{\frac{1}{p}}\geq\Theta
(2p)r^{\frac{1}{2p}}\geq\widetilde{y}_{\Theta}(r).
\]
Therefore, taking the infimum over all $p>p_{0}$, we obtain $y_{\Theta
}(r)\gtrsim\widetilde{y}_{\Theta}(r)$, completing the proof of \eqref{Claim1}.

Let $\varepsilon\in(\frac{1}{2p_{0}},\frac{1}{p_{0}})$. Observe that
\begin{equation}
r\in(1,\infty)\mapsto r^{-\varepsilon}\,\widetilde{y}_{\Theta}(r)=\inf
_{p>2p_{0}}\{\Theta(p)r^{-\varepsilon+1/p}\}\quad\text{is a decreasing
function}. \label{Claim2}%
\end{equation}
Let $t\in(0,1)$, then by \eqref{Claim1} and \eqref{Claim2}, we have
\begin{align*}
\int_{0}^{t}y_{\Theta}\bigg(\frac{1}{s}\bigg)^{p_{0}}\,ds  &  \approx\int%
_{0}^{t}\widetilde{y}_{\Theta}\bigg(\frac{1}{s}\bigg)^{p_{0}}\,ds=\int_{0}%
^{t}s^{-\varepsilon p_{0}}\bigg[\bigg(\frac{1}{s}\bigg)^{-\varepsilon
}\widetilde{y}_{\Theta}\bigg(\frac{1}{s}\bigg)\bigg]^{p_{0}}\,ds\\
&  \leq\bigg[\bigg(\frac{1}{t}\bigg)^{-\varepsilon}\widetilde{y}_{\Theta
}\bigg(\frac{1}{t}\bigg)\bigg]^{p_{0}}\int_{0}^{t}s^{-\varepsilon p_{0}}\,ds\\
&  \approx t^{1-\varepsilon p_{0}}\bigg[\bigg(\frac{1}{t}\bigg)^{-\varepsilon
}\widetilde{y}_{\Theta}\bigg(\frac{1}{t}\bigg)\bigg]^{p_{0}}\approx
ty_{\Theta}\bigg(\frac{1}{t}\bigg)^{p_{0}},
\end{align*}
this concludes the proof of (\ref{ClaimInitial}) and consequently the theorem
is proved.
\end{proof}

\subsection{Proof of Theorem \ref{Theorem2}\label{theorem2}}\label{Section22}

We shall assume that all function spaces are defined on $\Omega=\mathbb{R}%
^{2}$ (cf. Remark \ref{elcasoT} below for the modifications needed to deal
with $\Omega=\mathbb{T}^{2}$).

Let $p_{0}\in(1,\infty).$ Recall that CZOs act boundedly on $L^{p_{0}}\text{
and BMO}$ (cf. \cite{P66}), in particular, for $\mathcal{K}$ given by
\eqref{intro2} we have,
\[
\mathcal{K}:L^{p_{0}}\rightarrow L^{p_{0}}\qquad\text{and}\qquad
\mathcal{K}:\text{BMO}\rightarrow\text{BMO}.
\]
By interpolation we find that, with norm independent of $\theta$,
\begin{equation}\label{IntBMOInde}
\mathcal{K}:(L^{p_{0}},\text{BMO})_{\theta,p(\theta)}^{\blacktriangleleft
}\rightarrow(L^{p_{0}},\text{BMO})_{\theta,p(\theta)}^{\blacktriangleleft},
\end{equation}
where $p(\theta)=p_{0}/(1-\theta)$. This implies (cf. \eqref{intro2})
\begin{equation}
\Vert\nabla v\Vert_{(L^{p_{0}},\text{BMO})_{\theta,p(\theta)}%
^{\blacktriangleleft}}\lesssim\Vert\omega\Vert_{(L^{p_{0}},\text{BMO}%
)_{\theta,p(\theta)}^{\blacktriangleleft}} \label{JNEstimPr1.1}%
\end{equation}
with constant independent of $\theta$.

Next we show that
\begin{equation}
\Vert f-f_{\infty}\Vert_{(L^{p_{0}},L^{\infty})_{\theta,p(\theta
)}^{\blacktriangleleft}}\lesssim(1-\theta)^{-1}\Vert f\Vert_{(L^{p_{0}%
},\text{BMO})_{\theta,p(\theta)}^{\blacktriangleleft}}, \label{JNEstimPr1}%
\end{equation}
where $f_{\infty}:=\lim_{|Q|\rightarrow\infty}f_{Q}$ and $f_{Q}$ denotes the
\emph{integral average} of $f$ related to the cube $Q$, i.e.,  $f_{Q}:=\frac{1}{|Q|}\int_{Q}f$.

We start by reformulating the following quantitative version of the
John--Nirenberg embedding \cite{JN61}, which asserts that  BMO is locally
embedded into $e^{L}$ (cf. \cite[Proposition 8.10, p. 398]{BS88})
\begin{equation}
(f-f_{\infty})^{\ast\ast}(t)\lesssim\int_{t}^{\infty}f^{\#\ast}(s)\,\frac
{ds}{s}. \label{JNEstimPr}%
\end{equation}
Here $f^{\#}$ is the \emph{sharp maximal function of Fefferman--Stein}%
\footnote{Recall that $f^{\#}(x)=\sup_{Q\ni x}\frac{1}{|Q|}\int_{Q}|f-f_{Q}%
|$.} \cite{FS72}. Indeed, for our purposes, it is convenient to rewrite
(\ref{JNEstimPr}) in terms of $M^{\#}f$ rather than $f^{\#}.$ The connection
between these maximal functions is given by (cf. \cite[Lemma 3.4]{JT85})
\begin{equation}
f^{\#\ast}(t)\approx(M^{\#}f)^{\ast\ast}(t). \label{EstimMaxFunctions}%
\end{equation}
Therefore, \eqref{JNEstimPr} can be expressed as
\begin{equation}
(f-f_{\infty})^{\ast\ast}(t)\lesssim\int_{t}^{\infty}(M^{\#}f)^{\ast\ast
}(s)\,\frac{ds}{s}. \label{JNEstimPr2}%
\end{equation}
Applying $L^{p}$-norms on both sides of \eqref{JNEstimPr2} and estimating the
right-hand side using the pair of Hardy inequalities in \cite[Appendix A.4, page 272]{Stein},
we arrive at
\begin{align*}
\bigg\{\int_{0}^{\infty}[(f-f_{\infty})^{\ast\ast}(t)]^{p}\,dt\bigg\}^{1/p}
&  \lesssim\bigg\{\int_{0}^{\infty}\bigg[\int_{t}^{\infty}(M^{\#}f)^{\ast\ast
}(s)\,\frac{ds}{s}\bigg]^{p}\,dt\bigg\}^{1/p}\\
&  \leq p\,\bigg\{\int_{0}^{\infty}[(M^{\#}f)^{\ast\ast}(t)]^{p}%
\,dt\bigg\}^{1/p}\\
&  =p\,\bigg\{\int_{0}^{\infty}\bigg[\frac{1}{t}\int_{0}^{t}(M^{\#}f)^{\ast
}(s)\,ds\bigg]^{p}\,dt\bigg\}^{1/p}\\
&  \leq\frac{p^{2}}{p-1}\,\Vert f\Vert_{(L^{p})^{\#}}.
\end{align*}
Consequently, we have the following variant of the Fefferman--Stein inequality
(cf. \cite{BS88})
\begin{equation}\label{FSIneq}
\Vert f-f_{\infty}\Vert_{L^{p}}=\bigg\{\int_{0}^{\infty}[(f-f_{\infty})^{\ast
}(t)]^{p}\,dt\bigg\}^{1/p}\lesssim\frac{p^{2}}{p-1}\,\Vert f\Vert
_{(L^{p})^{\#}}.
\end{equation}
In particular, since $p_{0}>1,$ and $\frac{1}{p(\theta)}=\frac{1-\theta}%
{p_{0}},$
\[
\Vert f-f_{\infty}\Vert_{L^{p(\theta)}}\lesssim\frac{(1-\theta)^{-1}}%
{p_{0}-1+\theta}\,\Vert f\Vert_{(L^{p(\theta)})^{\#}}\lesssim(1-\theta
)^{-1}\Vert f\Vert_{(L^{p(\theta)})^{\#}}.
\]
In view of \eqref{IntroLp} and \eqref{ClaimEx1}, the previous estimate can be
rewritten as
\[
\Vert f-f_{\infty}\Vert_{(L^{p_{0}},L^{\infty})_{\theta,p(\theta
)}^{\blacktriangleleft}}\lesssim(1-\theta)^{-1}\Vert f\Vert_{(L^{p_{0}%
},\text{BMO})_{\theta,p(\theta)}^{\blacktriangleleft}},
\]
proving that \eqref{JNEstimPr1} holds.

Applying \eqref{JNEstimPr1} to $f=\nabla v$ (and noting that $(\nabla
v)_{\infty}=0$), combined with \eqref{JNEstimPr1.1}, yields
\[
\Vert\nabla v\Vert_{(L^{p_{0}},L^{\infty})_{\theta,p(\theta)}%
^{\blacktriangleleft}}\lesssim(1-\theta)^{-1}\Vert\omega\Vert_{(L^{p_{0}%
},\text{BMO})_{\theta,p(\theta)}^{\blacktriangleleft}}%
\]
uniformly with respect to $\theta\in(0,1)$. Multiplying both sides of the
above estimate by the factor $(1-\theta) \Theta\big(\frac{p_{0}}{1-\theta
}\big)^{-1}$ and taking the supremum over all $\theta\in(0,1)$, we
find
\begin{equation}
\Vert\nabla v\Vert_{\Delta_{\theta\in(0,1)}\Big\{\frac{(L^{p_{0}},L^{\infty
})_{\theta,p(\theta)}^{\blacktriangleleft}}{\Theta_{1}(\frac{p_{0}}{1-\theta
})}\Big\}}\lesssim\left\Vert \omega\right\Vert _{\Delta_{\theta\in
(0,1)}\Big\{\frac{(L^{p_{0}},\text{BMO})_{\theta,p(\theta)}%
^{\blacktriangleleft}}{\Theta(\frac{p_{0}}{1-\theta})}\Big\}}\approx
\Vert\omega\Vert_{Y_{p_{0}}^{\#\Theta}}, \label{JNEstimPr3}%
\end{equation}
where $\Theta_{1}$ was introduced in \eqref{intro:AuxTheta1} and we have also used \eqref{ClaimEx2} in the last step.

Now using \eqref{IntroIntSobDS}, we rewrite the left-hand side of
\eqref{JNEstimPr3} as%
\begin{equation*}
\Vert\nabla v\Vert_{\Delta_{\theta\in(0,1)}\Big\{\frac{(L^{p_{0}},L^{\infty
})_{\theta,p(\theta)}^{\blacktriangleleft}}{\Theta_{1}(\frac{p_{0}}{1-\theta
})}\Big\}}\approx\Vert v\Vert_{\Delta_{\theta\in(0,1)}\Big\{\frac{(\dot
{W}_{p_{0}}^{1},\dot{W}_{\infty}^{1})_{\theta,p(\theta)}^{\blacktriangleleft}%
}{\Theta_{1}(\frac{p_{0}}{1-\theta})}\Big\}}.
\end{equation*}
Furthermore, the Sobolev embedding theorem $\dot{W}_{p_{0}}^{1}\hookrightarrow
L^{\infty}$ (recall that $p_0 > 2$), and the reiteration property of $\Delta$-extrapolation (cf.
Appendix \ref{sec:A1}, (\ref{app1})) yield
\begin{align*}
\Delta_{\theta\in(0,1)}\bigg\{\frac{(\dot{W}_{p_{0}}^{1},\dot{W}_{\infty}%
^{1})_{\theta,p(\theta)}^{\blacktriangleleft}}{\Theta_{1}(\frac{p_{0}%
}{1-\theta})}\bigg\}  &  \hookrightarrow\Delta_{\theta\in(0,1)}\bigg\{\frac
{(L^{\infty},\dot{W}_{\infty}^{1})_{\theta,p(\theta)}^{\blacktriangleleft}%
}{\Theta_{1}(\frac{p_{0}}{1-\theta})}\bigg\}\\
&  =\Delta_{\theta\in(0,1)}\bigg\{\frac{(L^{\infty},\dot{W}_{\infty}%
^{1})_{\theta,\infty}^{\blacktriangleleft}}{\Theta_{1}(\frac{p_{0}}{1-\theta
})}\bigg\}.
\end{align*}
Consequently,
\begin{equation}
\left\Vert v\right\Vert _{\Delta_{\theta\in(0,1)}\Big\{\frac{(L^{\infty}%
,\dot{W}_{\infty}^{1})_{\theta,\infty}^{\blacktriangleleft}}{\Theta_{1}%
(\frac{p_{0}}{1-\theta})}\Big\}}\lesssim\Vert v\Vert_{\Delta_{\theta\in
(0,1)}\Big\{\frac{(\dot{W}_{p_{0}}^{1},\dot{W}_{\infty}^{1})_{\theta
,p(\theta)}^{\blacktriangleleft}}{\Theta_{1}(\frac{p_{0}}{1-\theta})}\Big\}}.
\label{JNEstimPr4}%
\end{equation}

Updating (\ref{JNEstimPr3}) via \eqref{JNEstimPr4} we arrive at
\begin{equation}
\left\Vert v\right\Vert _{\Delta_{\theta\in(0,1)}\Big\{\frac{(L^{\infty}%
,\dot{W}_{\infty}^{1})_{\theta,\infty}^{\blacktriangleleft}}{\Theta_{1}%
(\frac{p_{0}}{1-\theta})}\Big\}}\lesssim\Vert\omega\Vert_{Y_{p_{0}}^{\#\Theta
}}. \label{dersh}%
\end{equation}
Next we compute the norm of the extrapolation space on the left-hand side of
\eqref{dersh} and show that%
\begin{equation}
\Vert v\Vert_{\Delta_{\theta\in(0,1)}\Big\{\frac{(L^{\infty},\dot{W}_{\infty
}^{1})_{\theta,\infty}^{\blacktriangleleft}}{\Theta_{1}(\frac{p_{0}}{1-\theta
})}\Big\}}\approx\sup_{x,y }\,\frac{|v(x)-v(y)|}{\inf_{t>|x-y|}%
ty_{\Theta_{1}}(\frac{1}{t})}. \label{225new}%
\end{equation}
The argument was already outlined in Section
\ref{SectionIntroAPrior}. For the sake of completeness,  next we 
give full details. Indeed, in light of \eqref{intro:Auxy}, we have (cf. \eqref{la cara} and \eqref{4.5novel})
\begin{equation}
\left\Vert v\right\Vert _{\Delta_{\theta\in(0,1)}\Big\{\frac{(L^{\infty}%
,\dot{W}_{\infty}^{1})_{\theta,\infty}^{\blacktriangleleft}}{\Theta_{1}%
(\frac{p_{0}}{1-\theta})}\Big\}}     \approx\sup_{t\in(0,\infty)}\frac{K(t,v;L^{\infty},\dot{W}_{\infty}^{1}%
)}{ty_{\Theta_{1}}(\frac{1}{t})}. \label{alga}%
\end{equation}
Applying the characterization of the $K$-functional for $(L^{\infty},\dot
{W}_{\infty}^{1})$ (cf. \eqref{IntroKFunctModLinfty}) we find
\begin{align*}
\left\Vert v\right\Vert _{\Delta_{\theta\in(0,1)}\Big\{\frac{(L^{\infty}%
,\dot{W}_{\infty}^{1})_{\theta,\infty}^{\blacktriangleleft}}{\Theta_{1}%
(\frac{p_{0}}{1-\theta})}\Big\}}  &  \approx\sup_{t\in(0,\infty)}\frac{\sup_{|x-y|<t}|v(x)-v(y)|}%
{ty_{\Theta_{1}}(\frac{1}{t})}\nonumber\\
&  =\sup_{x,y }|v(x)-v(y)|\,\sup_{t>|x-y|}\frac{1}{ty_{\Theta_{1}%
}(\frac{1}{t})}. 
\end{align*}
Hence, the desired result \eqref{225new} follows.

Putting together \eqref{dersh} and \eqref{225new}, we
obtain
\[
\sup_{x,y\in\Omega}\,\frac{|v(x)-v(y)|}{\inf_{t>|x-y|}ty_{\Theta_{1}}(\frac
{1}{t})}\lesssim\Vert\omega\Vert_{Y_{p_{0}}^{\#\Theta}}.
\]
 In particular,
\[
|v(x)-v(y)|\lesssim|x-y|\,y_{\Theta_{1}}\bigg(\frac{1}{|x-y|}\bigg)\Vert
\omega\Vert_{Y_{p_{0}}^{\#\Theta}},
\]
thus completing the proof of Theorem \ref{Theorem2}. \qed

\begin{remark}
\label{elcasoT}The above proof can be easily adapted to deal with
$\Omega=\mathbb{T}^{2}$. In particular, the periodic counterpart of \eqref{JNEstimPr}
is given by (cf. \cite[Corollary 7.4, p. 379]{BS88}), for $t\in(0,\frac{1}%
{6})$,
\[
(f-f_{\mathbb{T}^{2}})^{\ast\ast}(t)\lesssim\int_{t}^{1}f^{\#\ast}%
(s)\,\frac{ds}{s},
\]
where $f_{\mathbb{T}^{2}}$ denotes the integral mean of $f$. Accordingly, the
fact that $(\nabla v)_{\infty}=0$ in the above proof is replaced
by\footnote{Using the Fundamental Theorem of Calculus together with the
continuity of $v$ (recall that $\omega\in L^{p}(\mathbb{T}^{2})$ with $p>2$),
we have $(\nabla v)_{\mathbb{T}^{2}}=0$.} $(\nabla v)_{\mathbb{T}^{2}}=0$.
\end{remark}

\subsection{Proof of Theorem \ref{Thm3Intro} }\label{SectionUY}

 Using the extrapolation techniques developed above, together with the Fefferman-Stein inequality, we present two slightly different proofs of Theorem \ref{Thm3Intro}. Our first proof is close in  spirit to the classical energy method of \cite{Y63}, involving $E(t)$ (cf. \eqref{EnergyDef}), while the  second proof may be viewed as a sharpening of the variant energy method of \cite{AB} involving $\dot{H}^{-1}$ norms of vorticities. The advantage of the second approach lies in its additional flexibility. Indeed, as we will show in detail in Section \ref{SectionUniqAS} below, this method  can be applied to derive uniqueness results not only for the $2$D Euler equations, but also a wide variety of  active scalar equations.

\vspace{2mm}

\textbf{A first approach to Theorem \ref{Thm3Intro}.} Assume that $v_1$ and $v_2$ are two solutions of  \eqref{intro1} with the same initial data $v_1(0) = v_2(0)$. Recall that (cf. \eqref{EnergyEstim})
\begin{equation}\label{13.0}
	\frac{d  E(t) }{dt}\leq \int_\Omega |\nabla v_1(t, x)| |v_1(t, x)-v_2(t, x)|^2 \, dx,
\end{equation}
where $E(t)$ is given by \eqref{EnergyDef}. A simple application of H\"older's inequality gives 
\begin{equation}\label{13.1}
	\int_\Omega |\nabla v_1(t, x)| |v_1(t, x)-v_2(t, x)|^2 \, dx \leq \|\nabla v_1(t)\|_{L^p(\Omega)} (2E(t))^{1-\frac{1}{p}}  \||v_1(t)-v_2(t)|^{\frac{2}{p}}\|_{L^\infty (\Omega)}. 
\end{equation}
From the classical Sobolev embedding, \eqref{intro2}, and the definition of $Y^{\# \Theta}_{p_0}$,   we derive that $v_1$ and $v_2$ are bounded. On the other hand, by virtue of the Fefferman-Stein inequality \eqref{FSIneq} and the extrapolation assertion \eqref{ClaimEx1} contained in Theorem \ref{ThmExtSh}, we can estimate
$$
	\|\nabla v_1(t)\|_{L^p(\Omega)} \lesssim  p \|\mathcal{K} \omega_1(t)\|_{(L^p(\Omega))^\#} \approx p \|\mathcal{K} \omega_1(t)\|_{(L^{p_0}(\Omega), \text{BMO}(\Omega))_{\theta, p}^{\blacktriangleleft}},
$$
where $p = p_0/(1-\theta)$. Moreover, by \eqref{IntBMOInde}, we get
\begin{equation}\label{17.3}
	\|\nabla v_1(t)\|_{L^p(\Omega)} \lesssim  p \|\omega_1(t)\|_{(L^{p_0}(\Omega), \text{BMO}(\Omega))_{\theta, p}^{\blacktriangleleft}}.
\end{equation}
 Hence \eqref{13.1} yields that
\begin{align*}
	\int_\Omega |\nabla v_1(t, x)| |v_1(t, x)-v_2(t, x)|^2 \, dx &\lesssim p \|\omega_1(t)\|_{(L^{p_0}(\Omega), \text{BMO}(\Omega))_{\theta, p}^{\blacktriangleleft}} E(t)^{1-\frac{1}{p}} \\
	& \lesssim  \Theta_1(p) \|\omega_1(t)\|_{Y^{\#\Theta}_{p_0}(\Omega)} E(t)^{1-\frac{1}{p}},
\end{align*}
where in the last step we have used \eqref{ClaimEx2} (recall also \eqref{intro:AuxTheta1}). Taking now the infimum over all $p > p_0$ in the last expression and using \eqref{13.0}, we arrive at 
\begin{equation}\label{18.1}
		\frac{d  E(t)}{dt} \lesssim \|\omega_1(t)\|_{Y^{\#\Theta}_{p_0}(\Omega)} E(t)   y_{\Theta_1} (E(t)^{-1}),
\end{equation}
i.e., the classical Yudovich estimate \eqref{KeyYudo} still holds for the larger class $Y^{\#\Theta}_{p_0}(\Omega)$. As a byproduct, the desired uniqueness assertion for $Y^{\#\Theta}_{p_0}(\Omega)$ follows from the Osgood condition on $y_{\Theta_1}$ (cf. \eqref{intro:osg}). \qed

\vspace{2mm}
\textbf{A second approach to Theorem \ref{Thm3Intro}.}
	Assume that $\omega_1$ and $\omega_2$ are two weak solutions of  \eqref{intro1} with $\omega_1(0)= \omega_2(0)$ and corresponding velocities $v_1$ and $v_2$. Consider the difference $\rho = \omega_1-\omega_2$ whose corresponding stream function $\psi$ is given by $\Delta \psi = \rho$ and velocity $v = k \ast \rho$. Then
	\begin{equation}\label{17.8} 
	\frac{1}{2} \frac{d}{dt} \int_\Omega |\nabla \psi|^2 \, dx = \int_\Omega \nabla \psi \cdot \rho v_1 \, dx + \int_\Omega \nabla \psi \cdot \omega_2 v \, dx.
	\end{equation}

It is convenient to express the first integral on the right-hand side of \eqref{17.8} as follows
\begin{align}
	 \int_\Omega \nabla \psi \cdot \rho v_1 \, dx &= \sum_{i, l} \int_\Omega  \partial_{ii} \psi  \partial_{l} \psi (v_1)_l \, dx \nonumber \\
	 & = \sum_{i, l} \int_\Omega [\partial_i \psi \partial_{li} \psi (v_1)_l + \partial_i \psi \partial_l \psi \partial_i (v_1)_l] \, dx \nonumber \\
	 & = \sum_{i, l}	\int_\Omega \bigg[ \partial_l \Big(\frac{1}{2} |\partial_i \psi|^2 \Big) (v_1)_l  + \partial_i \psi \partial_l \psi \partial_i (v_1)_l \bigg] \, dx.  \label{17.1}
\end{align}
In particular, applying H\"older's inequality for $p > 2$ (compare with \eqref{13.1})
\begin{equation*}
	 \int_\Omega \nabla \psi \cdot \rho v_1 \, dx \lesssim \int_\Omega |\nabla v_1| |\nabla \psi|^2 \, dx  \leq \|\nabla v_1\|_{L^p(\Omega)}  \|\nabla \psi\|_{L^2(\Omega)}^{2 (1-\frac{1}{p})}  \||\nabla \psi|^{\frac{2}{p}}\|_{L^\infty(\Omega)}. 
\end{equation*}
We estimate the right hand side invoking  \eqref{17.3} and Theorem \ref{ThmExtSh}, together with the fact that $\nabla \psi \in L^\infty$ (apply classical Sobolev inequality to $\psi \in \dot{W}^{2}_p$),
\begin{align}
	 \int_\Omega \nabla \psi \cdot \rho v_1 \, dx &\lesssim  p \|\omega_1\|_{(L^{p_0}(\Omega), \text{BMO}(\Omega))_{\theta, p}^{\blacktriangleleft}}  \|\nabla \psi\|_{L^2(\Omega)}^{2 (1-\frac{1}{p})} \nonumber \\
	 & \lesssim \Theta_1(p) \|\omega_1\|_{Y^{\# \Theta}_{p_0}(\Omega)} \|\nabla \psi\|_{L^2(\Omega)}^{2 (1-\frac{1}{p})}. \label{17.4}
\end{align}

On the other hand, applying H\"older's inequality twice and the Fefferman-Stein inequality \eqref{FSIneq} to $\omega_2$ (taking into account that $\int_\Omega \omega_2 =0$ if $\Omega= \T^2$), yields
\begin{align}
 \int_\Omega \nabla \psi \cdot \omega_2 v \, dx &\leq \|\omega_2\|_{L^p(\Omega)} \|\nabla \psi\|_{L^{\frac{2p}{p-1}}(\Omega)} \|v\|_{L^{\frac{2 p}{p-1}}(\Omega)} \nonumber \\
 & \lesssim p \|\omega_2\|_{(L^p(\Omega))^{\#}} \|\nabla \psi\|_{L^\infty(\Omega)}^{\frac{1}{p}} \|\nabla \psi\|_{L^2(\Omega)}^{1-\frac{1}{p}} \|v\|_{L^\infty(\Omega)}^{\frac{1}{p}} \|v\|_{L^2(\Omega)}^{1-\frac{1}{p}} \nonumber \\
 & \lesssim \Theta_1(p) \|\omega_2\|_{Y^{\# \Theta}_{p_0}(\Omega)} \|\nabla \psi\|_{L^2(\Omega)}^{1-\frac{1}{p}} \|v\|_{L^2(\Omega)}^{1-\frac{1}{p}} \nonumber \\
 & \approx  \Theta_1(p) \|\omega_2\|_{Y^{\# \Theta}_{p_0}(\Omega)} \|\nabla \psi\|_{L^2(\Omega)}^{2(1-\frac{1}{p})},  \label{17.5}
\end{align}
where the last estimate follows from the fact that $\|v\|_{L^2(\Omega)} \approx \|\nabla \psi\|_{L^2(\Omega)}$ (cf. \eqref{intro2}). 

From \eqref{17.8},  \eqref{17.4} and \eqref{17.5}, we achieve
$$
	\frac{d}{dt} \|\nabla \psi(t)\|_{L^2(\Omega)}^2 \lesssim \|\nabla \psi(t)\|_{L^2(\Omega)}^2 \Theta_1(p) \bigg(\frac{1}{\|\nabla \psi(t)\|_{L^2(\Omega)}^2} \bigg)^{\frac{1}{p}}. 
$$
Taking the infimum over all $p > p_0$ on the right-hand side of the last estimate yields  (compare with \eqref{18.1})
$$
	\frac{d}{dt} \|\nabla \psi(t)\|_{L^2(\Omega)}^2 \lesssim  \|\nabla \psi(t)\|_{L^2(\Omega)}^2 y_{\Theta_1} \bigg(\frac{1}{\|\nabla \psi(t)\|_{L^2(\Omega)}^2} \bigg).
$$
We are now in a position to apply Osgood's lemma (cf. \eqref{intro:osg}) and derive that $\|\nabla \psi(t)\|_{L^2(\Omega)}=0$ for all $t \in [0, T]$ and hence $\rho(t) \equiv 0$.  \qed

\subsection{Examples of vorticities in $Y^{\# \Theta}_{p_{0}}$ that are not in
$Y^{\Theta}_{p_{0}}$}\label{SectionCounterExamples}

To simplify the exposition, throughout this section we assume that
$\Omega=\mathbb{T}^{2}$. Recall that by construction, given a growth function
$\Theta$, $Y_{p_{0}}^{\#\Theta}$ is a bigger space than $Y_{p_{0}}^{\Theta}$
(cf. \eqref{IntroEmYshY}). Furthermore, in the special case $\Theta
(p)\approx1$, we have $Y_{p_{0}}^{\Theta}=L^{\infty}\subsetneq\text{BMO}%
=Y_{p_{0}}^{\#\Theta}$. It is of interest to understand better the
relationship between these spaces. In this section, we provide a method to
construct explicit examples of functions $\omega\in Y_{p_{0}}^{\#\Theta
}\backslash Y_{p_{0}}^{\Theta},$ for a variety of growths.

\begin{example}
[The case $\Theta(p) \approx p^{\alpha}, \, \alpha> 0$]\label{ExampleLinear}
Let
\[
\omega(x) = |\log|x||^{\alpha+ 1}.
\]
We will show that
\begin{equation}
\label{Goal1}\omega\in Y^{\#\Theta}, \qquad\omega\not \in Y^{\Theta}.
\end{equation}
Indeed, basic computations lead to
\begin{equation}
\label{Counter1}\omega^{*}(t) \approx(-\log t)^{\alpha+ 1}%
\end{equation}
and (recall that $\Omega= \mathbb{T}^{2}$)
\begin{equation}
\label{Counter1.1}|\nabla\omega|^{*}(t) \approx t^{-1/2} (-\log t)^{\alpha}.
\end{equation}

Next we are able to bound $\omega^{\#\ast}$ via the following pointwise
estimate obtained in \cite[Theorem 2.6]{DT22}:
\begin{equation}
\omega^{\#\ast}(t)\lesssim\sum_{l=0}^{1}\bigg[t^{-1/r_{0}}\bigg(\int_{0}%
^{t}(\xi^{1/r}\,|\nabla^{l}\omega|^{\ast}(\xi))^{r_{0}}\,\frac{d\xi}{\xi
}\bigg)^{1/r_{0}}+\sup_{t<\xi<1}\xi^{1/2}\,|\nabla^{l}\omega|^{\ast}%
(\xi)\bigg], \label{EstimSGrad}%
\end{equation}
where $r_{0}>2$ and $r=\frac{2r_{0}}{2+r_{0}}$. We treat the four terms
appearing on the right-hand side of the previous estimate. Assume first that
$l=0$ (i.e., $|\nabla^{0}\omega|^{\ast}=\omega^{\ast}$). By \eqref{Counter1},
we have
\[
\bigg(\int_{0}^{t}(\xi^{1/r}\,\omega^{\ast}(\xi))^{r_{0}}\,\frac{d\xi}{\xi
}\bigg)^{1/r_{0}}\approx\bigg(\int_{0}^{t}[\xi^{1/r}(-\log\xi)^{\alpha
+1}]^{r_{0}}\,\frac{d\xi}{\xi}\bigg)^{1/r_{0}}\approx t^{1/r}(-\log
t)^{\alpha+1}%
\]
and
\[
\sup_{t<\xi<1}\xi^{1/2}\,\omega^{\ast}(\xi)\approx\sup_{t<\xi<1}\xi
^{1/2}\,(-\log\xi)^{\alpha+1}\approx1.
\]
Putting these estimates together, we obtain
\begin{align}
t^{-1/r_{0}}\bigg(\int_{0}^{t}(\xi^{1/r}\,\omega^{\ast}(\xi))^{r_{0}}%
\,\frac{d\xi}{\xi}\bigg)^{1/r_{0}}+\sup_{t<\xi<1}\xi^{1/2}\,\omega^{\ast}%
(\xi)  & \approx t^{1/2}(-\log t)^{\alpha+1}+1\approx1.\label{Counter2}
\end{align}

Next we deal with the term on the right-hand side of \eqref{EstimSGrad} that
corresponds to $l=1$. Using \eqref{Counter1.1}, we get the following
estimates
\begin{align*}
\bigg(\int_{0}^{t}(\xi^{1/r}\,|\nabla\omega|^{\ast}(\xi))^{r_{0}}\,\frac{d\xi
}{\xi}\bigg)^{1/r_{0}}  &  \approx\bigg(\int_{0}^{t}(\xi^{1/r_{0}}(-\log
\xi)^{\alpha})^{r_{0}}\,\frac{d\xi}{\xi}\bigg)^{1/r_{0}}\\
&  \approx t^{1/r_{0}}(-\log t)^{\alpha}%
\end{align*}
and (since $\alpha>0$)
\[
\sup_{t<\xi<1}\xi^{1/2}\,|\nabla\omega|^{\ast}(\xi)\approx\sup_{t<\xi<1}%
(-\log\xi)^{\alpha}=(-\log t)^{\alpha}.
\]
Hence
\begin{equation}
t^{-1/r_{0}}\bigg(\int_{0}^{t}(\xi^{1/r}\,|\nabla\omega|^{\ast}(\xi))^{r_{0}%
}\,\frac{d\xi}{\xi}\bigg)^{1/r_{0}}+\sup_{t<\xi<1}\xi^{1/2}\,|\nabla
\omega|^{\ast}(\xi)\approx(-\log t)^{\alpha}. \label{Counter3}%
\end{equation}

Inserting \eqref{Counter2} and \eqref{Counter3} into \eqref{EstimSGrad}, we
obtain the following upper estimate for $\omega^{\#*}$,
\[
\omega^{\#*}(t) \lesssim(-\log t)^{\alpha},
\]
or equivalently (cf. \eqref{EstimMaxFunctions})
\begin{equation}
\label{Counter4}(M^{\#} \omega)^{**}(t) \lesssim(-\log t)^{\alpha}.
\end{equation}

Since we are working on $\Omega=\mathbb{T}^{2},$ without loss of generality we
may assume that $p_{0}=1$. Applying Theorem \ref{Thm2.13} and Example \ref{ExIntroNew} (with $\Theta
(p)\approx p^{\alpha}$) and using \eqref{Counter1} and \eqref{Counter4}, we
compute
\[
\Vert\omega\Vert_{Y^{\Theta}}\approx\sup_{t\in(0,e^{-1})}\frac{\omega
^{\ast\ast}(t)}{\Theta(-\log t)}\approx\sup_{t\in(0,e^{-1})}\frac{(-\log
t)^{\alpha+1}}{(-\log t)^{\alpha}}=\infty
\]
and
\[
\Vert\omega\Vert_{Y^{\#\Theta}}\approx\sup_{t\in(0,e^{-1})}\frac{(M^{\#}%
\omega)^{\ast\ast}(t)}{\Theta(-\log t)}\lesssim\sup_{t\in(0,e^{-1})}%
\frac{(-\log t)^{\alpha}}{(-\log t)^{\alpha}}=1.
\]
This concludes the proof of \eqref{Goal1}.
\end{example}

\begin{example}
[The case $\Theta(p)\approx\log p$]\label{Example4} This example is motivated
by the fact that $\log(1+|\log|x||)$ is a prototype of a function in
$Y^{\Theta},$ while $|\log|x||$ is a prototype of a function in $\text{BMO}$.
We consider their product, namely,
\[
\omega(x)=(1+|\log|x||)\log(1+|\log|x||),
\]
and show that
\begin{equation}
\omega\in Y^{\#\Theta},\qquad\omega\not \in Y^{\Theta}. \label{Goal2}%
\end{equation}
We follow closely the method of Example \ref{ExampleLinear}. Specifically, by
elementary manipulations we find
\begin{equation}
\omega^{\ast}(t)\approx(1-\log t)\log(1-\log t) \label{ExExtreme}%
\end{equation}
and
\[
|\nabla\omega|^{\ast}(t)\approx t^{-1/2}\,(1+\log(1-\log t)).
\]
Recall that $r=\frac{2r_{0}}{2+r_{0}}$, where $r_{0}>2$. Therefore
\[
t^{-1/r_{0}}\bigg(\int_{0}^{t}(\xi^{1/r}\,\omega^{\ast}(\xi))^{r_{0}}%
\,\frac{d\xi}{\xi}\bigg)^{1/r_{0}}+\sup_{t<\xi<1}\xi^{1/2}\,\omega^{\ast}%
(\xi)\approx t^{1/2}(1-\log t)\log(1-\log t)+1\approx1
\]
and
\[
t^{-1/r_{0}}\bigg(\int_{0}^{t}(\xi^{1/r}\,|\nabla\omega|^{\ast}(\xi))^{r_{0}%
}\,\frac{d\xi}{\xi}\bigg)^{1/r_{0}}+\sup_{t<\xi<1}\xi^{1/2}\,|\nabla
\omega|^{\ast}(\xi)\approx\log(1-\log t).
\]
Inserting these two estimates into \eqref{EstimSGrad}, we obtain
\begin{equation}
\omega^{\#\ast}(t)\lesssim\log(1-\log t). \label{ExExtreme2}%
\end{equation}

Invoking Theorem \ref{Thm2.13} and Example \ref{ExIntroNew} (with $\Theta(p) \approx\log p$) together with
\eqref{ExExtreme} and \eqref{ExExtreme2}, we get
\[
\|\omega\|_{Y^{\Theta}} \approx\sup_{t \in(0, e^{-1})} (-\log t) = \infty
\]
and
\[
\|\omega\|_{Y^{\#\Theta}} \lesssim\sup_{t \in(0, e^{-1})} \frac{\log(-\log
t)}{\log(-\log t)} = 1.
\]
Hence $\omega$ fulfils \eqref{Goal2}.
\end{example}

\begin{remark}
	It is possible to extend the methodology applied in Examples \ref{ExampleLinear} and \ref{Example4}  to deal with more general growths $\Theta$ of logarithmic type. Further details are left to the reader. 
\end{remark}

\section{The spaces $\dot{B}^{\beta}_{\Pi}$}

\label{SectionVishik}

Let $\mathcal{S}(\mathbb{R}^{d})$ denote the Schwartz space and $\mathcal{S}%
^{\prime}(\mathbb{R}^{d})$ the space of tempered distributions. We consider
the space $\dot{\mathcal{S}}(\mathbb{R}^{d})$ formed by $\varphi \in \mathcal{S}(\R^d)$ with $(D^\alpha \widehat{\varphi} )(0) = 0$ for any multi-index $\alpha \in \N_0^d$, where $\N_0 = \N \cup \{0\}$; this space carries the
natural Fr\'{e}chet topology inherited from $\mathcal{S}(\mathbb{R}^{d})$. Let
$\dot{\mathcal{S}}^{\prime}(\mathbb{R}^{d})$ be its dual space, which can be
identified with $\mathcal{S}^{\prime}(\mathbb{R}^{d})$ modulo polynomials.

Let $\varphi\in C_{c}^{\infty}(\mathbb{R}^{d})$ be a radial function with
\[
\text{supp }\varphi\subset\bigg\{\xi:\frac{3}{4}<|\xi|<\frac{7}{4}\bigg\},
\]%
\[
\varphi(\xi)=1\quad\text{for}\quad\frac{7}{8}<|\xi|<\frac{9}{8},
\]%
\[
\sum_{j\in\mathbb{Z}}\varphi(2^{-j}\xi)=1\quad\text{for all}\quad\xi\neq0.
\]
Then
\[
\widehat{\dot{\Delta}_{j}f}(\xi):=\varphi(2^{-j}\xi)\widehat{f}(\xi
),
\]
and%
\[\Delta_{j}:=\dot{\Delta}_{j} \quad \text{if} \quad j>0, \qquad 
\Delta_{0}:=\text{Id}-\sum_{j>0}\Delta_{j}.
\]

\begin{definition}
\label{DefVishik} Let $\beta\in\mathbb{R}$, and let $\Pi$ be a growth
function. Thus $\dot{B}_{\Pi}^{\beta}(\mathbb{R}^{d})$ will denote the
(homogeneous) \emph{Vishik space} formed by all $f\in\dot{\mathcal{S}}%
^{\prime}(\mathbb{R}^{d})$ such that
\[
\Vert f\Vert_{\dot{B}_{\Pi}^{\beta}(\mathbb{R}^{d})}:=\sup_{N\geq0}\,\frac
{1}{\Pi(N)}\,\sum_{j=-\infty}^{N}2^{j\beta}\,\Vert\dot{\Delta}_{j}%
f\Vert_{L^{\infty}(\mathbb{R}^{d})}<\infty.
\]
The inhomogeneous counterpart, $B_{\Pi}^{\beta}(\mathbb{R}^{d}),$ is the set
of all $f\in\mathcal{S}^{\prime}(\mathbb{R}^{d})$ such that
\[
\Vert f\Vert_{B_{\Pi}^{\beta}(\mathbb{R}^{d})}:=\sup_{N\geq0}\,\frac{1}%
{\Pi(N)}\,\sum_{j=0}^{N}2^{j\beta}\,\Vert\Delta_{j}f\Vert_{L^{\infty
}(\mathbb{R}^{d})}<\infty.
\]
The corresponding periodic spaces, $\dot{B}_{\Pi}^{\beta}(\mathbb{T}^{d})$ and
$B_{\Pi}^{\beta}(\mathbb{T}^{d}),$ can be defined analogously.
\end{definition}

\begin{remark}
\label{RemVishik}

\begin{enumerate}
\item[(i)] Standard properties of multipliers can be used to show that 
Vishik spaces do not depend (up to equivalence of norms) on the chosen
generator $\varphi$.

\item[(ii)] In the special case $\beta=0,$ $B_{\Pi}^{\beta}(\mathbb{R}^{d})$
coincides with the classical space $B_{\Pi}(\mathbb{R}^{d})$ (cf. \eqref{introVishikSpace}).

\item[(iii)] Assume $\Pi(N) \approx1$. Then $\dot{B}^{\beta}_{\Pi}(\mathbb{R}^{d})
= \dot{B}^{\beta}_{\infty, 1}(\mathbb{R}^{d})$, the classical Besov space  defined by
\begin{equation}
\label{BesovDef}\|f\|_{ \dot{B}^{\beta}_{\infty, 1}(\mathbb{R}^{d})} :=
\sum_{j \in\mathbb{Z}} 2^{j \beta} \, \|\dot{\Delta}_{j} f\|_{L^{\infty
}(\mathbb{R}^{d})}.
\end{equation}
Analogously,  $B^{\beta}_{\Pi}(\mathbb{R}^{d}) = B^\beta_{\infty, 1}(\R^d)$ with
$$
	\|f\|_{B^{\beta}_{\infty, 1}(\mathbb{R}^{d})} :=
\sum_{j \in\mathbb{N}_0} 2^{j \beta} \, \|\Delta_{j} f\|_{L^{\infty
}(\mathbb{R}^{d})}.
$$
The reader is referred to Appendix \ref{sec:retract} for further details on classical Besov spaces. 
\end{enumerate}
\end{remark}

\subsection{Characterizations}

\label{SectionCharacVishik}

The goal of this section is to provide characterizations of the Vishik spaces
as extrapolation spaces. We show analogs of the results obtained earlier for
the $Y_{p_{0}}^{\Theta}$ and $Y_{p_{0}}^{\#\Theta}$ spaces (cf. Section
\ref{Section2.1}).

We need to impose some natural restrictions on the growth functions used in
this section.

\begin{definition}
\label{DefPClass} Let $\kappa>0.$ We shall denote by $\mathcal{P}_{\kappa}$,
the set of all growth functions $\Pi$ satisfying the following conditions:

\begin{enumerate}
\item[(i)] $\Pi: [0, \infty) \to(0, \infty)$ is non-decreasing,

\item[(ii)] $\Pi$ is doubling,

\item[(iii)] $e^{1/p} \,\Pi(p)$ is quasi-decreasing,

\item[(iv)] $\sum_{j=N}^{\infty}2^{-j \kappa} \Pi(j) \lesssim2^{-N \kappa}
\Pi(N)$ for every $N \geq0$.
\end{enumerate}
\end{definition}

Clearly,  $\Pi(p)=(p+1)^{\alpha} (\log (p+e))^{\alpha_1} (\log_2 (p+e))^{\alpha_2} \cdots (\log_m (p+e))^{\alpha_m}$, where $\alpha,\alpha_i 
\geq0$, are examples of growth functions in $\mathcal{P}_{\kappa}.$

In dealing with interpolation of Besov spaces, we make use of the fact that the
computation of $K$-functionals for the Besov pair $(\dot{B}_{\infty
,1}^{\beta-\kappa}(\mathbb{R}^{d}),\dot{B}_{\infty,1}^{\beta}(\mathbb{R}^{d})$
can be reduced, via the method of retracts (cf. Appendix \ref{sec:retract} for further details),
to the computation of $K$-functionals for the vector-valued sequence spaces
$(\ell_{1}^{\beta-\kappa}(L^{\infty}(\mathbb{R}^{d}),\ell_{1}^{\beta}(L^{\infty
}(\mathbb{R}^{d}))$ (cf. \eqref{DefSeqSpaces}). 


Using the retract technique we proved Theorem \ref{teo:nec} below which, in
particular, implies%
\begin{equation}
(  \dot{B}_{\infty,1}^{\beta-\kappa}(\mathbb{R}^{d}),\dot{B}_{\infty
,1}^{\beta}(\mathbb{R}^{d}))_{\theta,1}^{\blacktriangleleft}=\dot
{B}_{\infty,1}^{\alpha}(\mathbb{R}^{d})\label{besov1}%
\end{equation}
with $\alpha=(1-\theta
)(\beta-\kappa)+\theta\beta$.  Here, the equivalence constant is independent of $\theta$.

\begin{theorem}
\label{ThmVE} Let $\Pi\in\mathcal{P}_{\kappa}$ and $\beta\in\mathbb{R}$.
Then
\begin{align}
\Vert f\Vert_{\dot{B}_{\Pi}^{\beta}(\mathbb{R}^{d})\cap\dot{B}_{\infty
,1}^{\beta-\kappa}(\mathbb{R}^{d})}  &  \approx\sup_{t\in(0,\infty)}%
\frac{K(t,f;\dot{B}_{\infty,1}^{\beta-\kappa}(\mathbb{R}^{d}),\dot{B}%
_{\infty,1}^{\beta}(\mathbb{R}^{d}))}{ty_{\Pi}(\frac{1}{t})}\nonumber\\
&  \approx\Vert f\Vert_{\Delta_{\theta\in(0,1)}\Big\{\frac{(\dot{B}_{\infty
,1}^{\beta-\kappa}(\mathbb{R}^{d}),\dot{B}_{\infty,1}^{\beta}(\mathbb{R}%
^{d}))_{\theta,\frac{1}{1-\theta}}^{\blacktriangleleft}}{\Pi(\frac{1}%
{1-\theta})}\Big\}}\label{VFM2}\\
&  \approx\sup_{\alpha\in(\beta-\kappa,\beta)}\frac{\Vert f\Vert_{\dot
{B}_{\infty,1}^{\alpha}(\mathbb{R}^{d})}}{\Pi(\frac{1}{\beta-\alpha})}.
\nonumber%
\end{align}
Here $y_{\Pi}$ is given by \eqref{intro:Auxy} (with $p_{0}=1$). 
\end{theorem}

\begin{proof}
According to \eqref{app1} and \eqref{la cara} (with $p_0= 1$), 
\begin{equation} \label{36novel}
\Vert f\Vert_{\Delta_{\theta\in(0,1)}\Big\{\frac{(\dot{B}_{\infty,1}%
^{\beta-\kappa}(\mathbb{R}^{d}),\dot{B}_{\infty,1}^{\beta}(\mathbb{R}%
^{d}))_{\theta,\frac{1}{1-\theta}}^{\blacktriangleleft}}{\Pi(\frac{1}%
{1-\theta})}\Big\}}   \approx\sup_{t\in(0,\infty)}\frac{K(t,f;\dot{B}_{\infty,1}^{\beta-\kappa
}(\mathbb{R}^{d}),\dot{B}_{\infty,1}^{\beta}(\mathbb{R}^{d}))}{ty_{\Pi}%
(\frac{1}{t})}.
\end{equation}
This shows the second equivalence in (\ref{VFM2}). 

We now prove the third equivalence in (\ref{VFM2}). We do this using the reiteration property for the $\Delta$-extrapolation method given in \eqref{app1} and \eqref{app1novel} (with $p=1$)
\[
\Delta_{\theta\in(0,1)}\bigg\{\frac{(\dot{B}_{\infty,1}^{\beta-\kappa
}(\mathbb{R}^{d}),\dot{B}_{\infty,1}^{\beta}(\mathbb{R}^{d}))_{\theta,\frac
{1}{1-\theta}}^{\blacktriangleleft}}{\Pi(\frac{1}{1-\theta})}\bigg\}=\Delta
_{\theta\in(0,1)}\bigg\{\frac{(\dot{B}_{\infty,1}^{\beta-\kappa}(\mathbb{R}%
^{d}),\dot{B}_{\infty,1}^{\beta}(\mathbb{R}^{d}))_{\theta,1}%
^{\blacktriangleleft}}{\Pi(\frac{1}{1-\theta})}\bigg\}
\]
combined with (\ref{besov1}), and a change of variables, to obtain%
\[
\left\Vert f\right\Vert _{\Delta_{\theta\in(0,1)}\Big\{\frac{(\dot{B}%
_{\infty,1}^{\beta-\kappa}(\mathbb{R}^{d}),\dot{B}_{\infty,1}^{\beta
}(\mathbb{R}^{d}))_{\theta,\frac{1}{1-\theta}}^{\blacktriangleleft}}{\Pi(\frac{1}{1-\theta}%
)}\Big\}} \approx \sup_{\alpha\in(\beta-\kappa,\beta)}\frac{\Vert f\Vert_{\dot
{B}_{\infty,1}^{\alpha}(\mathbb{R}^{d})}}{\Pi(\frac{1}{\beta-\alpha})}.
\]

Finally, we prove the first equivalence in \eqref{VFM2}. Note that in view of
the properties of $y_{\Pi},$ proved in Lemma \ref{Lemma1}, and the
monotonicity properties of $K$-functionals, we have
\begin{align}
\sup_{t\in(0,\infty)}\frac{K(t,f;\dot{B}_{\infty,1}^{\beta-\kappa}%
(\mathbb{R}^{d}),\dot{B}_{\infty,1}^{\beta}(\mathbb{R}^{d}))}{ty_{\Pi}%
(\frac{1}{t})} &  \approx\sup_{t\in(0,1/e)}\,\frac{K(t,f;\dot{B}_{\infty
,1}^{\beta-\kappa}(\mathbb{R}^{d}),\dot{B}_{\infty,1}^{\beta}(\mathbb{R}%
^{d}))}{t\Pi(-\log t)}\nonumber\\
&  \hspace{1cm}+\sup_{t\in(1/e,1)}\frac{K(t,f;\dot{B}_{\infty,1}^{\beta
-\kappa}(\mathbb{R}^{d}),\dot{B}_{\infty,1}^{\beta}(\mathbb{R}^{d}))}{ty_{\Pi
}(\frac{1}{t})}\nonumber\\
&  \hspace{1cm}+\sup_{t\in(1,\infty)}\,K(t,f;\dot{B}_{\infty,1}^{\beta-\kappa
}(\mathbb{R}^{d}),\dot{B}_{\infty,1}^{\beta}(\mathbb{R}^{d}))\nonumber\\
&  \approx I+II,\label{ComExtra2}%
\end{align}
where
\begin{equation*}
I:=\sup_{t\in(0,1/e)}\,\frac{K(t,f;\dot{B}_{\infty,1}^{\beta-\kappa
}(\mathbb{R}^{d}),\dot{B}_{\infty,1}^{\beta}(\mathbb{R}^{d}))}{t\Pi(-\log
t)}
\end{equation*}
and
\[
II:=\sup_{t\in(1,\infty)}K(t,f;\dot{B}_{\infty,1}^{\beta-\kappa}%
(\mathbb{R}^{d}),\dot{B}_{\infty,1}^{\beta}(\mathbb{R}^{d})).
\]

The desired equivalence in \eqref{VFM2} will then follow if we show the following
\begin{equation}
I\approx\Vert f\Vert_{\dot{B}_{\Pi}^{\beta}(\mathbb{R}^{d})}
\label{ComExtra3}%
\end{equation}%
and
\begin{equation}
II\approx\Vert f\Vert_{\dot{B}_{\infty,1}^{\beta-\kappa}(\mathbb{R}^{d})}.
\label{ComExtra4}%
\end{equation}
To prove these claims we revert to sequences of vector-valued functions using
the method of retracts. Then we see that
\begin{align*}
I  &  \approx\sup_{t\in(0,1/e)}\,\frac{1}{t\Pi(-\log t)}\,K(t,\{\dot{\Delta
}_{j}f\}_{j\in\mathbb{Z}};\ell_{1}^{\beta-\kappa}(L^{\infty}(\mathbb{R}%
^{d})),\ell_{1}^{\beta}(L^{\infty}(\mathbb{R}^{d})))\\
&  \approx\sup_{N\geq0}\,\frac{1}{2^{-N\kappa}\,\Pi(N)}\,K(2^{-N\kappa}%
,\{\dot{\Delta}_{j}f\}_{j\in\mathbb{Z}};\ell_{1}^{\beta-\kappa}(L^{\infty
}(\mathbb{R}^{d})),\ell_{1}^{\beta}(L^{\infty}(\mathbb{R}^{d}))),
\end{align*}
where in the last step we have used the monotonicity properties of $\Pi$ and
$K$-functionals. Applying known characterizations for $K$-functionals (cf. \eqref{42}), we derive
\begin{align*}
I  &  \approx\sup_{N\geq0}\,\frac{1}{2^{-N\kappa}\,\Pi(N)}\,\sum_{j=-\infty
}^{\infty}\min\{1,2^{(j-N)\kappa}\}\,2^{j(\beta-\kappa)}\Vert\dot{\Delta}%
_{j}f\Vert_{L^{\infty}(\mathbb{R}^{d})}\nonumber\\
&  \approx\sup_{N\geq0}\,\frac{1}{\Pi(N)}\,\sum_{j=-\infty}^{N}2^{j\beta}%
\Vert\dot{\Delta}_{j}f\Vert_{L^{\infty}(\mathbb{R}^{d})}\nonumber\\
&  \hspace{1cm}+\sup_{N\geq0}\,\frac{1}{2^{-N\kappa}\,\Pi(N)}\,\sum
_{j=N}^{\infty}2^{j(\beta-\kappa)}\Vert\dot{\Delta}_{j}f\Vert_{L^{\infty
}(\mathbb{R}^{d})}\nonumber\\
&  =:I_{1}+I_{2}. 
\end{align*}

%

We claim that
\begin{equation}
I_{2}\lesssim I_{1}, \label{ThmVishikLimProof3}%
\end{equation}
from where it follows that (cf. Definition \ref{DefVishik})
\begin{equation*}
I\approx I_{1}=\Vert f\Vert_{\dot{B}_{\Pi}^{\beta}(\mathbb{R}^{d})},
\end{equation*}
i.e., \eqref{ComExtra3} holds.

It remains to prove \eqref{ThmVishikLimProof3}. Let $N\geq0$. We use the
assumption (iv) in the definition of $\mathcal{P}_{\kappa}$ (cf. Definition \ref{DefPClass}), to estimate
\begin{align*}
\sum_{j=N}^{\infty}2^{j(\beta-\kappa)}\Vert\dot{\Delta}_{j}f\Vert_{L^{\infty
}(\mathbb{R}^{d})}  &  \leq\sup_{j\geq0}\,\{\Pi(j)^{-1}2^{j\beta}\Vert
\dot{\Delta}_{j}f\Vert_{L^{\infty}(\mathbb{R}^{d})}\}\,\sum_{j=N}^{\infty
}2^{-j\kappa}\Pi(j)\\
&  \approx2^{-N\kappa}\Pi(N)\,\sup_{j\geq0}\,\{\Pi(j)^{-1}2^{j\beta}\Vert
\dot{\Delta}_{j}f\Vert_{L^{\infty}(\mathbb{R}^{d})}\}\\
&  \leq2^{-N\kappa}\Pi(N)\,I_{1}.
\end{align*}
Taking the supremum over all $N\geq0$, we arrive at the desired estimate \eqref{ThmVishikLimProof3}.

We turn to the proof of (\ref{ComExtra4}). The estimate $\lesssim$ follows
trivially from the very definition of $K$-functional. The converse estimate
can be obtained from \eqref{42}. To be more precise, let $N\in\mathbb{N}$, then
\begin{align*}
K(2^{N\kappa},\{f_{j}\}_{j\in\mathbb{Z}};\ell_{1}^{\beta-\kappa}(L^{\infty
}(\mathbb{R}^{d})),\ell_{1}^{\beta}(L^{\infty}(\mathbb{R}^{d}))) &
\approx\sum_{j=-\infty}^{\infty}\min\{1,2^{(j+N)\kappa}\}\,2^{j(\beta-\kappa
)}\Vert f_{j}\Vert_{L^{\infty}(\mathbb{R}^{d})}\\
&  \geq\sum_{j=-N}^{\infty}2^{j(\beta-\kappa)}\Vert f_{j}\Vert_{L^{\infty
}(\mathbb{R}^{d})}.
\end{align*}
Hence
\[
\sup_{N\in\mathbb{N}}\,K(2^{N\kappa},\{f_{j}\}_{j\in\mathbb{Z}};\ell
_{1}^{\beta-\kappa}(L^{\infty}(\mathbb{R}^{d})),\ell_{1}^{\beta}(L^{\infty
}(\mathbb{R}^{d})))\gtrsim\sum_{j=-\infty}^{\infty}2^{j(\beta-\kappa)}\Vert
f_{j}\Vert_{L^{\infty}(\mathbb{R}^{d})}.
\]
Using once again the retraction method, we obtain (cf. \eqref{BesovDef})
\[
II\gtrsim\Vert f\Vert_{\dot{B}_{\infty,1}^{\beta-\kappa}(\mathbb{R}^{d})}.
\]
This ends the proof of \eqref{ComExtra4}.

Finally, putting together \eqref{ComExtra2}, \eqref{ComExtra3} and
\eqref{ComExtra4},
\[
\sup_{t\in(0,\infty)}\frac{K(t,f;\dot{B}_{\infty,1}^{\beta-\kappa}%
(\mathbb{R}^{d}),\dot{B}_{\infty,1}^{\beta}(\mathbb{R}^{d}))}{ty_{\Pi}%
(\frac{1}{t})}\approx\Vert f\Vert_{\dot{B}_{\Pi}^{\beta}(\mathbb{R}^{d}%
)\cap\dot{B}_{\infty,1}^{\beta-\kappa}(\mathbb{R}^{d})}.
\]
This concludes the proof of the theorem.
\end{proof}

\begin{remark}\label{Remark3.3}
The inhomogeneous counterpart of Theorem \ref{ThmVE} also holds. Note that, in
this case, the outcome is independent of $\kappa>0$ since
\[
B_{\Pi}^{\beta}(\mathbb{R}^{d})\cap B_{\infty,1}^{\beta-\kappa}(\mathbb{R}%
^{d})=B_{\Pi}^{\beta}(\mathbb{R}^{d}).
\]
Indeed, for any fixed $N\in\mathbb{N}$, we have
\begin{align*}
\Vert f\Vert_{B_{\infty,1}^{\beta-\kappa}(\mathbb{R}^{d})}  &  =\sum_{j=0}%
^{N}2^{j(\beta-\kappa)}\Vert\Delta_{j}f\Vert_{L^{\infty}(\mathbb{R}^{d})}%
+\sum_{j=N+1}^{\infty}2^{j(\beta-\kappa)}\Vert\Delta_{j}f\Vert_{L^{\infty
}(\mathbb{R}^{d})}\\
&  \leq\sum_{j=0}^{N}2^{j\beta}\Vert\Delta_{j}f\Vert_{L^{\infty}%
(\mathbb{R}^{d})}+\sum_{j=N+1}^{\infty}2^{j(\beta-\kappa)}\Vert\Delta
_{j}f\Vert_{L^{\infty}(\mathbb{R}^{d})}\\
&  \leq\Pi(N)\Vert f\Vert_{B_{\Pi}^{\beta}(\mathbb{R}^{d})}+\bigg(\sum
_{j=N+1}^{\infty}2^{-j\kappa}\Pi(j)\bigg)\,\Vert f\Vert_{B_{\Pi}^{\beta
}(\mathbb{R}^{d})}\\
&  \lesssim\Pi(N)\,\Vert f\Vert_{B_{\Pi}^{\beta}%
(\mathbb{R}^{d})},
\end{align*}
where the last step follows from the property (iv) in the definition of
$\mathcal{P}_{\kappa}$ (cf. Definition \ref{DefPClass}). 

One may also show that the characterizations provided by Theorem \ref{ThmVE} in the inhomogoneous setting are independent of $\kappa$ via
$$
	\sup_{\alpha\in(\beta-\kappa,\beta)}\frac{\Vert f\Vert_{
B_{\infty,1}^{\alpha}(\mathbb{R}^{d})}}{\Pi(\frac{1}{\beta-\alpha})} \approx \sup_{\alpha\in(\beta-\kappa_0,\beta)}\frac{\Vert f\Vert_{B_{\infty,1}^{\alpha}(\mathbb{R}^{d})}}{\Pi(\frac{1}{\beta-\alpha})}
$$
for $\kappa, \kappa_0 > 0$. The latter is an immediate consequence of the trivial embedding $B^\alpha_{\infty, 1}(\R^d) \hookrightarrow B^{\alpha - \varepsilon}_{\infty, 1}(\R^d)$ for any $\varepsilon > 0$. 

A third explanation to the independence of $\kappa$ in Theorem \ref{ThmVE}, now from an extrapolation point of view, may be found in \eqref{app2} below. 
\end{remark}

\subsection{A priori Biot-Savart estimates for a family of active scalar equations}\label{Section32}
 Recall the class of active scalar equations on $\mathbb{R}%
^{d}$ introduced in Section \ref{SVAS} (cf.  \eqref{Trans})
\begin{equation}
\omega_{t}+V\omega\cdot\nabla\omega=0,\label{Trans2}%
\end{equation}
where the assumption \eqref{E12n} on $V$ is now replaced by its homogeneous counterpart
\begin{equation}\label{20.1}
	\|V \dot{\Delta}_j f\|_{L^\infty} \lesssim 2^{j (\beta-1)} \|\dot{\Delta}_j f\|_{L^\infty}, \qquad \forall j \in \Z.  
\end{equation}
Important examples of  equations of type \eqref{Trans2} with \eqref{20.1} are $2$D Euler equations and aggregation equations ($\beta =0$), SQG and IPM equations ($\beta =1$), and generalized SQG equations ($\beta \in (0, 1)$). 

Next we obtain estimates for the modulus of continuity of $v = V \omega$ in terms of regularity properties  of $\omega$.

%

\begin{theorem}
\label{ThmEstimModVishikActive} 
Let $\beta \in \R$. Assume that the growth function $\Pi\in\mathcal{P}_{1},$ and let $\omega
\in\dot{B}_{\Pi}^{\beta}(\R^d)\cap\dot{B}_{\infty,1}^{\beta-1}(\R^d)$.
Then
\begin{equation*}
|v(x)-v(y)|\lesssim|x-y|\,y_{\Pi}\bigg(\frac{1}{|x-y|}\bigg)\,\Vert\omega
\Vert_{\dot{B}_{\Pi}^{\beta}(\R^d)\cap\dot{B}_{\infty,1}^{\beta-1}(\R^d)},
\end{equation*}
where $y_{\Pi}$ is the Yudovich function associated to $\Pi$ (cf.
\eqref{intro:Auxy} with $p_{0}=1$).
\end{theorem}

\begin{proof}

Since $v = V \omega$, it follows from \eqref{20.1} that
$$
	\|v\|_{\dot{B}^1_{\infty, 1}(\R^d)} \lesssim \|\omega\|_{\dot{B}^\beta_{\infty, 1}(\R^d)} \qquad \text{and} \qquad  \|v\|_{\dot{B}^0_{\infty, 1}(\R^d)} \lesssim \|\omega\|_{\dot{B}^{\beta-1}_{\infty, 1}(\R^d)} 
$$
and then, by interpolation, 
\begin{equation}
\Vert v\Vert_{(\dot{B}_{\infty,1}^{0}(\R^d),\dot{B}_{\infty,1}^{1}%
(\R^d))_{\theta,\frac{1}{1-\theta}}^{\blacktriangleleft}}\lesssim\Vert
\omega\Vert_{(\dot{B}_{\infty,1}^{\beta-1}(\R^d),\dot{B}_{\infty,1}^{\beta
}(\R^d))_{\theta,\frac{1}{1-\theta}}^{\blacktriangleleft}}.\label{316}%
\end{equation}

From the trivial embedding
\begin{equation}
\dot{B}_{\infty,1}^{0}(\R^d)\hookrightarrow L^{\infty}(\R^d) \label{317}%
\end{equation}
and the classical Bernstein inequality for entire functions of exponential
type (cf. \cite[p. 116]{N75}), we derive that
\begin{align*}
\Vert\nabla v\Vert_{L^{\infty}(\R^d)}  &  \lesssim\Vert\nabla v\Vert
_{\dot{B}_{\infty,1}^{0}(\R^d)}=\sum_{j=-\infty}^{\infty}\Vert\nabla
\dot{\Delta}_{j}v\Vert_{L^{\infty}(\R^d)}\\
&  \lesssim\sum_{j=-\infty}^{\infty}2^{j}\,\Vert\dot{\Delta}_{j}%
v\Vert_{L^{\infty}(\R^d)}=\Vert v\Vert_{\dot{B}_{\infty,1}^{1}(\R^d)}.
\end{align*}
In other words,
\begin{equation}
\dot{B}_{\infty,1}^{1}(\R^d)\hookrightarrow\dot{W}_{\infty}^{1}(\R^d).
\label{318}%
\end{equation}

According to \eqref{317} and \eqref{318},
\begin{equation*}
\Vert v\Vert_{(L^{\infty}(\R^d),\dot{W}_{\infty}^{1}(\R^d))_{\theta
,\frac{1}{1-\theta}}^{\blacktriangleleft}}\lesssim\Vert v\Vert_{(\dot
{B}_{\infty,1}^{0}(\R^d),\dot{B}_{\infty,1}^{1}(\R^d))_{\theta,\frac
{1}{1-\theta}}^{\blacktriangleleft}},
\end{equation*}
which implies (cf. \eqref{316})
\[
\Vert v\Vert_{(L^{\infty}(\R^d),\dot{W}_{\infty}^{1}(\R^d))_{\theta
,\frac{1}{1-\theta}}^{\blacktriangleleft}}\lesssim\Vert\omega\Vert_{(\dot
{B}_{\infty,1}^{\beta-1}(\R^d),\dot{B}_{\infty,1}^{\beta}(\R^d
))_{\theta,\frac{1}{1-\theta}}^{\blacktriangleleft}}.
\]
Multiplying this inequality by $1/\Pi((1-\theta)^{-1})$ and taking supremum
over all $\theta\in(0,1)$, we arrive at the extrapolation inequality
\begin{equation}
\Vert v\Vert_{\Delta_{\theta\in(0,1)}\Big\{\frac{(L^{\infty}(\R^d),\dot
{W}_{\infty}^{1}(\R^d))_{\theta,\frac{1}{1-\theta}}^{\blacktriangleleft}%
}{\Pi(\frac{1}{1-\theta})}\Big\}}\lesssim\Vert\omega\Vert_{\Delta_{\theta
\in(0,1)}\Big\{\frac{(\dot{B}_{\infty,1}^{\beta-1}(\R^d),\dot{B}_{\infty
,1}^{\beta}(\R^d))_{\theta,\frac{1}{1-\theta}}^{\blacktriangleleft}}%
{\Pi(\frac{1}{1-\theta})}\Big\}}.\label{320}%
\end{equation}

The right-hand side of \eqref{320} was computed in Theorem \ref{ThmVE}:
\begin{equation}
\Delta_{\theta\in(0,1)}\bigg\{\frac{(\dot{B}_{\infty,1}^{\beta-1}(\R^d
),\dot{B}_{\infty,1}^{\beta}(\R^d))_{\theta,\frac{1}{1-\theta}%
}^{\blacktriangleleft}}{\Pi(\frac{1}{1-\theta})}\bigg\}=\dot{B}_{\Pi}^{\beta
}(\R^d)\cap\dot{B}_{\infty,1}^{\beta-1}(\R^d).\label{321}%
\end{equation}
While the left-hand side of \eqref{320} was already characterized in
\eqref{225new} (replace $\Theta_{1}$ by $\Pi$ and let $p_{0}=1$) as
\begin{equation}
\Vert v\Vert_{\Delta_{\theta\in(0,1)}\Big\{\frac{(L^{\infty}(\R^d),\dot
{W}_{\infty}^{1}(\R^d))_{\theta,\frac{1}{1-\theta}}^{\blacktriangleleft}%
}{\Pi(\frac{1}{1-\theta})}\Big\}}\approx\sup_{x,y\in\R^d}\,\frac
{|v(x)-v(y)|}{\inf_{t>|x-y|}ty_{\Pi}(\frac{1}{t})}.\label{322}%
\end{equation}
Inserting \eqref{321} and \eqref{322} into \eqref{320}, we have
\[
|v(x)-v(y)|\lesssim|x-y|y_{\Pi}\bigg(\frac{1}{|x-y|}\bigg)\,\Vert\omega
\Vert_{\dot{B}_{\Pi}^{\beta}(\R^d)\cap\dot{B}_{\infty,1}^{\beta-1}(\R^d)}.
\]
\end{proof}

\subsection{Comparison between Theorem \ref{ThmEstimModVishikActive} with $\beta =0$ and \eqref{18c}}

\label{SectionCom} As already mentioned in Remark \ref{Remark5}, both
conditions \eqref{IntroPi} and \eqref{introOsB} are equivalent for $\Pi
\in\mathcal{P}_{1}$ (cf. also Lemma \ref{Lemma1}). In this section, we investigate the relationships between
the function spaces involved in \eqref{IntroViSpa} and Theorem \ref{ThmEstimModVishikActive} , namely,
\[
B_{\Pi}(\mathbb{R}^{2})\cap L^{p_{0}}(\mathbb{R}^{2}),\qquad p_{0}\in(1,2),
\]
and
\[
\dot{B}_{\Pi}(\mathbb{R}^{2})\cap\dot{B}_{\infty,1}^{-1}(\mathbb{R}^{2}),
\]
respectively.

\begin{proposition}
\label{PropoComp1} Let $p_{0} \in(1, d)$. Then
\[
B_{\Pi}(\mathbb{R}^{d}) \cap L^{p_{0}}(\mathbb{R}^{d}) \hookrightarrow\dot
{B}_{\Pi}(\mathbb{R}^{d}) \cap\dot{B}^{-1}_{\infty, 1}(\mathbb{R}^{d}).
\]

\end{proposition}

\begin{proof}
Assume $f\in B_{\Pi}(\mathbb{R}^{d})\cap L^{p_{0}}(\mathbb{R}^{d})$. The norm
of $f$ in $\dot{B}_{\Pi}(\mathbb{R}^{d})$ can be estimated as follows. For
each $N\in\mathbb{N}$, we can apply the classical Nikolskii's inequality for
entire functions of exponential type (cf. \cite[Theorem 3.3.5, p. 126]{N75})
and basic multiplier assertions in order to get
\begin{align*}
\sum_{j=-\infty}^{N}\Vert\dot{\Delta}_{j}f\Vert_{L^{\infty}(\mathbb{R}^{d})}
&  =\sum_{j=-\infty}^{0}\Vert\dot{\Delta}_{j}f\Vert_{L^{\infty}(\mathbb{R}%
^{d})}+\sum_{j=1}^{N}\Vert\dot{\Delta}_{j}f\Vert_{L^{\infty}(\mathbb{R}^{d}%
)}\\
&  \lesssim\sum_{j=-\infty}^{0}2^{jd/p_{0}}\Vert\dot{\Delta}_{j}%
f\Vert_{L^{p_{0}}(\mathbb{R}^{d})}+\Pi(N)\Vert f\Vert_{B_{\Pi}(\mathbb{R}%
^{d})}\\
&  \lesssim\Vert f\Vert_{L^{p_{0}}(\mathbb{R}^{d})}+\Pi(N)\Vert f\Vert
_{B_{\Pi}(\mathbb{R}^{d})}\\
&  \lesssim\Pi(N)\Vert f\Vert_{B_{\Pi}(\mathbb{R}^{d})\cap L^{p_{0}%
}(\mathbb{R}^{d})}.
\end{align*}
Hence
\begin{equation}
B_{\Pi}(\mathbb{R}^{d})\cap L^{p_{0}}(\mathbb{R}^{d})\hookrightarrow\dot
{B}_{\Pi}(\mathbb{R}^{d}) \label{414}%
\end{equation}
for any $p_{0}\in(1,\infty)$.

On the other hand, the $\dot{B}_{\infty,1}^{-1}$-norm of $f$ can be split
into
\begin{equation}
\Vert f\Vert_{\dot{B}_{\infty,1}^{-1}(\mathbb{R}^{d})}=\sum_{j=-\infty
}^{\infty}2^{-j}\Vert\dot{\Delta}_{j}f\Vert_{L^{\infty}(\mathbb{R}^{d})}=I+II,
\label{Split1}%
\end{equation}
where
\begin{equation}
I:=\sum_{j=-\infty}^{0}2^{-j}\Vert\dot{\Delta}_{j}f\Vert_{L^{\infty
}(\mathbb{R}^{d})} \label{I1}%
\end{equation}
and
\begin{equation}
II:=\sum_{j=1}^{\infty}2^{-j}\Vert\dot{\Delta}_{j}f\Vert_{L^{\infty
}(\mathbb{R}^{d})}. \label{l2}%
\end{equation}

We can estimate $I$ applying Nikolskii's inequality:
\begin{equation}
\label{Split2}I \lesssim\sum_{j=-\infty}^{0} 2^{-j (1-\frac{d}{p_{0}})}
\|\dot{\Delta}_{j} f\|_{L^{p_{0}}(\mathbb{R}^{d})} \lesssim\|f\|_{L^{p_{0}%
}(\mathbb{R}^{d})} \, \sum_{j=0}^{\infty}2^{j (1- \frac{d}{p_{0}})}
\lesssim\|f\|_{L^{p_{0}}(\mathbb{R}^{d})},
\end{equation}
where the last step follows from the assumption $p_{0} < d$.

To estimate $II$, we can argue as follows. The fact that $f\in B_{\Pi
}(\mathbb{R}^{d})$ implies that
\[
\Vert\Delta_{j}f\Vert_{L^{\infty}(\mathbb{R}^{d})}\leq\Pi(j)\Vert
f\Vert_{B_{\Pi}(\mathbb{R}^{d})},\qquad j\in\mathbb{N}.
\]
Therefore, taking into account that $\Pi\in\mathcal{P}_{1}$ (cf. item (iv) in
Definition \ref{DefPClass}), we obtain
\begin{equation}
II=\sum_{j=1}^{\infty}2^{-j}\Vert\Delta_{j}f\Vert_{L^{\infty}(\mathbb{R}^{d}%
)}\leq\Vert f\Vert_{B_{\Pi}(\mathbb{R}^{d})}\,\sum_{j=1}^{\infty}2^{-j}%
\Pi(j)\lesssim\Vert f\Vert_{B_{\Pi}(\mathbb{R}^{d})}. \label{416}%
\end{equation}
Inserting \eqref{Split2} and \eqref{416} into \eqref{Split1}, one achieves
\begin{equation}
B_{\Pi}(\mathbb{R}^{d})\cap L^{p_{0}}(\mathbb{R}^{d})\hookrightarrow\dot
{B}_{\infty,1}^{-1}(\mathbb{R}^{d}). \label{418}%
\end{equation}

The combination of \eqref{414} and \eqref{418} yields the desired embedding
\[
B_{\Pi}(\mathbb{R}^{d}) \cap L^{p_{0}}(\mathbb{R}^{d}) \hookrightarrow\dot
{B}_{\Pi}(\mathbb{R}^{d}) \cap\dot{B}^{-1}_{\infty, 1}(\mathbb{R}^{d}).
\]

\end{proof}

We close this section by showing that Theorem \ref{ThmEstimModVishikActive} with $\Pi (p) \approx p$  also improves the classical Vishik's estimate formulated in terms of  \eqref{VUC}. 

\begin{proposition}
\label{PropComp3} Let $p_{0} \in(1, d)$. Then
\begin{equation}
\label{330}\emph{bmo}(\mathbb{R}^{d}) \cap L^{p_{0}}(\mathbb{R}^{d})
\hookrightarrow\emph{BMO}(\mathbb{R}^{d}) \cap\dot{B}^{-1}_{\infty,
1}(\mathbb{R}^{d}).
\end{equation}

\end{proposition}

\begin{proof}
Assume $f \in\text{bmo}(\mathbb{R}^{d}) \cap L^{p_{0}}(\mathbb{R}^{d})$. The
Besov $\dot{B}^{-1}_{\infty, 1}$-norm can be expressed as in
\eqref{Split1}--\eqref{l2}, with $I$ satisfying \eqref{Split2}. Hence
\begin{equation}
\label{Spli5}\|f\|_{\dot{B}^{-1}_{\infty, 1}(\mathbb{R}^{d})} = I + II
\lesssim\|f\|_{L^{p_{0}}(\mathbb{R}^{d})} + II.
\end{equation}

Next we estimate $II$. The well-known embedding $\text{bmo}(\mathbb{R}^{d})
\hookrightarrow B^{0}_{\infty, \infty}(\mathbb{R}^{d})$ (cf. \eqref{intro112})
implies
\begin{equation}
\label{Split3}II = \sum_{j=1}^{\infty}2^{-j} \|\Delta_{j} f\|_{L^{\infty
}(\mathbb{R}^{d})} \lesssim\|f\|_{B^{0}_{\infty, \infty}(\mathbb{R}^{d})}
\lesssim\|f\|_{\text{bmo}(\mathbb{R}^{d})}.
\end{equation}

According to \eqref{Spli5} and \eqref{Split3},
\[
\text{bmo}(\mathbb{R}^{d}) \cap L^{p_{0}}(\mathbb{R}^{d}) \hookrightarrow
\dot{B}^{-1}_{\infty, 1}(\mathbb{R}^{d}).
\]
This together with the trivial embedding $\text{bmo}(\mathbb{R}^{d})
\hookrightarrow\text{BMO}(\mathbb{R}^{d})$ implies the desired result \eqref{330}.
\end{proof}

\subsection{Proof of Theorem \ref{Thm35a}}\label{SectionUniqAS}

	Let $\beta=1$, and let $\omega_1$, $\omega_2$ be two solutions of \eqref{Trans}  in $L^\infty([0, T]; B^1_\Pi)$, such that $\omega_1(0) = \omega_2(0)$. Let $\rho = \omega_1-\omega_2$. Then, by elementary manipulations, 
\begin{equation}
	\frac{1}{2} \frac{d}{d t} \|\rho(t)\|_{L^2}^2  = - \int (\rho V \omega_1 \cdot \nabla \rho + \rho V \rho \cdot \nabla \omega_2) \, dx= - \int \rho V \rho \cdot \nabla \omega_2 \, dx, \label{E13}
\end{equation}
where we have used the divergence free assumption of $V$ in the last step. 

Suppose that $\alpha_0 \in (0, 1)$ and let $\alpha \in (0, \alpha_0)$. Recall the well-known fact that the Besov pair\footnote{To avoid density  issues, the space  $B^s_{p, q}$ with $p=\infty$ and/or $q=\infty$ may be defined as the closure of $\mathcal{S}$ in $B^s_{p, q}$. See \cite[Section 2.11]{Tri83} for further details.} $(B^{\alpha}_{1, \infty}, B^{-\alpha}_{\infty, 1})$ is in duality pairing within $(\mathcal{S}, \mathcal{S}')$.
Applying Theorem \ref{ThmVE} (specially Remark \ref{Remark3.3}) with $\beta =1$, we obtain
\begin{align}
	\bigg| \int \rho V \rho \cdot \nabla \omega_2 \, dx \bigg|& \lesssim \|\nabla \omega_2(t)\|_{B^{-\alpha}_{\infty, 1}} \|\rho(t) V \rho(t)\|_{B^\alpha_{1, \infty}}  \nonumber  \\
	& \lesssim \|\omega_2(t)\|_{B^{1-\alpha}_{\infty, 1}} \|\rho(t) V \rho(t)\|_{B^\alpha_{1, \infty}}  \nonumber  \\
	& \lesssim \Pi \bigg(\frac{1}{\alpha} \bigg) \|\omega_2(t)\|_{B^1_\Pi}  \|\rho(t) V \rho(t)\|_{B^\alpha_{1, \infty}}. \label{E14}
\end{align}
On the other hand, according to the interpolation formula \eqref{IntBesovNormalized} we have (recall that the normalized interpolation constant in \eqref{IntroNormInt} is given by $c_{\theta, \infty} = 1$):
\begin{equation}\label{E6.13}
	\|\rho(t) V \rho(t)\|_{B^\alpha_{1, \infty}} \approx \|\rho(t) V \rho(t)\|_{(B^0_{1, \infty}, B^{\alpha_0}_{1, \infty})_{\frac{\alpha}{\alpha_0}}, \infty}. 
\end{equation}
At this point, we can invoke the following remarkable property (cf. \cite{M16}) of the normalized interpolation scales $\bar{A}_{\theta,p}^{\blacktriangleleft}$: the  normalization constant $c_{\theta, p}$  (cf. \eqref{IntroNormInt})  is  optimal  for the family of interpolation inequalities\footnote{In the special case $p=\infty$, the inequality \eqref{E6} is a simple consequence of the fact that $K(t, f; A_0, A_1) \leq \min\{\|f\|_{A_0}, t \|f\|_{A_1}\}$.}
\begin{equation}\label{E6}
	\|f\|_{\bar{A}_{\theta,p}^{\blacktriangleleft}} \leq \|f\|_{A_0}^{1-\theta} \|f\|_{A_1}^{\theta}, \qquad f \in A_0 \cap A_1. 
\end{equation}
In particular,  putting together \eqref{E6.13} and \eqref{E6}, we obtain
$$
	\|\rho(t) V \rho(t)\|_{B^\alpha_{1, \infty}}  \lesssim \|\rho(t) V \rho(t)\|_{B^0_{1, \infty}}^{1-\frac{\alpha}{\alpha_0}} \|\rho(t) V \rho(t)\|_{B^{\alpha_0}_{1, \infty}}^{\frac{\alpha}{\alpha_0}}.
$$
The trivial embedding $L^1 \hookrightarrow B^0_{1, \infty}$ and H\"older's inequality, now yield
\begin{align}
	\|\rho(t) V \rho(t)\|_{B^\alpha_{1, \infty}} & \lesssim \|\rho(t) V \rho(t)\|_{L^1}^{1-\frac{\alpha}{\alpha_0}} \|\rho(t) V \rho(t)\|_{B^{\alpha_0}_{1, \infty}}^{\frac{\alpha}{\alpha_0}} \nonumber\\
	& \leq (\|\rho(t)\|_{L^2} \|V \rho(t)\|_{L^2})^{1-\frac{\alpha}{\alpha_0}}  \|\rho(t) V \rho(t)\|_{B^{\alpha_0}_{1, \infty}}^{\frac{\alpha}{\alpha_0}} \nonumber \\
	& \lesssim \|\rho(t)\|_{L^2}^{2(1-\frac{\alpha}{\alpha_0})} \|\rho(t) V \rho(t)\|_{B^{\alpha_0}_{1, \infty}}^{\frac{\alpha}{\alpha_0}},  \label{E15}
\end{align}
where in the last estimate we used the assumption $V: L^2 \to L^2$ . 
Furthermore, we claim that 
\begin{equation}\label{Product}
\|\rho(t) V \rho(t)\|_{B^{\alpha_0}_{1, \infty}} < \infty.
\end{equation} To see this observe that from the trivial embeddings (recall that all function spaces are defined on $\Omega = \T^2$)
$$B^1_\Pi \hookrightarrow B^{\frac{1+\alpha_0}{2}}_{\infty, 1} \hookrightarrow B^{\alpha_0}_{1, \infty} \cap B^0_{\infty, 1} \hookrightarrow B^{\alpha_0}_{1, \infty}  \cap L^\infty$$
we have that $\rho(t) \in B^{\alpha_0}_{1, \infty}  \cap L^\infty$. In addition, $V \rho(t) \in B^{\alpha_0}_{1, \infty}  \cap L^\infty$ since $V \rho(t) \in B^1_\Pi$ (cf.  \eqref{E12n} with $\beta =1$). Therefore the desired result will follow if we can show that
\begin{equation}\label{ProductBounded}
	\|f g\|_{B^s_{p, q}} \lesssim \|f\|_{B^s_{p, q}} \|g\|_{L^\infty} + \|f\|_{L^\infty} \|g\|_{B^s_{p, q}}, \qquad s \in (0, 1). 
\end{equation}
The proof of  \eqref{ProductBounded} itself follows trivially\footnote{For the purposes of this paper, it is enough to apply \eqref{ProductBounded} with $s \in (0, 1)$, but we would like to mention that  \eqref{ProductBounded} is known to be true  for any arbitrary smoothness parameter $s > 0$, cf.  \cite[Corollary 2.86]{BCD11}. However, the general approach given in \cite{BCD11} is not elementary anymore and hinges upon Bony's decompositions.} writing
 \begin{equation*}
 fg(x+h)-fg(x) = f(x+h)(g(x+h)-g(x)) + g(x)(f(x+h)-f(x))
 \end{equation*}
  together with the well-known characterization of $B^s_{p, q}$ in terms of differences $$\|f\|_{B^s_{p, q}} \approx \|f\|_{L^p} + \bigg(\int_{|h| < 1} |h|^{-s q} \|f(\cdot + h)-f\|_{L^p}^q \, \frac{dh}{|h|^d} \bigg)^{1/q}, \qquad s \in (0, 1).
  $$

In view of \eqref{Product}, the estimate \eqref{E15} now reads as follows
\begin{equation*}
	\|\rho(t) V \rho(t)\|_{B^\alpha_{1, \infty}} \lesssim  \|\rho(t)\|_{L^2}^{2(1-\frac{\alpha}{\alpha_0})}.
\end{equation*}
Consequently, upon updating \eqref{E14}, we obtain
$$
	\bigg| \int \rho V \rho \cdot \nabla \omega_2 \, dx \bigg| \lesssim   \|\rho(t)\|_{L^2}^{2} \Pi \bigg(\frac{1}{\alpha} \bigg) \bigg(\frac{1}{\|\rho(t)\|^2_{L^2}} \bigg)^{\alpha},
$$
where the underlying equivalence constant is independent of $\alpha \in (0, \alpha_0)$. As a consequence   (cf. \eqref{intro:Auxy}) we have
$$
	\bigg| \int \rho V \rho \cdot \nabla \omega_2 \, dx \bigg| \lesssim  \|\rho(t)\|_{L^2}^{2} y_\Pi \bigg(\frac{1}{\|\rho(t)\|^2_{L^2}} \bigg)
$$
and, by \eqref{E13}, 
$$
	 \frac{d}{d t}  \|\rho(t)\|_{L^2}^2 \lesssim   \|\rho(t)\|_{L^2}^{2} y_\Pi \bigg(\frac{1}{\|\rho(t)\|^2_{L^2}} \bigg).
$$
An application of Osgood's lemma allows us to conclude that $\rho(t) =0$ for all $t$. \qed 

\begin{remark}\label{RemarkR}
	As already mentioned in Section \ref{SVAS}, the above proof can be easily adapted to work with $\Omega = \T^d, \R^d$. Note that the case $\Omega= \R^d$ requires the following additional (although very mild) assumption:  $\omega \in L^\infty([0, T]; B^{\alpha_0}_{1, \infty}(\R^d))$ where $\alpha_0 \in (0, 1)$ can be chosen arbitrary small. Moreover, in view of Theorem \ref{ThmEstimModVishikActive} and specially Section \ref{SectionCom}, it is clear that Theorem \ref{Thm35a} can be improved if one works with the homogeneous space $\dot{B}^1_\Pi$, however we prefer to avoid such a technical issue in this paper.  
\end{remark}

\subsection{Proof of Theorem \ref{IntroThm6}}
Assume that $\omega_1$ and $\omega_2$ are two weak solutions of \eqref{Trans} with $\beta =1$ in \eqref{E12n}. Assume further that  $\omega_1(0)= \omega_2(0)$. Let $\rho = \omega_1-\omega_2$ and $\Delta \psi = \rho$. According to \eqref{17.8}  and \eqref{17.1}, we have
	\begin{align}
	\frac{1}{2} \frac{d}{dt} \int |\nabla \psi|^2 \, dx &= - \int (A(\omega_1)-A(\omega_2)) (\omega_1-\omega_2) \, dx  \nonumber\\
	& \hspace{1cm} + \sum_{i, l}	\int \bigg[ \partial_l \Big(\frac{1}{2} |\partial_i \psi|^2 \Big) (V \omega_1)_l  + \partial_i \psi \partial_l \psi \partial_i (V \omega_1)_l \bigg] \, dx \label{16.8}  \\
	&\hspace{1cm}+ \int \nabla \psi \cdot \omega_2 V \rho \, dx. \nonumber
	\end{align}
	
	By the monotonicity assumption on $A$, it is clear that
	\begin{equation}\label{16.10}
		 - \int (A(\omega_1)-A(\omega_2)) (\omega_1-\omega_2) \, dx \leq 0. 
	\end{equation}
	
	Let $\alpha_0 \in (0, 1)$ be a fixed parameter. Given any $\alpha \in (0, \alpha_0)$, we can apply the duality theory of Besov spaces, the standing assumption on $V$ given by \eqref{E12n} and Theorem \ref{ThmVE} (with $\beta =0$) in order to get
	\begin{align}
		\bigg| \int \partial_l \Big(\frac{1}{2} |\partial_i \psi|^2 \Big) (V \omega_1)_l \, dx  \bigg| & \lesssim \|(\partial_i \psi)^2\|_{B^{\alpha}_{1, \infty}} \|\partial_l (V \omega_1)_l\|_{B^{-\alpha}_{\infty, 1}} \nonumber \\
		& \lesssim  \|(\partial_i \psi)^2\|_{B^{\alpha}_{1, \infty}} \|\omega_1\|_{B^{-\alpha}_{\infty, 1}} \nonumber\\
		& \lesssim \|(\partial_i \psi)^2\|_{B^{\alpha}_{1, \infty}}  \Pi \bigg(\frac{1}{\alpha} \bigg) \| \omega_1\|_{B_\Pi}. \label{16.2}
	\end{align}
	Moreover, it is well known that the following interpolation formula holds
	\begin{equation}\label{16.1}
		B^\alpha_{1, \infty} = (B^0_{1, \infty}, B^{\alpha_0}_{1, \infty})_{\frac{\alpha}{\alpha_0}, \infty}.
	\end{equation}
	However, for our purposes we require precise information on the equivalence constants with respect to $\alpha$. This can be achieved using the normalized interpolation constants in \eqref{IntroNormInt}. In particular, according to Remark \ref{RemarkA} from Appendix \ref{sec:A0} (note that $c_{\theta, \infty} = 1$), 
	$$
		\|f\|_{B^\alpha_{1, \infty}} \approx \|f\|_{(B^0_{1, \infty}, B^{\alpha_0}_{1, \infty})_{\frac{\alpha}{\alpha_0}, \infty}}
	$$
	uniformly  in $\alpha$. Then, by \eqref{E6}, 
	\begin{align}
		\|(\partial_i \psi)^2\|_{B^{\alpha}_{1, \infty}} & \approx \|(\partial_i \psi)^2\|_{(B^0_{1, \infty}, B^{\alpha_0}_{1, \infty})_{\frac{\alpha}{\alpha_0}, \infty}}  \nonumber\\
		& \leq \|(\partial_i \psi)^2\|_{B^0_{1, \infty}}^{1-\frac{\alpha}{\alpha_0}}  \|(\partial_i \psi)^2\|_{B^{\alpha_0}_{1, \infty}}^{\frac{\alpha}{\alpha_0}}  \nonumber \\
		& \lesssim \|(\partial_i \psi)^2\|_{L^1}^{1-\frac{\alpha}{\alpha_0}}  \|(\partial_i \psi)^2\|_{B^{\alpha_0}_{1, \infty}}^{\frac{\alpha}{\alpha_0}}  \nonumber \\
		& = \|\partial_i \psi\|_{L^2}^{2(1-\frac{\alpha}{\alpha_0})}  \|(\partial_i \psi)^2\|_{B^{\alpha_0}_{1, \infty}}^{\frac{\alpha}{\alpha_0}}. \label{16.3}  
	\end{align}
	Moreover, we claim that
	\begin{equation}\label{Claim16}
		\|(\partial_i \psi)^2\|_{B^{\alpha_0}_{1, \infty}} < \infty. 	
	\end{equation}
	Indeed, we  now show that $\nabla \psi \in L^\infty \cap B^{\alpha_0}_{1, \infty}$. By Theorem \ref{ThmVE} (see also Remark \ref{Remark3.3})
	\begin{equation*}
		\|\nabla \psi\|_{L^\infty}  = \|\nabla (-\Delta)^{-1} \rho \|_{L^\infty} \lesssim \|\nabla (-\Delta)^{-1} \rho\|_{B^0_{\infty, 1}} \approx \|\rho\|_{B^{-1}_{\infty, 1}} \lesssim \|\rho\|_{B_\Pi},
	\end{equation*}	
	and (note that the involved Besov spaces are defined over $\Omega = \T^2$)
	$$
		\|\nabla \psi\|_{B^{\alpha_0}_{1, \infty}} \approx  \| \rho \|_{B^{\alpha_0-1}_{1, \infty}} \lesssim \|\rho\|_{B^{\alpha_0-1}_{\infty, 1}}  \lesssim \|\rho\|_{B_\Pi}. 
	$$
	In light of \eqref{ProductBounded}, we have
	$$
		\|(\partial_i \psi)^2\|_{B^{\alpha_0}_{1, \infty}} \lesssim \|\nabla \psi\|_{L^\infty} \|\nabla \psi\|_{B^{\alpha_0}_{1, \infty}} \lesssim \|\rho\|_{B_\Pi}^2,
	$$
	hence \eqref{Claim16} holds. 
	
	Combining \eqref{16.2}, \eqref{16.3} and \eqref{Claim16} yields
	\begin{equation}\label{16.4}
		\bigg| \int \partial_l \Big(\frac{1}{2} |\partial_i \psi|^2 \Big) (V \omega_1)_l \, dx  \bigg| \lesssim  \Pi \bigg(\frac{1}{\alpha} \bigg)  \|\nabla \psi\|_{L^2}^{2(1-\frac{\alpha}{\alpha_0})}. 
	\end{equation}
	
	Notice that the same reasoning we used to achieve \eqref{16.4} also yields
	\begin{equation}\label{16.5}
		\bigg|\int \partial_i \psi \partial_l \psi \partial_i (V \omega_1)_l  \, dx \bigg| \lesssim   \Pi \bigg(\frac{1}{\alpha} \bigg)   \|\nabla \psi\|_{L^2}^{2(1-\frac{\alpha}{\alpha_0})}
	\end{equation}
	and
	\begin{equation}\label{16.6}
	 \bigg| \int \nabla \psi \cdot \omega_2 V \rho \, dx \bigg|  \lesssim \Pi \bigg(\frac{1}{\alpha} \bigg)  \|\nabla \psi \cdot V \rho\|_{L^1}^{1-\frac{\alpha}{\alpha_0}} \|\nabla \psi \cdot V \rho\|_{B^{\alpha_0}_{1, \infty}}^{\frac{\alpha}{\alpha_0}}. 
	\end{equation}
	
	Next we estimate the right-hand side of \eqref{16.6}. Since $V: \dot{H}^{-1} \to L^2$ and $-\Delta \psi = \rho$,
	$$
		\|\nabla \psi \cdot V \rho\|_{L^1} \leq \|\nabla \psi\|_{L^2} \|V \rho\|_{L^2} \lesssim \|\nabla \psi\|_{L^2} \|(-\Delta)^{1/2} \psi\|_{L^2} \approx \|\nabla \psi\|^2_{L^2}. 
	$$
	Moreover, $V \rho$ satisfies (cf. \eqref{E12n})
	$$
		\|V \rho\|_{L^\infty} \lesssim \|V \rho\|_{B^0_{\infty, 1}}  \lesssim \|\rho\|_{B^{-1}_{\infty, 1}} \lesssim \|\rho\|_{B_\Pi}
	$$
	and
	$$
		\|V \rho\|_{B^{\alpha_0}_{1, \infty}} \lesssim  \|V \rho\|_{B^{\alpha_0}_{\infty, 1}} \lesssim  \| \rho\|_{B^{\alpha_0-1}_{\infty, 1}} \lesssim \|\rho\|_{B_\Pi}, 
	$$
	therefore \eqref{ProductBounded} enables us to conclude that $\nabla \psi \cdot V \rho \in B^{\alpha_0}_{1, \infty}$. Hence \eqref{16.6} gives 
	\begin{equation}\label{16.7}
		 \bigg| \int \nabla \psi \cdot \omega_2 V \rho \, dx \bigg|  \lesssim \Pi \bigg(\frac{1}{\alpha} \bigg)  \|\nabla \psi\|_{L^2}^{2(1-\frac{\alpha}{\alpha_0})}. 
	\end{equation}
	
	Plugging \eqref{16.10}, \eqref{16.4},  \eqref{16.5} and  \eqref{16.7} into \eqref{16.8}, we see that the following estimate holds 
	$$
		\frac{d}{dt} \|\nabla \psi\|_{L^2}^2 \lesssim  \Pi \bigg(\frac{1}{\alpha} \bigg)  \|\nabla \psi\|_{L^2}^{2(1-\frac{\alpha}{\alpha_0})}
	$$
	for every $\alpha \in (0, \alpha_0)$. Then, by \eqref{intro:Auxy},
	$$
		\frac{d}{dt} \|\nabla \psi(t)\|_{L^2}^2 \lesssim \|\nabla \psi(t)\|_{L^2}^2 y_\Pi \bigg(\frac{1}{\|\nabla \psi(t)\|_{L^2}^2} \bigg)
	$$
	and an application of Osgood's lemma leads to $\|\nabla \psi(t)\|_{L^2} =0$ on $[0, T]$ and then $\psi(t) =0$ on $[0, T]$.   \qed
\begin{remark}
	The methodology of Remark \ref{RemarkR} can be easily adapted to obtain an analog of  Theorem \ref{IntroThm6}  for smoothing linear operators $V$ of order one in $\R^2$.  For this case, the following  additional (but rather natural) assumption is needed: $\omega \in L^\infty([0, T]; B^{\alpha_0-1}_{1, \infty}(\R^2))$ for some arbitrary small value of $\alpha_0 \in (0, 1)$. 
\end{remark}

%

\appendix

\section{Interpolation and Extrapolation: An Atlas\label{sec:A0}}

In order to help non experts in extrapolation, in this appendix we give a
summary of results used in the paper, with documentation, commentaries and
examples. We keep the notation and assumptions laid out in the previous
sections. For the sake of convenience we also recall the location of some
basic definitions.

\subsection{Interpolation of Besov spaces via retraction revisited\label{sec:retract}}

A common technique used in interpolation theory is to translate interpolation
of function spaces into equivalent interpolation problems for sequence spaces
(\textquotedblleft the method of retracts\textquotedblright). In particular,
in this paper the $\ell_{1}^{\beta}(L^{\infty}(\mathbb{R}^{d}))$ spaces of
sequences of vector-valued functions play an important role. We say that
$\{f_{j}\}_{j\in\mathbb{Z}}\in\ell_{1}^{\beta}(L^{\infty}(\mathbb{R}^{d}))$
if
\begin{equation}\label{DefSeqSpaces}
\Vert\{f_{j}\}_{j\in\mathbb{Z}}\Vert_{\ell_{1}^{\beta}(L^{\infty}%
(\mathbb{R}^{d}))}:=\sum_{j\in\mathbb{Z}}2^{j\beta}\Vert f_{j}\Vert
_{L^{\infty}(\mathbb{R}^{d})}<\infty.
\end{equation}
Their usefulness for us is that we can translate the interpolation of $\dot
{B}_{\infty,1}^{\beta}(\mathbb{R}^{d})$ spaces into interpolation of $\ell
_{1}^{\beta}(L^{\infty}(\mathbb{R}^{d})).$ This is effected by showing that
the map
\[
f\mapsto\{\dot{\Delta}_{j}f\}_{j\in\mathbb{Z}}%
\]
defines a retract from $\dot{B}_{\infty,1}^{\beta}(\mathbb{R}^{d})$ onto
$\ell_{1}^{\beta}(L^{\infty}(\mathbb{R}^{d}))$ (cf. \cite[Definition 6.4.1 
and Theorem 6.4.3, pages 150-152]{BL76}). In particular, $K$-functionals relative to vector-valued sequence spaces can be explicitly computed, cf. \cite[p. 120]{BL76}. 

\begin{example}
\label{ejemplo:markao} We have
\begin{equation}\label{42}
K(t,\{f_{j}\}_{j\in\mathbb{Z}};\ell_{1}^{\beta-\kappa}(L^{\infty}%
(\mathbb{R}^{d})),\ell_{1}^{\beta}(L^{\infty}(\mathbb{R}^{d})))\approx
\sum_{j=-\infty}^{\infty}\min\{1,2^{j\kappa}t\}\,2^{j(\beta-\kappa)}\,\Vert
f_{j}\Vert_{L^{\infty}(\mathbb{R}^{d})}.
\end{equation}
In fact, more general statements are available in the literature. For example, \eqref{42} holds true with $L^\infty(\R^d)$ replaced by any Banach space $X$. However, for our purposes, it is enough to restrict attention to \eqref{42}. 
\end{example}

The next result is a (possible) slight improvement of a known result (cf. \cite[Section 6.4]{BL76}) since we
provide sharp constants.

\begin{theorem}
\label{teo:nec}Let $-\infty<s_{0}<s_{1}<\infty,\,\theta\in(0,1)$, and
$s=(1-\theta)s_{0}+\theta s_{1}$. Then
\begin{equation*}
(\dot{B}_{\infty,1}^{s_{0}}(\mathbb{R}^{d}),\dot{B}_{\infty,1}^{s_{1}%
}(\mathbb{R}^{d}))_{\theta,1}^{\blacktriangleleft}=\dot{B}_{\infty,1}%
^{s}(\mathbb{R}^{d})
\end{equation*}
with underlying equivalence constants independent of $\theta$.
\end{theorem}

\begin{proof}
We use the retraction method of interpolation,
which implies that
\[
\Vert f\Vert_{(\dot{B}_{\infty,1}^{s_{0}}(\mathbb{R}^{d}),\dot{B}_{\infty
,1}^{s_{1}}(\mathbb{R}^{d}))_{\theta,1}}\approx\Vert\{\dot{\Delta}%
_{j}f\}_{j\in\mathbb{Z}}\Vert_{(\ell_{1}^{s_{0}}(L^{\infty}(\mathbb{R}%
^{d})),\ell_{1}^{s_{1}}(L^{\infty}(\mathbb{R}^{d})))_{\theta,1}}%
\]
with constants independent of $\theta$. On the other hand, it follows from
\eqref{42} that, for $\{f_{j}\}_{j\in\mathbb{Z}}\in\ell_{1}^{s_{0}%
}(L^{\infty}(\mathbb{R}^{d}))+\ell_{1}^{s_{1}}(L^{\infty}(\mathbb{R}^{d}))$,
\begin{align*}
\Vert\{f_{j}\}_{j\in\mathbb{Z}}\Vert_{(\ell_{1}^{s_{0}}(L^{\infty}%
(\mathbb{R}^{d})),\ell_{1}^{s_{1}}(L^{\infty}(\mathbb{R}^{d})))_{\theta,1}}
&  \approx\sum_{\nu=-\infty}^{\infty}2^{-\theta\nu(s_{1}-s_{0})}K(2^{\nu
(s_{1}-s_{0})},\{f_{j}\}_{j\in\mathbb{Z}};\ell_{1}^{s_{0}}(L^{\infty
}(\mathbb{R}^{d})),\ell_{1}^{s_{1}}(L^{\infty}(\mathbb{R}^{d})))\\
&  \hspace{-4cm}\approx\sum_{\nu=-\infty}^{\infty}2^{-\theta\nu(s_{1}-s_{0}%
)}\sum_{j=-\infty}^{\infty}\min\{2^{js_{0}},2^{js_{1}}2^{\nu(s_{1}-s_{0}%
)}\}\,\Vert f_{j}\Vert_{L^{\infty}(\mathbb{R}^{d})}\\
&  \hspace{-4cm}=\sum_{j=-\infty}^{\infty}2^{js_{0}}\Vert f_{j}\Vert
_{L^{\infty}(\mathbb{R}^{d})}\bigg(2^{j(s_{1}-s_{0})}\sum_{\nu=-\infty}%
^{-j}2^{\nu(s_{1}-s_{0})(1-\theta)}+\sum_{\nu=-j}^{\infty}2^{-\theta\nu
(s_{1}-s_{0})}\bigg)\\
&  \hspace{-4cm}\approx(\theta(1-\theta))^{-1}\sum_{j=-\infty}^{\infty}%
2^{js}\Vert f_{j}\Vert_{L^{\infty}(\mathbb{R}^{d})}.
\end{align*}
As a consequence, \
\begin{align*}
\Vert f\Vert_{(\dot{B}_{\infty,1}^{s_{0}}(\mathbb{R}^{d}),\dot{B}_{\infty
,1}^{s_{1}}(\mathbb{R}^{d}))_{\theta,1}^{\blacktriangleleft}}  &
\approx(\theta(1-\theta))\Vert\{\dot{\Delta}_{j}f\}_{j\in\mathbb{Z}}%
\Vert_{(\ell_{1}^{s_{0}}(L^{\infty}(\mathbb{R}^{d})),\ell_{1}^{s_{1}%
}(L^{\infty}(\mathbb{R}^{d})))_{\theta,1}}\\
& \hspace{-2cm} \approx\sum_{j=-\infty}^{\infty}2^{js}  \Vert \dot{\Delta}_j f\Vert_{L^{\infty
}(\mathbb{R}^{d})} = \left\Vert f\right\Vert _{\dot{B}_{\infty,1}^{s}(\mathbb{R}^{d})},
\end{align*}
as we wished to show.
\end{proof}

\begin{remark}\label{RemarkA}
	For $s \in \R$ and $p, q \in [1, \infty]$, the (homogeneous) Besov space $\dot{B}^s_{p, q}(\R^d)$ (cf. \cite{BCD11, BL76, Tri83}) can be defined as the space of tempered distributions modulo polynomials endowed with  the seminorm:
		\begin{equation*}
		\|f\|_{ \dot{B}^{s}_{p, q}(\mathbb{R}^{d})} := \bigg(
\sum_{j \in\mathbb{Z}} 2^{j s q} \, \|\dot{\Delta}_{j} f\|^q_{L^{p
}(\mathbb{R}^{d})} \bigg)^{1/q},
\end{equation*}
where the usual modification is made if $q=\infty$. With obvious modifications, the proof of Theorem \ref{teo:nec} can be adapted to establish  the following more general statement: If $1 \leq p, q \leq \infty, -\infty < s_0 < s_1 < \infty, \theta \in (0, 1)$, and $s = (1-\theta) s_0 + \theta s_1$, then
	\begin{equation}\label{IntBesovNormalized}
(\dot{B}_{p,q}^{s_{0}}(\mathbb{R}^{d}),\dot{B}_{p,q}^{s_{1}%
}(\mathbb{R}^{d}))_{\theta,q}^{\blacktriangleleft}=\dot{B}_{p, q}%
^{s}(\mathbb{R}^{d})
\end{equation}
with equivalence constants independent of $\theta$.  Furthermore, the inhomogeneous counterpart of \eqref{IntBesovNormalized} also holds. 
\end{remark}

\subsection{Results from Section \ref{SectionIntro1.3}}\label{sec:A1}
\subsubsection{Characterizations of $\Delta$-extrapolation spaces}

A very useful result for the computation of $\Delta$-extrapolation spaces is
the formula (\ref{form:universal}), which for convenience we reproduce here
\begin{equation}
\Delta_{\theta\in(0,1)}\bigg\{\frac{(A_{0},A_{1})_{\theta,p(\theta
)}^{\blacktriangleleft}}{\Theta(\frac{1}{1-\theta})}\bigg\}=\Delta_{\theta
\in(0,1)}\bigg\{\frac{(A_{0},A_{1})_{\theta,\infty}^{\blacktriangleleft}%
}{\Theta(\frac{1}{1-\theta})}\bigg\}. \label{app1}%
\end{equation}
The result is used implicitly in \cite{JM91} and a detailed proof can be found in
\cite[Theorem 21, page 44]{MLNM}. Further results and generalizations can be found in \cite{KM}. The import of (\ref{app1}) relies on the fact that the norm of the space on the right-hand side is a double
supremum, which allows us to apply \textquotedblleft Fubini\textquotedblright\ as follows
\begin{align}
\left\Vert f\right\Vert _{\Delta_{\theta\in(0,1)}\Big\{\frac{(A_{0}%
,A_{1})_{\theta,\infty}^{\blacktriangleleft}}{\Theta(\frac{1}{1-\theta}%
)}\Big\}}  &  =\sup_{\theta \in (0, 1)}\sup_{t>0}\frac{t^{-\theta} K(t,f;A_{0},A_{1})%
}{ \Theta(\frac{1}{1-\theta})}\nonumber\\
&  =\sup_{t>0}\frac{K(t,f;A_{0},A_{1})}{t}\sup_{\theta \in (0, 1)}\frac{t^{1-\theta}%
}{\Theta(\frac{1}{1-\theta})}\nonumber\\
&  =\sup_{t>0}\frac{K(t,f;A_{0},A_{1})}{t\inf_{\theta \in (0, 1)}\{\Theta(\frac
{1}{1-\theta})t^{-(1-\theta)}\}}\nonumber\\
&  =\sup_{t>0}\frac{K(t^{1/p_0},f;A_{0},A_{1})}{t^{1/p_0} \varphi_{\Theta}(\frac{1}{t})},
\label{la cara}%
\end{align}
where
\[
\varphi_{\Theta}(t)=\inf_{\theta \in (0, 1)} 	\bigg\{\Theta \bigg(\frac{1}{1-\theta} \bigg) \, t^{\frac{1-\theta}{p_0}} \bigg\}, \qquad p_0 > 0. 
\]
The $\varphi_{\Theta}$ functions thus played a major role in \cite{JM91}.
Remarkably, as we pointed out in Section \ref{SectionIntro1.3}, although coming from completely different considerations they essentially coincide with the original
Yudovich functions $y_{\Theta}$ (cf. \eqref{intro:Auxy})
\begin{equation}\label{4.5novel}
\varphi_{\Theta}(t) \approx y_{\Theta}(t);
\end{equation}
see also \eqref{laphidy}. 

In this paper we make a strong use of both \eqref{app1} and (\ref{la cara}). In particular, to give an explicit
characterization of the Yudovich spaces $Y^{\Theta}$ (cf. (\ref{IntroYExtra})
and (\ref{explica2})) and the sharp Yudovich spaces $Y^{\#\Theta}$ (cf. Theorem \ref{Thm2.2}, specially \eqref{212newnovel} and \eqref{212newnovel2}), to establish various characterizations of the Vishik spaces $B^\beta_\Pi$ (cf. Theorem \ref{ThmVE}, in particular, \eqref{36novel}), and to get estimates for the modulus of smoothness (cf. \eqref{alga}).

\subsubsection{Reiteration}

 The scale $\{(A_{0},A_{1})_{\theta,p(\theta
)}\}_{\theta \in (0, 1)}$ is an example of \emph{interpolation scale of exact order $\theta$}. We refer to \cite[p. 7]{JM91} for the precise definition.  Another important example is given by $\{(A_0, A_1)^{\blacktriangleleft}_{\theta, p}\}_{\theta \in (0, 1)}$ for a fixed $p \in [1, \infty]$ (cf. \eqref{IntroNormInt}). It turns out that the extrapolation formula \eqref{app1} is only a special case of a more general phenomenon based on interpolation scales of exact order $\theta$. In particular, the formula is still valid when $\{(A_{0},A_{1})_{\theta,p(\theta
)}\}_{\theta \in (0, 1)}$ is replaced by $\{(A_0, A_1)^{\blacktriangleleft}_{\theta, p}\}_{\theta \in (0, 1)}$, namely, 
\begin{equation}
\Delta_{\theta\in(0,1)}\bigg\{\frac{(A_{0},A_{1})_{\theta,p}^{\blacktriangleleft}}{\Theta(\frac{1}{1-\theta})}\bigg\}=\Delta_{\theta
\in(0,1)}\bigg\{\frac{(A_{0},A_{1})_{\theta,\infty}^{\blacktriangleleft}%
}{\Theta(\frac{1}{1-\theta})}\bigg\}. \label{app1novel}%
\end{equation}
See \cite[Theorem 21, page 44]{MLNM}.

 \emph{The ordered case $A_{1}\hookrightarrow A_{0}$}. In this case we can
achieve further simplifications for the computation of the $\Delta
$-extrapolation spaces. In particular, we claim that
\begin{equation}
\Delta_{\theta\in(0,1)}\bigg\{\frac{(A_{0},A_{1})_{\theta,p(\theta
)}^{\blacktriangleleft}}{\Theta(\frac{1}{1-\theta})}\bigg\}=\Delta_{\theta
\in(\theta_{0},1)}\bigg\{\frac{(A_{0},A_{1})_{\theta,p(\theta)}%
^{\blacktriangleleft}}{\Theta(\frac{1}{1-\theta})}\bigg\} \label{app2}%
\end{equation}
for every $\theta_0 \in (0, 1)$. This result may be considered as the extrapolation version of  the fact that the definitions of  $Y^{\Theta}_{p_0}(\Omega)$ and $Y^{\#\Theta}_{p_0}(\Omega)$ with $|\Omega| < \infty$ are independent of $p_0 \in [1, \infty)$. Indeed, recall that  $Y^{\Theta}_{p_0}(\Omega)$ and $Y^{\#\Theta}_{p_0}(\Omega)$ are the $\Delta$-extrapolation spaces relative to the pairs $\left(  L^{p_{0}}(\Omega),L^{\infty}(\Omega)\right)$ and $(L^{p_0}(\Omega), \text{BMO}(\Omega))$, respectively; cf. \eqref{IntroLp}-\eqref{IntroYExtra} and \eqref{ClaimEx1}-\eqref{ClaimEx2}. 

\begin{proof}[Proof of \eqref{app2}]

The non trivial part of the statement is the embedding $\hookleftarrow$. This can be derived as follows. We have %
\[
\left\Vert f\right\Vert _{\Delta_{\theta\in(0,1)}\Big\{\frac{(A_{0}%
,A_{1})_{\theta,p(\theta)}^{\blacktriangleleft}}{\Theta(\frac{1}{1-\theta}%
)}\Big\}}\leq\left\Vert f\right\Vert _{\Delta_{\theta\in(0,\theta_{0}%
)}\Big\{\frac{(A_{0},A_{1})_{\theta,p(\theta)}^{\blacktriangleleft}}%
{\Theta(\frac{1}{1-\theta})}\Big\}}+\left\Vert f\right\Vert _{\Delta
_{\theta\in(\theta_{0},1)}\Big\{\frac{(A_{0},A_{1})_{\theta,p(\theta
)}^{\blacktriangleleft}}{\Theta(\frac{1}{1-\theta})}\Big\}}.%
\]
To estimate the first term on the right-hand side of the previous inequality, we will make use of the fact that
\begin{equation}\label{48new}
(A_{0}%
,A_{1})_{\theta_0,p(\theta_0)} \hookrightarrow	(A_{0}
,A_{1})_{\theta,p(\theta)}^{\blacktriangleleft} \qquad \text{if} \qquad \theta < \theta_0,
\end{equation}
where the embedding constant is independent of $\theta$. This result can be shown by using similar techniques as in \cite[Lemma 4.5(ii)]{DM22}. To make the presentation self-contained, we next provide full details. Using that $A_1 \hookrightarrow A_0$, it is plain to see that
$$
	K(t, f; A_0, A_1) \approx \|f\|_{A_0} \approx K(1, f; A_0, A_1), \qquad \text{for} \qquad t > 1. 
$$
Then, by H\"older's inequality (noting that $p(\theta) = \frac{1}{1-\theta} < \frac{1}{1-\theta_0} = p (\theta_0)$) and monotonicity properties of $K$-functionals, 
\begin{align*}
	\|f\|_{(A_{0}
,A_{1})_{\theta,p(\theta)}}^{p(\theta)} & = \int_0^1 [t^{-\theta} K(t, f; A_0, A_1)]^{p(\theta)} \, \frac{dt}{t} +  \int_1^\infty [t^{-\theta} K(t, f; A_0, A_1)]^{p(\theta)} \, \frac{dt}{t} \\
& \approx \int_0^1 [t^{-\theta} K(t, f; A_0, A_1)]^{p(\theta)} \, \frac{dt}{t} + \frac{1}{\theta p(\theta)} \,  \|f\|_{A_0}^{p(\theta)} \\
& \leq \bigg( \int_0^1 [t^{-\theta_0} K(t, f; A_0, A_1)]^{p(\theta_0)} \, \frac{dt}{t} \bigg)^{p(\theta)/p(\theta_0)} + \frac{1}{\theta p(\theta)} \,  \|f\|_{A_0}^{p(\theta)} \\
& \lesssim \frac{1}{\theta} \,  \bigg( \int_0^1 [t^{-\theta_0} K(t, f; A_0, A_1)]^{p(\theta_0)} \, \frac{dt}{t} \bigg)^{p(\theta)/p(\theta_0)}  \\
& \leq \frac{1}{\theta} \, \|f\|_{(A_0, A_1)_{\theta_0, p(\theta_0)}}^{p(\theta)},
\end{align*}
where we have also used that $\theta \in (0, \theta_0)$ (and so $p(\theta) \approx 1$) in the penultimate estimate. The proof of \eqref{48new} is finished. 

It follows from \eqref{48new} that
\begin{align*}
\sup_{\theta\in(0,\theta_{0})}   \frac{\left\Vert f\right\Vert _{(A_{0}%
,A_{1})_{\theta,p(\theta)}^{\blacktriangleleft}}} {\Theta (\frac{1}{1-\theta})}
&  \lesssim \left\Vert f\right\Vert _{(A_{0},A_{1})_{\theta_{0},p(\theta_{0}%
)}}\sup_{\theta\in(0,\theta_{0})}\frac{1}{\Theta(\frac
{1}{1-\theta})}\\
&  = \frac{\left\Vert f\right\Vert
_{(A_{0},A_{1})_{\theta_{0},p(\theta_{0})}}}{\Theta(1)} \\
&  \lesssim \sup_{\theta\in(\theta_{0}, 1)}   \frac{\left\Vert f\right\Vert _{(A_{0}%
,A_{1})_{\theta,p(\theta)}^{\blacktriangleleft}}} {\Theta (\frac{1}{1-\theta})}.
\end{align*}
This completes the proof of \eqref{app2}. 

\end{proof}

\subsection{Computability of $K$-functionals}\label{SectionA3}

A central issue in interpolation theory is to find explicit expressions for $K$-functionals. Next we list some well-known examples of characterizations for $K$-functionals (see also Example \ref{ejemplo:markao}) and we refer the interested reader to \cite{BS88, BL76} for further examples.

\begin{example}\label{Ex1.a} Let $(L^{p_0}(\Omega),L^{\infty}(\Omega)), \, p_0 \in (0, \infty)$. Then (cf. \cite[Theorem 5.2.1, page 109]{BL76} and 
\cite[Theorem 1.6, page 298]{BS88})
\begin{equation}\label{A9new}
K(t,f;L^{p_0}(\Omega),L^{\infty}(\Omega)) \approx \bigg(\int_{0}^{t^{p_0}}(f^{\ast}(\xi))^{p_0} \, d \xi \bigg)^{1/p_0}. 
\end{equation}
Equality holds in \eqref{A9new} if $p_0 = 1$. 
\end{example}

\begin{example}\label{Ex2.a}Let $(L^{p_0}(\R^d), \text{BMO}(\R^d)), \, p_0 \in (0, \infty)$. Then (cf. \cite[Corollary 3.3]{JT85})
\begin{equation}\label{A10novel}
K(t, f; L^{p_{0}}(\R^d), \text{BMO}(\R^d)) \approx\bigg(\int_{0}^{t^{p_{0}}}
[(M^{\#}_{\R^d} f)^{*}(\xi)]^{p_{0}} \, d \xi\bigg)^{1/p_{0}},
\end{equation}
where $M^{\#}_{\R^d}$ is the maximal function given in \eqref{IntroDefMaxS}. The local counterpart for function spaces defined on cubes also holds true. In order to suitably interpret
$(L^{p_{0}}(\R^d), \text{BMO}(\R^d))$ as an interpolation pair, it is
necessary to factor out constant functions. Then $L^{p_{0}}(\R^d)$ can be
identified with $L^{p_{0}}(\R^d)/\mathbb{C}$, where $\|f\|_{L^{p_{0}}%
(\R^d) / \mathbb{C}} = \inf_{c \in\mathbb{C}} \|f-c\|_{L^{p_{0}}(\R^d)}$.
\end{example}

 \begin{example}\label{Ex3.a}The $K$-functional for the pair $(L^{\infty}(\R^d),\dot{W}^{1}_\infty(\R^d))$ plays an important role in our work  (cf. (\ref{IntroKFunctModLinfty})) and  can be characterized as  (cf. \cite[(4.42), p. 341]{BS88} and \cite[Theorem 1]{JS77})
\begin{equation}\label{A11novel}
K(t,f;L^{\infty}(\R^d),\dot{W}_{\infty}^{1} (\R^d))\approx\sup_{\left\vert x-y\right\vert
\leq t}\left\vert f(x)-f(y)\right\vert. 
\end{equation}
The corresponding formula for periodic functions  is also true. 
\end{example}

A natural assumption for $(A_{0},A_{1})$ is to be \emph{Gagliardo closed} in the sense  that
\begin{equation*}
\left\Vert f\right\Vert _{A_{0}} \approx \sup_{t > 0} K(t,f;A_{0}%
,A_{1}) \qquad \text{and} \qquad \left\Vert f\right\Vert _{A_{1}}\approx \sup_{t > 0}%
\frac{K(t,f;A_{0},A_{1})}{t}. 
\end{equation*}
See \cite[p. 320]{BS88}. 
This condition is easily verified for many classical pairs of spaces, in particular, for the pairs given in Examples \ref{Ex1.a}--\ref{Ex3.a}.

\begin{proposition}
\label{ejemplomarkao}Suppose that $(A_{0},A_{1})$ is Gagliardo closed, and let
$\Theta(p)\approx 1$. Then%
\[
\Delta_{\theta\in(0,1)}\bigg\{\frac{(A_{0},A_{1})_{\theta,p(\theta
)}^{\blacktriangleleft}}{\Theta(\frac{1}{1-\theta})}\bigg\}=A_{0}\cap A_{1}.
\]

\end{proposition}

\begin{proof}
Using \eqref{app1} and the monotonicity properties of $K$-functionals (note that $K(t,f;A_{0},A_{1})$ increases and $\frac{K(t,f;A_{0},A_{1})}{t}$ decreases), we have
\begin{align*}
\|f\|_{\Delta_{\theta\in(0,1)}\Big\{\frac{(A_{0},A_{1})_{\theta,p(\theta
)}^{\blacktriangleleft}}{\Theta(\frac{1}{1-\theta})}\Big\}} & \approx \sup_{\theta \in (0, 1)}\sup_{t>0} t^{-\theta} K(t,f;A_{0},A_{1}) \\
 &  \approx\sup_{\theta \in (0, 1)
}\sup_{t>1}t^{-\theta} K(t,f;A_{0},A_{1})+\sup_{\theta \in (0, 1)}\sup_{t<1} t^{-\theta}
K(t,f;A_{0},A_{1})\\
&  =\sup_{t>1}K(t,f;A_{0},A_{1})\sup_{\theta \in (0, 1)} t^{-\theta}+\sup_{t<1}%
\frac{K(t,f;A_{0},A_{1})}{t}\sup_{\theta \in (0, 1)}t^{1-\theta}\\
&  =\sup_{t>1}K(t,f;A_{0},A_{1})+\sup_{t<1}\frac{K(t,f;A_{0},A_{1})}{t}\\
&  \approx \left\Vert f\right\Vert _{A_{0}}+\left\Vert f\right\Vert _{A_{1}}.
\end{align*}

\end{proof}

The previous result can be used to readily justify the assertions made in
(\ref{introYLinfty}) and (\ref{introYLinftysharp}) for $Y^\Theta_{p_0}(\Omega)$ and $Y^{\#\Theta}_{p_0}(\Omega)$ via the extrapolation formulae \eqref{IntroYExtra} and \eqref{ClaimEx2}. For Vishik spaces, we have $\dot{B}_{\Pi}^{\beta}(\mathbb{R}^{d})=\dot
{B}_{\infty,1}^{\beta}(\mathbb{R}^{d}),$ $\Pi(p)\approx1$, and one could apply
Proposition \ref{ejemplomarkao} to obtain the corresponding result that appears in \eqref{VFM2}.

%
%
%

\end{document}